\documentclass[12pt]{article}
\usepackage[utf8]{inputenc}
\usepackage{amsmath}
\usepackage{amsthm}
\usepackage{amsfonts}
\usepackage{amssymb}
\usepackage{mathabx}
\usepackage{mathtools}
\usepackage{thm-restate}
\usepackage{enumerate}
\usepackage[normalem]{ulem}

\usepackage{tikz}
\usepackage[cmtip,arrow]{xy}
\usepackage{pb-diagram, pb-xy}
\usetikzlibrary{arrows}
\usepackage{tikz-cd}
\usepackage[margin=1in]{geometry}
\usepackage{graphicx}
\usepackage{caption}
\usepackage{subcaption}
\usepackage{appendix}
\usepackage[export]{adjustbox}

\usepackage[nodayofweek]{datetime}
\usepackage{hyperref}
\usepackage{hyphenat}

\usepackage{color}
\usepackage{transparent}

\definecolor{Green}{RGB}{30, 150, 30}

\newtheorem{thm}{Theorem}[section]
\newtheorem{lem}[thm]{Lemma}
\newtheorem{prop}[thm]{Proposition}
\newtheorem{cor}[thm]{Corollary}

\theoremstyle{definition}
\newtheorem{defn}[thm]{Definition}
\newtheorem{claim}[thm]{Claim}
\newtheorem{example}[thm]{Example}
\newtheorem{remark}[thm]{Remark}
\newtheorem{notation}[thm]{Notation}

\newenvironment{subproof}[1][\proofname]{%
	\begin{proof}[#1]%
	}{%
	\end{proof}%
}

\newcommand{\C}{\mathcal{C}}

\newcommand{\base}{\operatorname{base}}

\newcommand{\mcg}{\operatorname{MCG}}

\newcommand{\nest}{\sqsubseteq}
\newcommand{\propnest}{\sqsubsetneq}
\newcommand{\trans}{\pitchfork}

\newcommand{\diam}{\operatorname{diam}}
\newcommand{\mc}[1]{\mathcal{#1}}
\newcommand{\mf}[1]{\mathfrak{#1}}
\newcommand{\Push}{\operatorname{Push}}

\newcommand{\isom}{\operatorname{Isom}}

\newcommand{\stab}{\operatorname{Stab}}
\newcommand{\sat}{\operatorname{Sat}}
\newcommand{\markcurve}{\operatorname{markarc}}
\newcommand{\unmark}{\operatorname{unmark}}
\newcommand{\ball}{\operatorname{ball}}

\newcommand{\st}{\operatorname{st}}
\newcommand{\lk}{\operatorname{lk}}

\newcommand{\bsat}{\operatorname{BSat}}
\newcommand{\blk}{\operatorname{lk}_{B}}
\newcommand{\xlk}{\operatorname{lk}_{X}}

\newcommand{\cl}{\operatorname{cl}}

\newcommand{\supp}{\operatorname{supp}}

\newcommand{\Sdot}{S^z}

\newcommand{\tsh}[1]{\left\{\kern-.7ex\left\{#1\right\}\kern-.7ex\right\}}

\title{Extensions of Multicurve Stabilizers are Hierarchically Hyperbolic}
\author{Jacob Russell}

\begin{document}

\maketitle

\begin{abstract}
	\noindent For a  closed and orientable surface $S$ with genus at least 2, we prove the $\pi_1(S)$-extensions of the stabilizers of  multicurves on $S$ are hierarchically hyperbolic groups. This answers a question of Durham, Dowdall, Leininger, and Sisto. We also include an appendix that employs work of Charney, Cordes, and Sisto to characterize the Morse boundaries of  hierarchically hyperbolic groups  whose largest acylindrical action on a hyperbolic space is on a quasi-tree.
\end{abstract}

\tableofcontents

\section{Introduction}

For an orientable closed surface $S$  and any group $G$, a $\pi_1(S)$-extension of $G$ is any  group $E$ that fits into  the short exact sequence \[ 1 \to \pi_1(S) \to E \to G \to 1.\] Topologically, $\pi_1(S)$-extensions arise as the  fundamental groups of surface bundles. 

When $S$ has genus at least 2, the Birman exact sequence tells us that the mapping class group of $S$ and the mapping class group fixing a marked point $z \in S$ fit into the  short exact sequence \[ 1 \to \pi_1(S) \to \mcg(S;z) \to \mcg(S) \to 1 \] where $\pi_1(S)$ is identified with the point pushing subgroup of $\mcg(S;z)$.  This allows us to create a $\pi_1(S)$-extension of  any subgroup $G < \mcg(S)$ by taking the full preimage of $G$ in $\mcg(S;z)$.

Because the mapping class group of $S$ is naturally isomorphic to the outer automorphism group of $\pi_1(S)$, every $\pi_1(S)$-extension of an abstract group $G$  produces a monodromy homomorphism  $G \to \mcg(S)$.  The extension groups arising from the Birman exact sequence are therefore the extensions whose monodromies have trivial kernel. Since a $\pi_1(S)$-extension is determined up to isomorphism by the monodromy into $\mcg(S)$,  understanding how properties of the monodromy influence the properties of the $\pi_1(S)$-extension group is an important problem in mapping class groups.

The first examples of $\pi_1(S)$-extensions arise from 3-manifolds that fiber over the circle with fiber $S$. In these examples, $G \cong Z$ and,  by Thurston's Geometerization Theorem,
any such $3$-manifold admits a hyperbolic structure if  and only if the image of the monodromy is generated by a pseudo-Anosov element of $\mcg(S)$; see \cite{Thurston_fibering_over_circle}. Work of Farb and Mosher \cite{Farb_Mosher_convex_cocompact} plus Hamenst\"adt \cite{Hamenstadt_extensions_of_surface_groups}  expanded Thurston's result to prove that a $\pi_1(S)$-extension of a group $G$ is \emph{Gromov hyperbolic} if and only if the monodromy has finite kernel and \emph{convex cocompact} image; see also \cite{Mj_Pranab_Metric_bundles}. 	

In the present work, we study the geometry of $\pi_1(S)$-extensions when $G$ is the stabilizer in $\mcg(S)$ of a multicurve on $S$. While the presence of Dehn twists in $G$ prevents these extensions from ever being Gromov hyperbolic, we prove they are \emph{hierarchically hyperbolic}. This gives a positive answer to a question of  Dowdall, Durham, Leininger, and Sisto \cite[Question 1.12]{DDLS_Veech_groups_II} asked as part of a search for a robust definition of geometric finiteness in the mapping class group; see Section \ref{sec:geom_finite} for further discussion of this motivation.

\begin{thm}\label{intro_thm:E_alpha_is_HHG}
	Let $S$ be a closed orientable surface with genus at least 2. Let $\alpha$ be a multicurve on $S$ and $G_\alpha$ be the stabilizer of $\alpha$ in $\mcg(S)$. If $E_\alpha$ is the full preimage of $G_\alpha$ in $\mcg(S;z)$, then $E_\alpha$ is a hierarchically hyperbolic group.
\end{thm}

Hierarchical hyperbolicity was introduced by Behrstock, Hagen, and Sisto to axiomatize the coarse geometric structure of the mapping class group arising from the machinery of Masur and Minsky \cite{BHS_HHSI,BHS_HHSII}. Masur and Minsky proved that the curve complex is Gromov hyperbolic  and   used  projections of $\mcg(S)$ onto the curve complexes of subsurfaces of $S$ to greatly illuminate the geometry of the mapping class group \cite{MMI,MMII}. 
 A number of subsequent results  built on this work, producing a beautify theory of how the geometry of $\mcg(S)$ can be decoded from a combination of the geometry of these curve complexes  and combinatorial information about the subsurfaces of $S$, e.g., \cite{Behrstock_Thesis,behrstock_Drutu_Mosher_Thickness,BKMM_Rigidity,BM_MCG_Rank,bowditch_uniform_hyperbolicity}. Hierarchical hyperbolicity axiomatizes this theory,  describing a class of spaces whose coarse geometry is encoded in a collection of projections onto hyperbolic metric spaces that are organized by a set of combinatorial relations.  Remarkably, the class of hierarchically hyperbolic spaces encompasses a variety of groups beyond the mapping class group including the fundamental group of most $3$-manifolds \cite{BHS_HHSII}, many cocompactly cubulated groups \cite{BHS_HHSI,Hagen_Susse_factor_system}, Artin groups of extra large type \cite{Extra_Large_Artin}, and several combinations of hyperbolic groups \cite{BR_Combination,Robbio_Spriano_Combination,BR_graph_products}. Hierarchical hyperbolicity also describes the coarse geometry of a number of other groups and spaces associated to surfaces such as Teichm\"uller space with both the Teichm\"uller and Weil--Peterson metrics \cite{BHS_HHSI,MMI,MMII, BKMM_Rigidity,Brock_WeilPetersson, Durham_Combinatorial_Teich,Rafi_Combinatorial_Teich, EMR_Teich_Rank}, the genus 2 handlebody group \cite{Chesser_genus_2_handlebody}, the $\pi_1(S)$-extensions of lattice Veech groups \cite{DDLS_Veech_groups_II}, certain quotients of the mapping class group \cite{BHS_HHS_AsDim,BHMS_cHHS}, and a wide variety of graphs built from curves on surfaces \cite{vokes}.

Hierarchical hyperbolicity produces a large number of geometric and algebraic consequences, e.g., \cite{BHS_HHSI,BHS_HHSII,BHS_HHS_Quasiflats,HHP_Coarse_Helly_and_HHG,RST_Quasiconvexity}. The following corollary states some salient examples of the new results that are gained automatically as a result of the hierarchical hyperbolicity of $E_\alpha$.

\begin{cor}
	Let $S$ be a closed orientable surface with genus at least 2. Let $\alpha$ be a multicurve on $S$ and let $E_\alpha$ be the full preimage in $\mcg(S;z)$ of the stabilizer of $\alpha$ in $\mcg(S)$.
	\begin{enumerate}
		\item $E_\alpha$ has quadratic Dehn function.
		\item $E_\alpha$ has finitely many conjugacy classes of finite order subgroups.
		\item $E_\alpha$ is semi-hyperbolic and hence has solvable conjugacy problem.
	\end{enumerate}
\end{cor}

Our proof of Theorem \ref{intro_thm:E_alpha_is_HHG} uses the recent \emph{combinatorial hierarchical hyperbolicity} machinery of Behrstock, Hagen, Martin, and Sisto \cite{BHMS_cHHS}. The key to this approach is to construct a hyperbolic simplicial complex $X_\alpha$ for $E_\alpha$ that is analogous to the curve complex for $\mcg(S)$. To define $X_\alpha$, let $\Sdot$ denote the surface obtained from $S$ by adding $z \in S$ as a marked point, and let $\Pi \colon \Sdot \to S$ be the map given by forgetting that $z$ is a marked point.  The vertices of $X_\alpha$  are all isotopy classes of curves $c$ on $\Sdot$ so that $\Pi(c)$ is  either contained in or disjoint from $\alpha$. There are edges between two vertices of $X_\alpha$ if the isotopy classes have disjoint representatives.  

We prove that $X_\alpha$ is not only hyperbolic, but  in fact a quasi-tree. Combining this with results on hierarchical hyperbolicity  from the literature produces several additional properties of $E_\alpha$.

\begin{thm}\label{intro_thm:quasi-tree}
	The graph $X_\alpha$ is uniformly quasi-isometric to a tree and the group $E_\alpha$ has the following properties.
	\begin{itemize}
		\item $E_\alpha$ acts acylindrically on $X_\alpha$ and this action is largest in the sense of \emph{\cite{ABO_largest_action}}.
		\item A subgroup $H <E_\alpha$ is stable if and only if the orbit map of $H$ in $X_\alpha$ is a quasi-isometric embedding. In particular, every stable subgroup of $E_\alpha$ is virtually free.
		\item The Morse boundary of $E_\alpha$ is an $\omega$-Cantor space.
	\end{itemize}
\end{thm}

The proof that the Morse boundary of $E_\alpha$ is an $\omega$-Cantor space uses  a technique developed by Charney, Cordes, and Sisto \cite{CCS_omega_Morse_boundary}. In the appendix, we work out this application not just for $E_\alpha$, but for all  hierarchically hyperbolic groups whose largest acylindrical action on a hyperbolic space is on a quasi-tree.

\subsection{Geometric finiteness in $\mcg(S)$}\label{sec:geom_finite}
There is a longstanding and fruitful analogy between discrete subgroups of $\isom(\mathbb{H}^n)$ and the subgroups of $\mcg(S)$; see \cite{Kent_Leininger_Survey} for an detailed overview. In the case of $\isom(\mathbb{H}^n)$, the best behaved discrete subgroups are the geometrically finite subgroups and the convex cocompact subgroups. A discrete subgroup $\Gamma <\isom(\mathbb{H}^n)$ is \emph{geometrically finite} if $\Gamma$ acts with finite co-volume on the convex hull of its limit set, and it is \emph{convex cocompact} if it instead acts cocompactly.

Building on this analogy, Farb and Mosher defined a \emph{convex cocompact} subgroup of $\mcg(S)$ as a subgroup whose orbit in the Teichm\"uller space of $S$ is quasiconvex \cite{Farb_Mosher_convex_cocompact}. Convex cocompact subgroups have since seen immense study, yielding a variety of  characterizations \cite{Kent_Leininger_convex_cocompactness,Kent_Leininger_convergence,DT_stability,BBKL_undistorted_pA}. Notable to the study of $\pi_1(S)$-extensions, a subgroup of $\mcg(S)$ is convex cocompact if and only if the full preimage in $\mcg(S;z)$ is Gromov hyperbolic \cite{Farb_Mosher_convex_cocompact,Hamenstadt_extensions_of_surface_groups,Mj_Pranab_Metric_bundles}.

While convex cocompactness has been well studied in the mapping class group, the lack of negative curvature in the Teichm\"uller space of $S$ has prevented the formulation of a robust analogue of  geometric finiteness for $\mcg(S)$. Despite this, there are several  subgroups of $\mcg(S)$ that one could naturally consider as geometrically finite. The most notable examples are the Veech groups, all of $\mcg(S)$ itself, and the stabilizers of multicurves on $S$. These subgroup  are all not convex cocompact, but each acts with finite co-volume on a well-behaved subset of the Teichm\"uller space of $S$. 

Dowdall, Durham, Leininger, and Sisto proposed that  geometric finiteness (however it is defined) in $\mcg(S)$ should be characterized by some sort of hyperbolicity of the extension group that encompasses Gromov hyperbolicity in the convex cocompact case \cite{DDLS_Veech_groups_II}. In support of this idea,  they proved that the full preimage in $\mcg(S;z)$ of a lattice Veech group is hierarchically hyperbolic and asked if the same could be proven for other candidates of geometrically finite subgroups \cite[Question 1.12]{DDLS_Veech_groups_II}. Theorem \ref{intro_thm:E_alpha_is_HHG} affirmatively answers this question for the stabilizers of multicurves. Combining these results with the observation that the full preimage of the entire $\mcg(S)$ is $\mcg(S;z)$---which is hierarchically hyperbolic---we now know that all three of the most naturally geometrically finite subgroups have hierarchically hyperbolic extension groups.

\subsection{Outline of the paper}

Section \ref{sec:background} reviews background material and notation on simplicial complexes, coarse geometry, group actions on graphs, curves on surfaces, and combinatorial hierarchically hyperbolic spaces. In Section \ref{sec:sketch}, we provide a sketch of the hierarchical hyperbolicity of $E_\alpha$ when $\alpha$ is a single curve.
Section \ref{sec:tree guided} presents the criteria we will use to prove the graph  $X_\alpha$ is a quasi-tree. 

The proof that $E_\alpha$ is a hierarchically hyperbolic group is spread over two sections.  In Section \ref{sec:non-annular_case}, we construct the quasi-tree $X_\alpha$ and use it to build a combinatorial HHS that $E_\alpha$ acts on non-properly. The source of this non-properness is the fact that Dehn twists in $E_\alpha$ will stabilize simplices in $X_\alpha$. In Section \ref{sec:annuli}, we address this lack of properness by constructing a new combinatorial HHS using a ``blow-up'' of the graph $X_\alpha$. This blow-up of $X_\alpha$ will prevent Dehn twists from fixing simplices by recording the action of twists around curves in $X_\alpha$. This in-turn causes $E_\alpha$ to act metrically properly on the resulting  combinatorial HHS. 

Section \ref{sec:structure} is largely expository, elaborating on the hierarchically hyperbolic structure of $E_\alpha$ and some of its consequences. We start by giving a description  of the hierarchically hyperbolic group structure on $E_\alpha$ using the topology of curves and subsurfaces. We then use this description to show that  a minor modification of the HHG structure produces an HHG structure for $E_\alpha$ that has an additional property introduced by Abbott, Behrstock, and Durham called \emph{unbounded products}. This allows us to  use results from the literature plus the fact that $X_\alpha$ is a quasi-tree to achieve the conclusions of Theorem \ref{intro_thm:quasi-tree}.

We include an appendix that shows how a technique developed by Charney, Cordes, and Sisto to understand the Morse boundary of right-angled Artin groups and graph manifolds groups can be extended to certain hierarchically hyperbolic groups. This extension is required to conclude the Morse boundary of $E_\alpha$ is an $\omega$-Cantor space.

\subsection*{Acknowledgments} The author is very grateful to Chris Leininger for abundant conversations on the topic and for detailed comments on drafts of this paper. The author also thanks Alessandro Sisto and Mark Hagen for useful conversations about combinatorial HHSs. The author is grateful for the work of an anonymous referee whose many comments improved the exposition of this paper.

\section{Background and notation}\label{sec:background}

\subsection{Simplicial Complexes}

Throughout, if $X$ is a simplicial complex or graph, then $X^0$ will denote the set of vertices of $X$. If the notation for the simplicial complex or graph depends on some parameter, for example $\C(S)$, then we will instead insert the superscript $0$ between the identifying symbol and the parameter, for example, $\C^0(S)$ is the vertices of the complex $\C(S)$. We equip each graph with the path metric coming from declaring each edge to have length 1. When considering a simplicial complex as a metric space we will implicitly be referring to the metric space that is the 1-skeleton of the simplicial complex equipped with this path metric.

We will make frequent use of the following notions of join, link, and star.

\begin{defn}[Join, link, and star]
	Let $X$ be a flag simplicial complex. If $Y$ and $Z$ are disjoint flag subcomplexes of $X$ so that every vertex of $Y$ is joined by an edge to every vertex of  $Z$, then the \emph{join} of $Y$ and $Z$, $Y \bowtie Z$, is the subcomplex of $X$ spanned by $Y$ and $Z$. Given a subcomplex $Y$ of $X$, the \emph{link} of $Y$, $\lk(Y)$, is the subcomplex of $X$ spanned by the vertices of $X$ that are joined by an edge to all the vertices of $Y$. The \emph{star} of $Y$, $\st(Y)$, is the join $Y \bowtie \lk(Y)$. 
\end{defn}

We say a simplex $\Delta$ of a simplicial complex $X$ is \emph{maximal} if $\lk(\Delta) = \emptyset$. This is equivalent to saying $\Delta \subseteq \Delta'$ implies $\Delta = \Delta'$  for any simplex $\Delta'$. Note that the maximal simplices of $X$ need not all have the same number of vertices. Finally, if $Y \subseteq X^0$, then we let $X-Y$ denote the subcomplex of $X$ spanned by the vertices of $X$ that are not in $Y$.

\subsection{Coarse geometry and hyperbolic spaces}
Throughout, when we measure the distances between sets in a metric space (or between a set a point), we are taking the minimum distance between the two sets.

A map $f \colon X \to Y$ between two metric spaces is a \emph{$(\lambda,\epsilon)$-quasi-isometric embedding} if $\lambda \geq 1$ and $\epsilon \geq 0$ and  for all $x_1,x_2 \in X$ we have \[ \frac{1}{\lambda} d_X(x_1,x_2) -\epsilon \leq d_Y(f(x_1),f(x_2)) \leq \lambda d_X(x_1,x_2) + \epsilon.\]  A $(\lambda,\epsilon)$-quasi-isometric embedding $f\colon X \to Y$ is a \emph{$(\lambda,\epsilon)$-quasi-isometry} if  for each $y \in Y$, there exists $x \in X$ with $d_Y(y,f(x)) \leq \epsilon$. Because we will largely be interested in geometry up to quasi-isometry, we will  often work with coarse maps instead of genuine functions. A \emph{coarse map} between metric spaces $X$ and $Y$ is a function $f\colon X \to 2^Y$ so that $f(x)$ is a non-empty subset of $Y$ with uniformly bounded diameter for all $x \in X$. By abuse of notation, we will often use the notation $f \colon X \to Y$ to signify coarse maps. We say a coarse map is a $(\lambda,\epsilon)$-quasi-isometric embedding (resp. a quasi-isometry) if it satisfies the same inequalities given above (with distances between sets being the minimum distance). A coarse map $f \colon X \to Y$ is \emph{$(\lambda,\epsilon)$-coarsely Lipschitz} if \[d_Y(f(x_1),f(x_2)) \leq \lambda d_X(x_1,x_2) + \epsilon\] for each $x_1,x_2 \in X$. 

When $X$ is a graph, a coarse map defined on the vertices $X^0$ induces a coarse map on the entire graph by sending points on edges to the  union of  the images of the vertices of that edge. In this case, we freely use the following lemma to ensure that the induced map is a coarsely Lipschitz coarse map.

\begin{lem}\label{lem:checking_coarsely_lipschitz}
	Let $X$ be a connected graph and $f \colon X^0 \to Y$ be a coarse map.  If there exist $\lambda\geq 1$ so that $\diam_Y(f(x)) \leq \lambda$ for each $x$ and  $d_Y(f(x_1),f(x_2)) \leq \lambda$ for every $x_1,x_2 \in X^0$ that are joined by an edge, then the map $f \colon X \to Y$ induced by the map on the vertices is $(3\lambda,0)$-coarsely Lipschitz coarse map.
\end{lem}

\begin{proof}
	Extend $f$ to points on edges by defining the image of an edge point to be the union of the images of the vertices of that edge. This is a coarse map by the hypotheses. Let $x,x' \in X$ and $x_0,x_1,\dots,x_n$ be the vertices on a geodesic in $X$ connecting $x$ to $x'$. Since we are measuring distance between sets as the minimum distance, we have
	\begin{align*}
	d_Y(f(x),f(x')) \leq d_Y(f(x_0),f(x_n)) \leq& \sum_{i=1}^n \diam_Y(f(x_{i-1}) \cup f(x_i)) \\
\leq& 3\lambda n = 3\lambda d_X(x_0,x_n) \leq 3\lambda d_X(x,x'). \qedhere
	\end{align*}  
\end{proof}

\noindent We will also frequently use the following criteria to prove that subsets of a graph are quasi-isometrically embedded.

\begin{lem}\label{lem:coarsel_lipschitz}
	Let $H$ and $Y$ be connected graphs equipped with the path metric. Suppose that $H$ is a subgraph of $Y$. If there exists a coarse map $\psi \colon Y \to H$ so that
	\begin{itemize}
		\item $\psi$ is the identity on $H^0$;
		\item  $\psi$ is $(\lambda,\epsilon)$-coarsely Lipschitz;
	\end{itemize}
	then, the inclusion of $H$ into $Y$ is a $(\lambda,\epsilon+2\lambda)$-quasi-isometric embedding.
\end{lem}

\begin{proof}
	 Let $h_1,h_2 \in H^0$. Since $H$ is a connected subgraph, we have $d_Y(h_1,h_2) \leq d_H(h_1,h_2)$. Since $\psi(h_1) = h_1$ and $\psi(h_2) = h_2$ we have $d_H(h_1,h_2 ) \leq \lambda d_Y(h_1,h_2) + \epsilon$, hence the inclusion is a $(\lambda,\epsilon)$-quasi-isometry on the vertices. Since each point on an edge is at most 1 from two vertices, this implies the inclusion is a $(\lambda,\epsilon+2\lambda)$-quasi-isometric embedding on all of $H$.
\end{proof}

We say that a geodesic metric space is \emph{$\delta$-hyperbolic} if for every geodesic triangle in the space, the $\delta$-neighborhood of any two sides covers the third side. If $X$ and $Y$ are quasi-isometric geodesic metric spaces and $X$ is $\delta$-hyperbolic, then $Y$ is $\delta'$-hyperbolic for some $\delta'$ depending on $\delta$ and the constants  of the quasi-isometry $X \to Y$. A connected graph $X$ is $0$-hyperbolic  if and only if $X$ is a tree. We say that a geodesic metric space is a \emph{quasi-tree} if it is quasi-isometric to a tree.

If $X$ is a $\delta$-hyperbolic graph and $H$ is a connected subgraph so that the inclusion $H \to X$ is a $(\lambda,\epsilon)$-quasi-isometric embedding, then for each $x \in X$ the set \[\{h \in H : d_X(x,h) = d_X(x,H)\}\] is uniformly bounded in terms of $\delta$, $\lambda$, and $\epsilon$.  We therefore have a coarse map $\mf{p} \colon X \to H$ given by  $\mf{p}(x) = \{h \in H : d_X(x,h) = d_X(x,H)\}$ (note, $\mf{p}$ is the identity on $H$). We call this map the \emph{coarse closest point projection} onto $H$. The coarse closest point projection is always coarsely Lipschitz \cite[Proposition III.$\Gamma$.3.11]{BH}.

\subsection{Groups acting on graphs}
Let $G$ be a finitely generated group acting on a connected simplicial graph $X$ by graph automorphisms. Let $\ball_r(x)$ denote the closed ball of radius $r \geq 0$ in $X$ around the vertex $x \in X^0$. The action of $G$ on $X$ is \emph{metrically proper} if for every $r \geq 0$ and $x \in X^0$, the set $$\{g \in G : g \cdot\ball_r(x) \cap \ball_r(x) \neq \emptyset\}$$ is finite. The action of $G$ on $X$ is \emph{cocompact} if the quotient of $X$ by the action of $G$ is compact. Similarly, the action is \emph{cobounded} if the quotient has finite diameter.

\subsection{Surfaces, curves, and mapping class groups} \label{sec:surfaces}
Let $S_{g,p}^n$ denote a connected orientable surface  with genus $g$, $p$ punctures, and $n$ marked points. The complexity of $S_{g,p}^n$ is $\xi(S_{g,p}^n) = 3g-3+p+n$.  Given a non-marked point $z$ on the surface $S \cong S_{g,p}^n$, let $\Sdot$ denote the surface obtained from $S$ by adding $z$ as an additional mark point. 

By a \emph{curve on $S$}, we mean an isotopy class of a simple closed curve on $S$ that is essential on $S$ minus the marked points of $S$.  Here, two curves are considered isotopic if they are isotopic on $S$ minus the marked points. We say  two curves are \emph{disjoint} if they are not equal and have disjoint representatives. A \emph{multicurve} on $S$  is a set of pairwise disjoint  curves on $S$. For two curves $c$ and $c'$,  we   define the \emph{intersection number},  $i(c,c')$, to be the minimal number of intersection points between representatives of $c$ and $c'$.  Note, $i(c,c')=0$ means either $c = c'$ or $c$ and $c'$ are disjoint.

By a \emph{subsurface} of $S$, we mean an isotopy class of an open essential subsurface of $S$, where  two subsurfaces are considered isotopic if they are isotopic on $S$ minus the marked points. Subsurfaces of $S$ are not required to be connected. For a subsurface $U \subseteq S$, we use $\partial U$ to denote the isotopy class of the curves on $S$ that are the boundary curves for any representative of $U$.  If $A$ is an annular subsurface of $S$, then the \emph{core curve} of $A$ is the isotopy class of $\partial A$.   

The set of all subsurface of $S$ has a partial order  denoted $\subseteq$ and a difference  denoted $-$. For  connected non-annular subsurfaces $U$ and $V$, we write  $U \subseteq V$ if  $U$ has a representative that is contained in a representative of $V$.  In this case, we  can subtract $U$ from $V$ as follows: fix a complete hyperbolic metric on $S$ minus the marked points and pick representatives for $V$ and $U$ whose boundary curves are geodesic. The subsurface $V - U$ is the open subsurface that is the interior of the difference of these representatives of $V$  and $U$.
For an annular subsurface $A$ and a connected subsurface $V$, we write $A \subseteq V$ if either  $V =A$ or $V$ is not an annulus and there is a representative of $A$ that is a non-peripheral annulus on a representative of $V$.  If $V = A$, then $V-A =\emptyset$. If $A \subsetneq V$, then there is a representative of $\partial A$ that is contained in a representative of $V$. Define $V-A$ to be the isotopy class of the subsurface obtained by subtracting this representative of $\partial A$ from this representative of $V$.  If $U$ and $V$ are possibly disconnected surfaces, then $U \subseteq V$ if every component of $U$ is contained in some component of $V$, where containment is defined in the above case of connected subsurfaces. Similarly, $V-U$ is defined for disconnected surfaces  with $U \subseteq V$ by extending the difference defined above over connected components.

There is also a containment relation between multicurves and subsurfaces as well as a difference operation whenever a subsurface contains a multicurve. If $U$ is a connected non-annular subsurface and a curve $c$ has a representative that is a non-peripheral curve on some representative of $U$, then we say $c$ is \emph{contained} in $U$ and write $c \subseteq U$.  If $A$ is an annulus, then we declare that the only curve that is contained in $A$ is the core curve of $A$ and we write $c \subseteq A$ only if $c$ is the core curve of $A$. If every component of a multicurve $\mu$ is contained in a (possibly disconnected) subsurface $U$, then we say $\mu$ is \emph{contained} in $U$ and write $\mu \subseteq  U$.  We say a collection of multicurves \emph{fills} a subsurface $U$ if $U$ is the smallest isotopy class of subsurfaces that contains the collection of multicurves. Any collection of curves fills a unique subsurface which we call the \emph{fill} of that collection. Note, since we are considering an annulus to contain its core curve,  the fill of a single curve  $c$ is the annulus with core curve $c$.  If a curve $c$ is contained in $U$,  and $U$ is the annulus with core curve $c$, then $U -c = \emptyset$. Otherwise, $U- c$  is  the isotopy class of the subsurface obtained by deleting the representative of $c$ that is contained in a representative of $U$. We extend this to define $U -\mu$ for a multicurve that  is contained in the subsurface $U$.

 Two  subsurfaces are \emph{disjoint} if they have no connected components in common and have disjoint representative.  If the subsurfaces $U$ and $V$ are not disjoint and neither contains the other, we say $U$ and $V$ \emph{overlap}. 
 A curve $c$ and subsurface $U$ are \emph{disjoint} if the annulus with core curve $c$ is disjoint from $U$. This extends to define the disjointness of a multicurve and a subsurface. If a multicurve $\mu$ is not disjoint from a subsurface $U$, then we say $\mu$ and $U$ \emph{intersect}.  If the multicurve $\mu$ intersects a subsurface $U$, but is disjoint from $\partial U$, then some component of $\mu$ is contained in $U$ and we let $\mu \cap U$  denote this subset of components of $\mu$.

A \emph{pants subsurface} is any subsurface homeomorphic to $S_{0,p}^n$ with $p+n =3$. Given a subsurface $U$ of $S$ we call a multicurve that is contained in $U$ and has maximal cardinality among multicurves on $U$ a \emph{pants decomposition} of $U$. Note, when $U$ has any annular components, then every pants decomposition of $U$ contains the core curve for each of the annular components of $U$.

The \emph{mapping class group}, $\mcg(S)$, of a surface $S$ is the group of diffeomorphism of $S$ modulo isotopies that fix the marked points. 
If $z \in S$ is not a marked point, then $\mcg(S;z)$ is the subgroup of $\mcg(S^z)$ that fixes the marked point $z$. When $S \cong S_{g,p}^0$, then $\mcg(S;z) = \mcg(\Sdot)$. The \emph{point pushing subgroup}, $\Push(S;z)$, is the subgroup of $\mcg(S;z)$ comprised of all pushes of the point $z$ along loops in $S$ based at $z$. The Birman exact sequence tells us that $\Push(S;z)$ is naturally isomorphic to $\pi_1(S)$ and is the kernel of the surjective map $\mcg(S;z) \to \mcg(S)$ induced by forgetting the marked point $z$.

\subsection{Complexes of curves and subsurface projection} Let $S \cong S_{g,p}^n$. Given a non-annular subsurface $U \subseteq S$, the  \emph{curve complex}, $\C(U)$, of $U$ is the flag simplicial complex whose vertices are  curves on $U$ and with an edge between two distinct curves if they are disjoint. If $U$ is connected and has complexity 1, then all curves of $U$ intersect, so we modify the definition of the edges to allow for an edge between any two curves that intersect minimally on $U$---twice if $(g,p+n) = (0,4)$ and once if $(g,p+n) = (1,1)$.

When $U$ is an annulus, then $\C(U)$ has an alternative description and is often called the \emph{annular complex} instead of the curve complex. To define $\C(U)$ when $U$ is an annulus, fix a complete hyperbolic metric on $S$ where the marked points are viewed as punctures. Let $\widetilde{U}$ be the annular cover of $S$ corresponding to the annulus $U$ equipped with the lift of this hyperbolic metric on $S$. Let $\widehat{U}$ denote the compactification of $\widetilde{U}$ to a closed annulus obtained from the usual compactification of the hyperbolic plane.  The vertex set of $\C(U)$ is the set of isotopy classes (relative the boundary at infinity of $\widehat{U}$)  of  arcs that connect the two boundary components of $\widehat{U}$. Two vertices of the annular complex $\C(U)$ are then joined by an edge if they have representatives with disjoint interiors.

Masur and Minsky proved that $\C(U)$ is hyperbolic in all cases where it is non-empty \cite{MMI}. A number of other authors have since proved that the constant of hyperbolicity is independent of the surface $S$ or subsurface $U$ \cite{Aougab_Uniform_hyperbolicity,bowditch_uniform_hyperbolicity,CRS_Uniform_hyperbolicity,HPW_Uniform_Hyperbolicity}.

For each subsurface $U \subseteq S$, Masur and Minsky defined a subsurface projection map $\pi_U \colon \C^0(S) \to 2^{\C(U)}$. We direct the reader to \cite{MMII} for the full definition of the projection maps and instead recall the properties we will need for the present work.

\begin{lem}[{\cite[Lemma 2.2, 2.3]{MMII}}]\label{lem:subsurface_projection_coarsely_Lipschtiz}
	Let $c$ and $c'$ be curves on a surface $S$ and $U$ be a non-pants subsurface of $S$.
	\begin{itemize}
		\item $\pi_U(c) \neq \emptyset$ if and only if $c$ intersects $U$ and $c \neq \partial U$ when $U$ is an annulus.
		\item $\diam(\pi_{U}(c)) \leq 3$ whenever $\pi_U(c)\neq \emptyset$.
		\item If $c$ and $c'$ are disjoint,  then $d_{\C(U)}(\pi_U(\alpha),\pi_U(\beta) ) \leq 2$ whenever $\pi_U(c)$ and $\pi_U(c')$ are both non-empty.  If additionally $c' \subseteq U$, then every curve in $\pi_U(c)$ will not intersect $c'= \pi_U(c')$.
	\end{itemize}	
\end{lem}

\begin{notation}
	When $c$ and $c'$ are two curves with  $\pi_U(c) \neq \emptyset$ and $\pi_U(c') \neq \emptyset$,  we use $d_U(c,c')$ to denote $d_{\C(U)}( \pi_U(c),\pi_U(c'))$. As before this is the minimum distance between the sets $\pi_U(c)$ and $\pi_U(c')$.
\end{notation}

Distances in the curve complex are known to be bounded above in terms of the intersection number of curves; for example \cite{Bowditch_Guessing_Geodesic}. This means we can modify the edge relation in $\C(U)$ to allow for a uniformly bounded number of intersections between curves without changing the quasi-isometry type of $\C(U)$. In our case, it will sometimes be convenient to work with  the following variant of the curve complex that is quasi-isometric to $\C(U)$ and hence uniformly hyperbolic.

\begin{defn}[Modified curve graph]\label{defn:C'(S)}
	For a  non-annular, non-pants subsurface $U \subseteq S$, define $\C'(U)$ to be the graph with has vertex set $\C^0(U)$ and  with an edge between curves if they are either disjoint or intersect at most 4 times inside of a connected complexity 1 subsurface of $S$. 
\end{defn}

When $S$ is a closed surface and $z \in S$, the map $\Pi \colon S^z \to S$ induced by forgetting the marked point $z$ has the property that the image of every essential curve on $S^z$ is still essential on $S$. Thus, $\Pi$ induces a map $\Pi \colon \C^0(\Sdot) \to \C^0(S)$. A theorem of Kent, Leininger, and Schleimer proves that the preimage of a curve  under $\Pi$ is a tree.

\begin{thm}[{\cite[Theorem 7.1]{KLS_Trees_and_MCG}}]\label{thm:fibers_are_trees}
	Let $S \cong S_{g,0}^0$ with $g \geq 2$. For each $\alpha \in \C^{0}(S)$, the subset of $\C(\Sdot)$ spanned by  $\Pi^{-1}(\alpha)$ is a tree.
\end{thm}		

\begin{figure}[ht]
	\centering
	\def\svgscale{.7}
	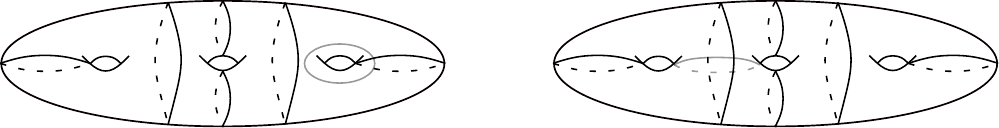
	\caption{Replacing the gray curve in each picture with the black curve it intersects gives an example of the two different types of flip moves in the pants graph.}
	\label{fig:flip}
\end{figure}

The \emph{pants graph}, $\mc{P}(S)$, of a surface $S$ is the graph whose vertices are all pants decomposition of $S$ and where two pants decompositions $\rho$ and $\rho'$ are joined by an edge if there exists curves $c \in \rho$ and $c' \in \rho'$ such that $(\rho-c) \cup c' = \rho'$ and $c$ and $c'$ are joined by an edge in the curve complex of the complexity 1 component of $S - (\rho - c)$. Pictorially, edges in $\mc{P}(S)$ correspond to one of the two ``flip'' moves shown in Figure \ref{fig:flip}. 

A \emph{marking}\footnote{What we call a marking, Masur and Minsky called a complete marking. As we will have no need for incomplete markings, we forgo the distinction.} $\mu$ of a surface $S$ is a set $\{(c_1,t_1),\dots,(c_k,t_k)\}$ where $c_1,\dots,c_k$ are curves that make a pants decomposition of $S$ and each $t_i$ is a vertex of the annular complex for the annulus whose core curve is $c_i$.
We use $\base(\mu)$ to denote the set  $\{c_1,\dots, c_k\}$ and $\markcurve(\mu)$ to denote $\{t_1,\dots,t_k\}$. A  marking $\mu$ is \emph{clean} if for each $i \in \{1,\dots, k\}$, the arc $t_i$  is the lift of a curve $c'_i$ such that $c'_i$ and $c_i$ fill a complexity 1 subsurface $U$ and $d_U(c'_i,c_i) = 1$.  
Given a marking $\mu$ and a non-pants subsurface $U$, Masur and Minsky defined the subsurface projection of $\mu$ as follows:
\begin{itemize}
	\item If $U$ is an annulus with core curve $c_i \in \base(\mu)$, then $\pi_U(\mu) = t_i$.
	\item If $U$ is any other subsurface, then  $\pi_U(\mu) = \pi_U(\base(\mu))$.
\end{itemize}

Masur and Minsky  defined the \emph{marking graph} $\mc{M}(S)$. This graph has the set of all clean markings of $S$ as vertices, with edges defined by a ``twist move'' and a ``flip plus clean-up'' move. The next result summarizes the properties of the marking graph that we will need.

\begin{thm}[\cite{MMII}]\label{thm:Marking_complex}
	The marking graph $\mc{M}(S)$ is a proper connected graph with the following properties:
	\begin{enumerate}
		\item \label{item:geometric} The action of $\mcg(S)$  on $S$ induces a  metrically proper and cocompact action of $\mcg(S)$ on $\mc{M}(S)$ by isometries.
		\item \label{item:distance_formula}There is a function $f \colon [0,\infty) \to [0,\infty)$ depending only on $S$ so that for any two clean markings $\mu,\nu \in \mc{M}(S)$ if $d_U(\mu,\nu)\leq r$ for each non-pants subsurface $U \subseteq S$, then $d_{\mc{M}(S)}(\mu,\nu) \leq f(r)$.
		\item \label{item:cleaning}  There is $D>0$ depending only on $S$ so that for each (not necessarily clean) marking $\mu$, there exists a set of clean markings $\cl(\mu)$  satisfying:
		\begin{itemize}
			\item  $\base(\mu') = \base(\mu)$  for each $\mu' \in \cl(\mu)$;
			\item if $\mu' \in \cl(\mu)$, then for each $(c_i,t_i)\in \mu$ and $(c_i,t_i') \in \mu'$, we have $d_{U_i}(t_i,t_i') \leq 3$ where $U_i$ is the annulus with core curve $c_i$;
			\item  $\cl(\mu)$ has diameter at most $D$ in $\mc{M}(S)$.
		\end{itemize}
	\end{enumerate}
\end{thm}

 The set of clear markings $\cl(\mu)$ defined in Item \ref{item:cleaning}  of Theorem \ref{thm:Marking_complex} are called the \emph{clean markings compatible with $\mu$}. Item \ref{item:distance_formula} is a special case of Masur and Minsky's celebrated distance formula for the marking graph; see \cite[Theorem 6.12]{MMII}.

\subsection{Combinatorial HHSs}\label{sec:cHHS}

Behrstock, Hagen, and Sisto's original definition of a hierarchically hyperbolic space required the construction of a large number of hyperbolic spaces and the verification of nine axioms \cite{BHS_HHSII}. Recently, Behrstock, Hagen, Martin, and Sisto introduced \emph{combinatorial hierarchically hyperbolic spaces} which reduces the construction of a hierarchically hyperbolic space to the construction of a pair of spaces $(X,W)$ \cite{BHMS_cHHS}. As we will not need to work directly with the full definition of a hierarchically hyperbolic space, we will forgo giving the definition and instead describe the combinatorial HHS machinery.

The starting place for a combinatorial hierarchically hyperbolic space is a pair of spaces $(X,W)$ where $X$ is a flag simplicial complex and $W$ is an graph whose vertices are all the maximal simplices of $X$. We call such a graph $W$ an \emph{$X$-graph}. To help illuminate the definitions we will maintain a running example where $X$ will be the curve complex $\C(S)$ and $W$ is the pants graph $\mc{P}(S)$; $\mc{P}(S)$ is a $\C(S)$-graph as pants decompositions of $S$ are exactly the  vertices of the maximal simplices in $\C(S)$.

The definition of a combinatorial HHS  $(X,W)$ includes a number of properties that the pair $(X,W)$ need to satisfy. Stating these properties requires a bit of set up. First, we need the following augmented version of $X$.

\begin{defn}[Augmented graph]
	If $W$ is an $X$-graph for the flag simplicial complex $X$, we define the \emph{$W$-augmented graph $X^{+W}$} as the graph with the same vertex set as $X$ and with two types of edges: 
	\begin{enumerate}
		\item ($X$-edge) If two vertices $x_1,x_2 \in X^0$ are joined by an edge in $X$, then $x_1$ and $x_2$ are joined by an edge in $X^{+W}$.
		\item ($W$-edge) If $\Delta_1$ and $\Delta_2$ are maximal simplices of $X$ that are joined by an edge in $W$, then each vertex of $\Delta_1$ is joined by an edge to each vertex of $\Delta_2$ in $X^{+W}$.
	\end{enumerate}
\end{defn}

\begin{example}
	For the case of $(\C(S),\mc{P}(S))$, the augmented graph $\C(S)^{+\mc{P}(S)}$ is a copy of the 1-skeleton of $\C(S)$ with additional edges between curves that intersect minimally in a complexity 1 subsurface.
\end{example}

Next we define an equivalence relation among simplices of $X$.

\begin{defn}
	Let $\Delta$ and $\Delta'$ be simplices of the flag simplicial complex $X$. We write $\Delta \sim \Delta'$ if $\lk(\Delta) = \lk(\Delta')$.
	We define the \emph{saturation} of $\Delta$, $\sat(\Delta)$, to be the set of vertices of $X$ contained in a simplex in the $\sim$-equivalence class of $\Delta$. That is $x \in \sat(\Delta)$  if and only if there exists $\Delta'  \sim \Delta$  so that $x$ is  a vertex of $\Delta'$.
\end{defn}

\begin{example}\label{ex:saturation}
	If $\Delta$ is a simplex of $\C(S)$, then $\Delta$ defines a subsurface $U_\Delta$ that is the disjoint union of all the non-pants components of $S - \Delta$.  The link of $\Delta$ in $\C(S)$ is then spanned by all curves on $U_\Delta$, that is, $\lk(\Delta) = \C(U_\Delta)$. Two simplices $\Delta$ and $\Delta'$  will therefore have equal links if and only if $U_\Delta = U_{\Delta'}$. Further, the saturation of $\Delta$ will be the set of curves of on $S - U_\Delta$ plus the curves in $\partial U_\Delta$, that is, $\sat(\Delta) = \C^0(S - U_\Delta) \cup \partial U_\Delta$. This is because every non-peripheral curve on $S-U_\Delta$ is part of a pants decomposition of $S-U_\Delta$, and the join of this pants decomposition with $\partial U_\Delta$ produces a simplex $\Delta'$ with $U_{\Delta'} = U_\Delta$.
\end{example}

Finally we define two  subspaces of $X^{+W}$ that are associated to every simplex of $X$.

\begin{defn}
	Let $X$ be a flag simplicial complex and let $W$ be an $X$-graph. For each simplex $\Delta$ of $X$, define
	\begin{itemize}
		\item  $Y_\Delta$  to be  $X^{+W} - \sat(\Delta)$;
		\item  $H(\Delta)$ to be the subgraph of $Y_\Delta$ spanned by the vertices in $\lk(\Delta)$.
	\end{itemize}
 For $H(\Delta)$, we are taking the link in $X$, not in $X^{+W}$, and then considering the subgraph of $Y_\Delta$ induced by those vertices. We give both $Y_\Delta$ and $H(\Delta)$ their intrinsic path metrics. By construction, we have $Y_\Delta = Y_{\Delta'}$ and  $H(\Delta) = H(\Delta')$ whenever $\Delta \sim \Delta'$.
\end{defn}

\begin{example}
	As discussed in Example \ref{ex:saturation}, for a simplex $\Delta$ of $\C(S)$,  the saturation of $\Delta$ is the set of curves on $S - U_\Delta$ plus the multicurve $\partial U_\Delta$. Thus $\C(S) - \sat(\Delta)$ is spanned by the set of all curves that  intersect $U_\Delta$.  The space $Y_\Delta$ is then $\C(S) - \sat(\Delta)$  with extra edges between curves that intersect minimally in a  complexity 1 subsurface of $S$. Because $\lk(\Delta) = \C(U_\Delta)$, the space $H(\Delta)$ is  a copy of $\C(U_\Delta)$ with additional edges between pairs of curves that intersect minimally in a complexity 1 subsurface of $U_\Delta$.
\end{example}

We can now state the definition of a combinatorial hierarchically hyperbolic space.

\begin{defn}[Combinatorial HHS]\label{defn:CHHS}
	Let $\delta \geq 0$, $X$ be a flag simplicial complex, and $W$ be an $X$-graph. The pair $(X,W)$  is a \emph{$\delta$-combinatorial HHS} if the following are satisfied.
	\begin{enumerate}
		\item \label{CHHS:finite_complexity} If $\Delta_1,\dots,\Delta_n$ are non-maximal simplices of $X$ with  $\lk(\Delta_1) \subsetneq \lk(\Delta_2) \subsetneq \dots \subsetneq \lk(\Delta_n)$, then $n \leq \delta$.
		\item \label{CHHS:hyp_X} The $W$-augmented graph $X^{+W}$ is connected and $\delta$-hyperbolic.
		
		\item \label{CHHS:geom_link_condition} For each non-maximal simplex $\Delta \subseteq X$, the space $H(\Delta)$ is $\delta$-hyperbolic and the inclusion $H(\Delta) \to Y_\Delta$ is a $(\delta,\delta)$-quasi-isometric embedding.
		
		\item \label{CHHS:comb_nesting_condition} Whenever $\Delta_1$ and $\Delta_2$  are non-maximal simplices for which there exists a non-maximal simplex $\Lambda$ such that $\lk(\Lambda) \subseteq \lk(\Delta_1) \cap \lk(\Delta_2)$ and $\diam(H(\Lambda)) \geq \delta$, then there exists a simplex $\Omega$ in the link of $\Delta_2$ so that $\lk(\Delta_2 \bowtie \Omega) \subseteq \lk(\Delta_1)$ and all simplices $\Lambda$ as above satisfy $\lk(\Lambda) \subseteq \lk(\Delta_2 \bowtie \Omega)$.
		
		\item \label{CHHS:C=C_0_condition} For each non-maximal simplex $\Delta$ and $x,y \in \lk(\Delta)$, if $x$ and $y$ are not joined by an $X$-edge of $X^{+W}$ but are joined by a $W$-edge of $X^{+W}$, then there exits simplices $ \Lambda_x, \Lambda_y$ so that $\Delta \subseteq \Lambda_x, \Lambda_y$ and $x \bowtie \Lambda_x$, $y \bowtie \Lambda_y$ are vertices of $W$ that are joined by an edge of $W$.
	\end{enumerate}
\end{defn}

\begin{remark}
	We will colloquially refer to conditions (\ref{CHHS:finite_complexity}), (\ref{CHHS:comb_nesting_condition}), and (\ref{CHHS:C=C_0_condition}) as the combinatorial conditions and conditions (\ref{CHHS:hyp_X}) and (\ref{CHHS:geom_link_condition}) as the geometric conditions.
\end{remark}

\begin{example}
	For $(\C(S),\mc{P}(S))$, the spaces $H(\Delta)$ are hyperbolic because they are quasi-isometric to $\C(U_\Delta)$. This is because distances in $\C(U_\Delta)$ are bounded above by the intersection number of curves. Because the space $Y_\Delta$ is spanned by curves that intersect the subsurface $U_\Delta$, subsurface projection gives a coarsely Lipschitz map $\pi_{U_\Delta} \colon Y_\Delta \to H(\Delta)$. This makes $H(\Delta)$ quasi-isometrically embedded in $Y_\Delta$ by Lemma \ref{lem:coarsel_lipschitz}.
	
	For the combinatorial conditions, observe that $\lk(\Delta) \subseteq \lk(\Delta')$ if and only if $U_\Delta \subseteq  U_{\Delta'}$. Thus, the finite complexity of $S$ gives a uniform bound on chains of properly nested links of simplices. For Item \ref{CHHS:comb_nesting_condition}, it suffices to verify the condition for simplices $\Delta_1$ and $\Delta_2$ where $U_{\Delta_1}$ and $U_{\Delta_2}$ overlap.  Let $Q$ be the largest subsurface contained in both $U_{\Delta_1}$ and $U_{\Delta_2}$, then let $\Gamma$ be a simplex so that $\Delta_2 \bowtie \Gamma$ is a pants decomposition of   $S -Q$, The desired simplex $\Omega$ will be $\Gamma$ joined with the set of boundary curves of $Q$  that are not curves of $\Delta_2$. For Item \ref{CHHS:C=C_0_condition},  if $x$ and $y$ are two curves joined by a $\mc{P}(S)$-edge but not a $\C(S)$-edge, then $x$ and $y$ fill a complexity 1 subsurface $V$. If $x$ and $y$ are in the link of $\Delta$, then $V\subseteq U_\Delta$. Hence, we can find $\Lambda_x$ and $\Lambda_y$ by taking the join of $\Delta$ with $\partial U_\Delta$ and pants decompositions of $U_\Delta$ that contain $x$ and $y$ respectively.
\end{example}

When $(X,W)$ is a combinatorial HHS, the space $W$ is a hierarchically hyperbolic space. Further, a group $G$ will be hierarchically hyperbolic if it acts as described below on both $X$ and $W$. In Section \ref{sec:structure}, we provide a brief summary of the salient parts of the definition of hierarchical hyperbolicity as well as a description of the hierarchically hyperbolic structure imparted on combinatorial  HHSs.

\begin{thm}[{\cite[Theorem 1.18, Remark 1.19]{BHMS_cHHS}}]\label{thm:CHHS_are_HHS}
	If $(X,W)$ is a $\delta$-combinatorial HHS, then $W$ is connected and a hierarchically hyperbolic space. Further, a finitely generated group $G$ will be a hierarchically hyperbolic group if 
	\begin{enumerate}
		\item $G$ acts on $X$ by simplicial automorphisms with finitely many orbits of links of simplices;
		\item the action of $G$ on maximal simplices of $X$ induces a metrically proper and cobounded action of $G$ on $W$ by isometries.
	\end{enumerate}
\end{thm}

\begin{remark}
	We emphasis that there is no a priori requirement that either $X$ or $W$ is connected. In the proof of Theorem \ref{thm:CHHS_are_HHS}, it is shown that the definition of a combinatorial HHS implies $W$ is connected;  see Sections 1 and 5.2 of \cite{BHMS_cHHS}. We will take advantage of this in Section \ref{sec:non-annular_case} to skip a direct proof of connectedness for our combinatorial HHS.
\end{remark}

\section{A sketch of the proof}\label{sec:sketch}
Before embarking on the work for our main result, we give an outline of the proof in the case where $\alpha$ is a single  curve on $S$. This case captures all of the key ideas of our proof while avoiding some technicalities that arise from $\alpha$ having multiple components. Recall, $\Sdot$ is the surface obtained from $S$ by adding $z$ as a marked point of $S$. Let $G_\alpha$ be the stabilizer of $\alpha$ in $\mcg(S)$ and $E_\alpha$ be the preimage of $G_\alpha$ in $\mcg(S;z)$.  Recall, $\Pi \colon \Sdot \to S$ is the map that forgets the marked point, and it induces a map $\Pi \colon \C^0(\Sdot) \to \C^0(S)$ because $S$ is closed.

Our first step is to build a combinatorial HHS $(X_\alpha,W_\alpha)$ where the vertex stabilizers of $E_\alpha$ on $W_\alpha$ are each generated by Dehn twists about disjoint curves.  
Let $F_\alpha$ be the subset of $\C(\Sdot)$ spanned by $\Pi^{-1}(\alpha)$. Theorem \ref{thm:fibers_are_trees} established that $F_\alpha$ is a tree. While $F_\alpha$ is a hyperbolic space that $E_\alpha$ acts on, this action only ``sees'' the action of the $\pi_1(S) \cong \Push(S;z)$ subgroup of $E_\alpha$; we need our hyperbolic space $X_\alpha$ to ``see'' the action of all of $E_\alpha$. To accomplish this, we define $X_\alpha$ to be the subset of $\C(\Sdot)$ spanned by the preimages of all curves in the star of $\alpha$ in $\C(S)$. That is $X^0_\alpha = \{c \in \C^0(\Sdot) : \Pi(c) \text{ is disjoint from or equal to } \alpha\}$. 

The rationale for this choice of $X_\alpha$ is that the star of $\alpha$ in $\C(S)$ plus the pants graph on the surface $S-\alpha$ make a combinatorial HHS that sees all of the action of $G_\alpha$ except for Dehn twists. Thus, by including the preimage of the whole star of $\alpha$, we capture both the $\pi_1(S)$ portion and the $G_\alpha$ portion of $E_\alpha$.

The key lemma we prove about $X_\alpha$ is that it is \emph{$F_\alpha$-guided}. This means  for each  vertex $x \in X_\alpha^0 - F_\alpha$, the set $\lk(x) \cap F_\alpha$ is a subtree $T_x$ of $F_\alpha$ and if $x,y$ are adjacent vertices of $X_\alpha$, then $T_x \cap T_y$ is non-empty. This allows for points in $X_\alpha$ to be connected by well behaved paths in $F_\alpha$ instead of  $X_\alpha$, and allows us to prove that  $X_\alpha$ is a quasi-tree using a variant of Manning's bottleneck criteria.

The maximal simplices of $X_\alpha$ are pants decomposition of $\Sdot$ that use only curves in $X_\alpha$. Our $X_\alpha$-graph $W_\alpha$ is therefore defined analogously to the pants graph: the vertices are pants decomposition using only curves in $X_\alpha$ with an edge between two pants decompositions that differ by flipping one curve to another curve that intersects it at most 4 times. The choice of intersection number 4 to define the edges is because $4$ is the smallest number of intersections that two distinct element of $F_\alpha$ can have inside of a 4-punctured sphere subsurface.

For a multicurve $\mu \subseteq X_\alpha$, we define a subsurface $U_\mu$ to be the subsurface filled by the curves in the $X_\alpha$-link of $\mu$. Unlike in the case of the pants graph, $U_\mu$ can have annular components. Because the edges of $X_\alpha$  correspond to disjointness of curves, the nesting of links in $X_\alpha$ can then be translated into the nesting of these $U_\mu$ subsurfaces in $\Sdot$. Thus, the combinatorial conditions for $(X_\alpha,W_\alpha)$ to be a combinatorial HHS follow from the topology of the $U_\mu$ subsurfaces similarly to how they are proved in the case of the pants graph.

For the geometric conditions, we split the links of multicurves into two cases. If $\mu$ is a simplex for $X_\alpha$ so that $\alpha \subseteq \Pi(\mu)$, then every curve on the subsurface $\Sdot -\mu$ is a curve in $X_\alpha$. Thus, the link of $\mu$ in $X_\alpha$ will be the  entire curve complex $\C(\Sdot-\mu)$. This means we can prove the quasi-isometric embedding and hyperbolicity conditions using the subsurface projection map $\pi_{\Sdot-\mu}$. On the other hand, if $\alpha \not \subseteq \Pi(\mu)$, then not every curve on $\Sdot-\mu$ will be a curve in $X_\alpha$ and  the link  of $\mu$ in $X_\alpha$ is some proper subcomplex of $\C(\Sdot-\mu)$.  This makes using the subsurface projection map $\pi_{\Sdot-\mu}$ problematic, as we  might get a curve that is not in $X_\alpha$ when we project a curve in $X_\alpha - \sat(\mu)$ to $\Sdot-\mu$.
Instead, we show that  $X_\alpha - \sat(\mu)$ is still guided by the fiber $F_\alpha$ because $\alpha \not \subseteq \Pi(\mu)$. This implies  $\lk(\mu)$ will be a quasi-tree like $X_\alpha$, and allows us to build a coarsely Lipschitz map onto $\lk(\mu)$ by first mapping a curve $x \in X^0_\alpha - \sat(\mu)$  into the fiber $F_\alpha$ and then taking a closest point projection in $F_\alpha$ to $\lk(\mu) \cap F_\alpha$.

Once we have proved $(X_\alpha,W_\alpha)$ is a combinatorial HHS, we need to modify both spaces to produce a combinatorial HHS where $E_\alpha$ has a metrically proper action. Since the $E_\alpha$-stabilizer of a vertex of $W_\alpha$ is generated by the Dehn twists around the curves in the pants decomposition, we need to modify $X_\alpha$ to account for these Dehn twists. We do this by ``blowing up'' each vertex of $X_\alpha$ as follows: for each curve $c \in X_\alpha$, we add as additional vertices the set $B(c)$ of vertices of the annular complex for the annulus with core curve $c$. We  add an edge between $c$ and each element of $B(c)$ as well as between each element of $B(c)$ and every element of $B(c') \cup \{c'\}$ whenever $c$ and $c'$ are disjoint. We use $B(X_\alpha)$ to denote this blow-up of $X_\alpha$. 

The maximal simplices of $B(X_\alpha)$ are now (not necessarily clean) markings on $\Sdot$, and we define the graph $B(W_\alpha)$ to have vertices all the markings that are maximal simplices of $B(X_\alpha)$ with edges corresponding to both twist and flip moves. Since the vertices of $B(W_\alpha)$ are markings instead of pants decompositions, the action of $E_\alpha$ on $B(W_\alpha)$ will be metrically proper. At almost all steps, the proof that $(B(X_\alpha),B(W_\alpha))$ is a combinatorial HHS reduces to the proof that $(X_\alpha,W_\alpha)$ is an combinatorial HHS. The only truly new conditions to be verified are the geometric conditions for the links of simplices $\mu$ where $\mu$ is a pants decomposition with all but one curve $c$ marked. The hyperbolic space associated to such a link is quasi-isometric to the annular complex of the annulus with core curve $c$, thus we can use the subsurface projection onto this annulus to verify the quasi-isometric embedding condition.

\section{Tree guided spaces}\label{sec:tree guided}
We now introduce the technique that we will use to verify that a graph or simplicial complex is a quasi-tree.

\begin{defn}\label{defn:tree_adjacent}
	Let $X$ be a graph and $F \subseteq X$ be a connected subgraph.  For each $x \in X^0$ define $L(x) = \{x\}$ if $x \in F$ and   $L(x) = \lk(x) \cap F$ if $x\not \in F$. We say $X$ is \emph{$F$-guided} if 
	\begin{enumerate}
		\item \label{tree_adj:L(x)_connected} for each $x \in X^0$, $L(x)$ is a non-empty connected subset of $F$;
		\item \label{tree_adj:chaining} if $x,y \in X^0$ are joined by an edge, then either $x$ and $y$ are joined by an edge of $F$ or $L(x) \cap L(y)$ is a non-empty connected subset of $F$.
	\end{enumerate}
	When $F$ is a tree and $X$ is $F$-guided, then we say $X$ is \emph{tree guided}.
\end{defn}

The key property imparted by  Definition \ref{defn:tree_adjacent} is that every path in the graph $X$ produces a connected subset of the guiding subgraph $F$ connecting the end points of the path.

\begin{defn}\label{defn:T_sequence}
Let $X$ be a graph and $F \subseteq X$ be a connected subgraph so that $X$ is $F$-guided. If $\gamma$ is an edge path in $X$ with vertices  $x_1,\dots,x_n$, the sequence $L(x_1),\dots, L(x_n)$ is called the \emph{$F$-sequence} for $\gamma$. Definition \ref{defn:tree_adjacent} implies that for any edge path of $X$, the union of the $L(x_i)$ in the $F$-sequence span a connected subset of $F$.
\end{defn}

The main result we will need is that tree guided graphs are quasi-isometric to trees.

\begin{lem}\label{lem:tree_adj_is_quasi-tree}
Let $X$ be a graph and $F \subseteq X$ be a connected subgraph. If $X$ is $F$-guided and $F$ is a tree, then $X$ is uniformly  quasi-isometric to a tree.
\end{lem}
To prove Lemma \ref{lem:tree_adj_is_quasi-tree}, 	we use the following variant of Manning's bottleneck criteria \cite{ManningBottleneck}. A proof of the this variant was given in Proposition 2.5 of \cite{DDLS_Veech_groups_II}.

\begin{prop}[\cite{ManningBottleneck,DDLS_Veech_groups_II}]\label{prop:bottleneck_variant}
	Let $X$ be a graph. Suppose that there exists a constant $R$ with the following property: for each pair of vertices $w, w' \in X$ there exists an edge path $\eta(w,w')$  from $w$ to $w'$ so that for any vertex $v$ of $\eta(w,w')$, any path from $w$ to $w'$ in $X$ intersects $\ball_R(v)$. Then $X$ is quasi-isometric to a tree, with quasi-isometry constants depending on $R$ only.
\end{prop}

\begin{proof}[Proof of Lemma \ref{lem:tree_adj_is_quasi-tree}]

	For each $x,y \in X^0$, let $\rho(x,y)$ be the shortest path in $F$ connecting $L(x)$ to $L(y)$ in $F$. Define $\eta(x,y)$ to be the edge path spanned by $x \cup \rho(x,y) \cup y$.  Let  $\gamma$ be any edge path of $X$ connecting $x$ and $y$. The  $F$-sequence for $\gamma$ gives a connected subset of $F$ that contains the end points of  $\rho(x,y)$. Since $F$ is a tree, the $F$-sequence of $\gamma$ must then contain all of $\rho(x,y)$. In particular, for each $v \in \rho(x,y)$, the path $\gamma$ must intersect the ball of radius $1$ around $v$. Hence, $X$ is  uniformly quasi-isometric to a tree by Proposition \ref{prop:bottleneck_variant}.	
\end{proof}

We will also need the following lemma that will allow us to preserve $F$-guidedness when removing a set of vertices that are distinct from $F$.

\begin{lem}\label{lem:T-guide_after_removal}
	Let $X$ be a graph and $F \subseteq X$ be a connected subgraph so that $X$ is $F$-guided. If $Z \subseteq X^0$ with $Z \cap F = \emptyset$, then $X-Z$ is $F$-guided 
\end{lem}

\begin{proof}
	Because $Z \cap F =\emptyset$, the set $L(x)$ is unchanged for each vertex $x \in X-Z$. Further, if $x_1$ and $x_2$ are adjacent vertices of $X-Z$, then they are also adjacent vertices of $X$ and  either $x_1,x_2 \in F$  or $L(x_1) \cap L(x_2)$ is a non-empty connected subset of $F$ as desired.
\end{proof}

We conclude by recording a more general version of Lemma \ref{lem:tree_adj_is_quasi-tree} that will we will not need, but may be useful for future work.
\begin{lem}
	Let $X$ be a graph and $F \subseteq X$ be a connected subgraph. If $X$ is $F$-guided, then $X$ is $(2,0)$-quasi-isometric to the graph obtained from $F$ by adding a vertex $v_x$ for each  $x \in X^0-F^0$ and connecting $v_x$ to every vertex of $L(x)$. Moreover, if $F$ is hyperbolic and each $L(x)$ is uniformly quasiconvex in $F$, then $X$ is hyperbolic.
\end{lem}

\begin{proof}
	Let $\widehat{F}$ denote the graph obtained from $F$ by adding a vertex $v_x$ for each $x \in X^0-F^0$ and connecting $v_x$ to every vertex of $L(x)$. Let $f \colon \widehat{F} \to X$ be given by $f(w) = w$ for $w \in F$ and $f(v_x) = x$ for $x \in X^0 - F^0$. 
	
	We first prove that $f$ is a quasi-isometry. Let $z_1=w_0,w_1,\dots,w_n =z_2$ be the vertices of the $\widehat{F}$-geodesic connecting $z_1,z_2 \in \widehat{F}^0$. Since $X$ is $F$-guided,  $f(w_0), f(w_1), \dots, f(w_n)$ are the vertices of a path in $X$ connecting $f(z_1)$ and $f(z_2)$. Thus, $d_X(f(z_1),f(z_2)) \leq d_{\widehat{F}}(z_1,z_2)$ for all $z_1, z_2 \in \widehat{F}$. For the other inequality, assume $f(z_1)$ and $f(z_2)$ are joined by an edge of $X$ while $z_1, z_2$ are not joined by an edge of $\widehat{F}$. This only occurs if $z_1 =v_x$ and $z_2=v_y$ for $x,y \in X^0 - F^0$. Since  $f(z_1)=v_x$ and $f(z_2)=v_y$ are joined by an edge of $X$, the second condition in the definition of a $F$-guided space ensures that $L(x) \cap L(y) \neq \emptyset$. But this means $d_{\widehat{F}}(x,y) = 2$. Hence, we have $ d_{\widehat{F}}( z_1,z_2) \leq 2 d_X(f(z_1),f(z_2))$ for all $z_1, z_2 \in \widehat{F}$.
	
	When $F$ is hyperbolic and each $L(x)$ is uniformly quasiconvex in $F$, Proposition 2.6 of \cite{Kapovich_Rafi_Hyp_Free_Factors} shows that  $\widehat{F}$ is hyperbolic with constant depending only on the hyperbolicity of $F$ and the quasiconvexity constant of the $L(x)$'s. Since $X$ is uniformly quasi-isometric to $\widehat{F}$, $X$ must be hyperbolic as well.
\end{proof}

\section{A combinatorial HHS without annuli}\label{sec:non-annular_case}
Fix a closed surface $S \cong S_{g,0}^0$ with $g \geq 2$ and let $\Sdot$ be the surface obtained from $S$ by  making $z \in S$ a marked point.  Fix a multicurve $\alpha$ on $S$ and let $G_\alpha$ be the stabilizer of $\alpha$ in $\mcg(S)$. Let $E_\alpha$ be the full preimage in $\mcg(S;z)$ of $G_\alpha$.   Recall that the map $\Pi \colon \C^0(\Sdot) \to \C^0(S)$ is the map on curves induced by removing the marked point of $\Sdot$. 

In this section, we  build a combinatorial HHS $(X_\alpha,W_\alpha)$ on which $E_\alpha$ acts with large vertex stabilizers. In Section \ref{sec:annuli}, we will modify $(X_\alpha,W_\alpha)$ to produce a combinatorial HHS denoted $(B(X_\alpha),B(W_\alpha))$ where $E_\alpha$ does have a metrically proper action on $B(W_\alpha)$. As we can reduce the proof for  $(B(X_\alpha),B(W_\alpha))$ to the case of $(X_\alpha,W_\alpha)$, this intermediate step allows us to present a simpler and more transparent proof. Further,  the spaces $X_\alpha$ and $W_\alpha$ are of intrinsic interest as they are the most direct analogues of the curve complex and pants graph for $E_\alpha$.

Hence forth,  we say that constants are \emph{uniform} if they do not depend on the surface $S$, the multicurve $\alpha$, or a simplex of $X_\alpha$. A \emph{uniform quasi-tree} is thus a quasi-tree where the quasi-isometry constants to a tree does not depend on any of these quantities.

\subsection{Definition of $X_\alpha$ and $W_\alpha$}

We define the hyperbolic space $X_\alpha$ by taking the preimage of the star of the multicurve $\alpha$ under the marked point forgetting map $\Pi \colon \C^0(\Sdot) \to \C^0(S)$.

\begin{defn}\label{defn:X}
	Define $X_\alpha$ to be the subgraph of $\C(\Sdot)$ spanned by the curves $c$ so that $\Pi(c)$ is either disjoint from or contained in $\alpha$ on $S$. If $\alpha_1, \dots, \alpha_m$ are the curves comprising $\alpha$, then define $F_i$ to be the subset of $X_\alpha \subseteq \C(\Sdot)$ spanned by $\Pi^{-1}(\alpha_i)$. We call $F_i$ the \emph{fiber over $\alpha_i$}. By Theorem \ref{thm:fibers_are_trees}, $F_i$ is a tree for each $i \in \{1,\dots,m\}$.
\end{defn}

The   simplices of $X_\alpha$ are precisely the multicurves on $\Sdot$ whose components are all curves in $X_\alpha$. Thus, we will freely go between a multicurve  of curves in $X_\alpha$ and the simplex its vertices span.

We will exploit the fact that the fibers $F_i$ are trees to prove $X_\alpha$ is itself a quasi-tree. Our key tool is the fact that whenever a simplex $\mu$ does not intersect a fiber $F_i$, then the link of $\mu$ intersects $F_i$ in a subtree of $F_i$. 

\begin{lem}\label{lem:links_are_subtrees_of_fiber}
	Let $\mu$ be a simplex in $X_\alpha$, and suppose there exists a fiber $F_i$ so that $\mu \cap F_i = \emptyset$. For all $a_1,a_2 \in \lk(\mu) \cap F_i$,  there is a path in $\lk(\mu) \cap F_i$ connecting $a_1$ to $a_2$.
\end{lem}

\begin{proof}
	Let $\mu \subseteq X_\alpha$ be a simplex and $F_i$ be a fiber so that $\mu \cap F_i = \emptyset$. This means $\Pi(\mu)$ is disjoint from $\alpha_i$ on $S$ ensuring that $\lk(\mu) \cap F_i \neq \emptyset$. Let $a_1,a_2 \in \lk(\mu) \cap F_i$. We prove by induction on the intersection number $i(a_1,a_2)$ that $a_1$ is connected to $a_2$ in $\lk(\mu) \cap F_i$.
	
	The base case of $i(a_1,a_2)=0$  is automatic since $i(a_1,a_2) =0$ implies $a_1$ is joined to $a_2$ by an edge of $F_i$ and this edge is in $\lk(\mu) \cap F_i$ since $a_1,a_2 \in \lk(\mu)$.
	
	Assume that the lemma holds for curves of $\lk(\mu)\cap F_i$ that intersect strictly less than $n$ times and suppose $0<i(a_1,a_2) = n$. Since $a_1$ and $a_2$ are not disjoint and are both elements of $F_i$, we know $a_1 \cup a_2$ must form a bigon that contains the marked point of $\Sdot$. In particular, there exists an innermost bigon of $a_1 \cup a_2$ around the marked point. We can surger $a_2$ across this inner most bigon as shown in Figure \ref{fig:surger_across_puncture} to produce a curve $a_2'$ so that $a_2'$ is disjoint from $a_2$ and $i(a_1,a_2') < n$. 
	\begin{figure}[ht]
		\centering
		\def\svgscale{1}
		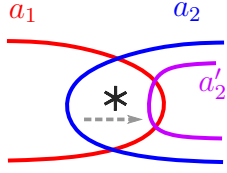
		\caption{The surgery of a curve across a bigon containing the marked point.}
		\label{fig:surger_across_puncture}
	\end{figure}
	Now, $a_2'$ is disjoint from $\mu$ on $\Sdot$, since $a_2$ is disjoint from $\mu$ on $\Sdot$. Further, $a_2' \in F_i$ since $a_2'$ would be isotopic to $a_2$ if the marked point of $\Sdot$ was removed. Thus, $a_2' \in \lk(\mu) \cap F_i$ with $i(a_1,a_2')<n$ and the induction hypothesis implies the path in $F_i$ from $a_1$ to $a_2'$ is contained in $\lk(\mu)$. Since $i(a_2,a_2') = 0$, this gives a path in $\lk(\mu) \cap F_i$ from $a_1$ to $a_2$.
\end{proof}

Using Lemma \ref{lem:links_are_subtrees_of_fiber}, we can prove $X_\alpha$ is a quasi-tree by proving it is $F_i$-guided for each fiber $F_i$.

\begin{prop}\label{prop:X_is_tree_adj}
	For each $i \in \{1,\dots,m\}$, the graph $X_\alpha$ is $F_i$-guided. In particular, $X_\alpha$ is a uniform quasi-tree.
\end{prop}

\begin{proof}
	Since each fiber $F_i$ is a tree, the final clause follows from the first clause by  Lemma \ref{lem:tree_adj_is_quasi-tree}.
	
	Fix a fiber $F_i$ over $\alpha_i$. For $x \in X^0_\alpha$, define $L_i(x) = \{x\}$ if $x \in F_i$ and $L_i(x) = \lk(x) \cap F_i$ if $x \not \in F_i$. By definition of $X_\alpha$, $L_i(x)$ is a non-empty subset of $F_i$ for all $x \in X^0_\alpha$, and by Lemma \ref{lem:links_are_subtrees_of_fiber}, $L_i(x)$ is connected.
	We now verify the second requirement of being $F_i$-guided: that $x$ being adjacent to $y$ in $X_\alpha$ implies $x,y \in F_i$ or $L_i(x) \cap L_i(y)$ is  non-empty and connected.
	
	Let $x,y \in X^0_\alpha$ be joined by an edge in $X_\alpha$. Thus, $x$ and $y$ are disjoint curves on $\Sdot$. Without loss of generality, assume $x \not \in F_i$. If $y \in F_i$, then $y \in L_i(x)$ by definition and we are done because $L_i(x) \cap L_i(y) =\{y\}$. We can therefore assume $y \not \in F_i$ as well. Similar to the proof of Lemma \ref{lem:links_are_subtrees_of_fiber}, we will use induction on intersection number to prove that there exists a curve $a$ in $L_i(x) \cap L_i(y)$. Since $F_i$ is a tree and $L_i(x)$ and $L_i(y)$ are already connected, this suffices to prove $L_i(x) \cap L_i(y)$ is connected.
	
	Let $a$ be any curve in $L_i(x)$. If $i(a,y) = 0$, then $a \in L_i(y)$ by definition, so suppose $i(a,y) >0$. Since $a \in F_i$ and $y \not \in F_i$, we must have that $\Pi(y)$ and $\Pi(a)$ are disjoint on $S$. Thus, $y$ and $a$ must form an innermost bigon  that contains the marked point of $\Sdot$. We can surger $a$ across the marked point, to produce a curve $a'$ with $i(a',y) < i(a,y)$ and $a'$ disjoint from $x$ (this surgery is identical to the surgery in Figure \ref{fig:surger_across_puncture} with $y= a_1$, $a_2 = a$ an $a_2'=a'$). Since $a'$ is obtained from $a$ by surgery across the marked point, $\Pi(a) = \Pi(a')$. This means $a' \in F_i$ and we have $a' \in L_i(x)$. Since  $i(a',y) < i(a,y)$, we can now induct on the intersection number to find $a'' \in L_i(x)$ with $i(a'',y) =0$. Such a curve $a''$ will be in $L_i(y)$ by the definition of edges in $X_\alpha$.
\end{proof}

The next lemma verifies that the maximal simplices of $X_\alpha$ are exactly the pants decomposition of $\Sdot$ consisting of curves in $X^0_\alpha$. We  give a more technical version of this fact that we will need later.

\begin{lem}\label{lem:Extending_To_Pants_decomposition}
	Let $Q$ be any subsurface of $\Sdot$ and $\gamma$ be a (possibly empty) multicurve on $Q$.
	If both $\gamma$ and $\partial Q$ are simplices of $X_\alpha$, then there exists a pants decomposition $ \tau$ of $Q$ so that $\gamma \subseteq \tau \subseteq X^0_\alpha$. 
\end{lem}

\begin{proof}
	First we show that we can assume $Q$ has no annular components without losing any generality. Since $\partial Q \subseteq X_\alpha$, the core curve of every annular component of $Q$ is a curve of $X_\alpha^0$. Thus, if  $Q_0$ is the union of the non-annular components of $Q$ and $\tau_0 \subseteq X_\alpha$ is a pants decomposition of $Q_0$ containing $\gamma \cap Q_0$, then we can define $\tau$ to be the union of $\tau_0$ and the core curves of the annular components of $Q$. Hence, we can assume $Q$ has no annular component.
	
In addition to assuming $Q$ has no annular component, we can also assume $Q$ has no pants components as pairs of pants do not contain any curves. These two assumptions mean that $\Pi(Q)$ does not contain any annular components.
	
	If  $\Pi(Q) \subseteq S -\alpha$, then every curve on $Q$ must be a curve in $X_\alpha$. This means any pants decomposition of $Q$ that contains $\gamma$ is a simplex of $X_\alpha$.
	
	Assume then that $\Pi(Q) \not \subseteq S-\alpha$. Since $\partial Q$ is a simplex of $X_\alpha$, no component of $\alpha$ will intersect $\partial \Pi(Q) = \Pi(\partial Q)$. Because $\Pi(Q)$ has no annular components, the only way for  $\Pi(Q) \not \subseteq S-\alpha$ without $\alpha$ intersecting $\partial \Pi(Q)$ is for a curve of $\alpha$ to be contained in $\Pi(Q)$. 
	
	If $Q$ does not contain the marked point $z$, then $\Pi$ restricted to $Q$ is a homeomorphism. Hence, for each curve  $\alpha_i \in \alpha \cap \Pi(Q)$ there is exactly one curve $a_i$ on $Q$ with $\Pi(a_i) = \alpha_i$. Thus, there is a multicurve $\eta$ that contains all of the curves of $\Pi^{-1}(\alpha)$ that are on $Q$. Moreover,  every curve of $Q - \eta$ is a curve in $X_\alpha$. Therefore, any  pants decomposition of $Q$ that contains $\eta$ will be a simplex of $X_\alpha$. Since $\gamma \subseteq X_\alpha$, $\gamma$ cannot intersect $\eta$. Hence any pants decomposition $\tau$ that contains $\gamma \cup \eta$ will suffice.

	Finally, assume  $Q$ does contain the marked point $z$. If $Q -\gamma$ has a pants component $P$ that contains the marked point, then $\Pi$ restricted to $Q - P$ is a homeomorphism. Hence, the same reasoning as the previous paragraph yields a pants decomposition $\tau'$ of $Q-P$ that is a simplex of $X_\alpha$ and contains $\gamma \cap (Q-P)$. The desired simplex is then $\tau = \tau' \cup \gamma$.  If $Q -\gamma$ does not have a pants component that contains the marked point, then for some $\alpha_i \in \alpha \cap \Pi(Q)$  there exist $a_1,a_2 \in \Pi^{-1}(\alpha_i)$ so that $a_1,a_2 \subseteq Q$ and $Q - (a_1\cup a_2)$ has a pants component $P$ that   contains  the marked point. We can additionally choose $a_1$ and $a_2$ so that  every curve of $\gamma$ is either equal to or disjoint from each $a_i$. As before, there is a pants decomposition $\tau'$ of $Q-P$ that is a simplex of $X_\alpha$ and contains $\gamma \cap (Q-P)$. In this case $\tau = \tau' \cup a_1 \cup a_2$ is the desired pants decomposition.
\end{proof}

 We now define the $X_\alpha$-space $W_\alpha$. 

\begin{defn} \label{defn:W}
	Define $W_\alpha$ to be the the following graph.
	\begin{itemize}
		\item Vertices: Maximal simplices of $X_\alpha$.
		\item Edges: Distinct vertices $\rho$ and $\rho'$ are joined by an edge if there exists curves $c\in \rho$ and $d \in \rho'$ so that $i(c,d) \leq 4$ and $\rho = (\rho'- d) \cup c$.
	\end{itemize}
\end{defn}

By Lemma \ref{lem:Extending_To_Pants_decomposition}, the vertices of $W_\alpha$ are  all pants decomposition of $\Sdot$ that are built from curves in $X_\alpha$. Two such pants decomposition are joined by an edge precisely when they differ by  at most one curve that intersect at most 4 times. Hence, $W_\alpha$ is an analogue of the pants graph for the group $E_\alpha$.

\subsection{Combinatorial conditions for $(X_\alpha,W_\alpha)$}
This subsection is devoted to verifying that $(X_\alpha,W_\alpha)$ satisfies the combinatorial properties of the definition of combinatorial HHS (Items  \ref{CHHS:finite_complexity}, \ref{CHHS:comb_nesting_condition}, and \ref{CHHS:C=C_0_condition} of Definition \ref{defn:CHHS}). We start by explaining how the combinatorics of links of simplices of $X_\alpha$ are encoded in the topology of subsurfaces of $\Sdot$. Through out this subsection, $\lk(\cdot)$ will denote links in the simplicial complex $X_\alpha$.

\begin{defn}[Subsurface for a multicurve]\label{defn:U_mu}
	Given a simplex $\mu$ of $X_\alpha$, define $U_\mu$ to be the (possibly disconnected) subsurface filled by  the curves of $X_\alpha$ that are disjoint from $\mu$. That is, $U_\mu$ is the subsurface filled by the curves in $\lk(\mu)$.
\end{defn}

\begin{figure}[ht]
	\centering
\begingroup%
  \makeatletter%
  \providecommand\color[2][]{%
    \errmessage{(Inkscape) Color is used for the text in Inkscape, but the package 'color.sty' is not loaded}%
    \renewcommand\color[2][]{}%
  }%
  \providecommand\transparent[1]{%
    \errmessage{(Inkscape) Transparency is used (non-zero) for the text in Inkscape, but the package 'transparent.sty' is not loaded}%
    \renewcommand\transparent[1]{}%
  }%
  \providecommand\rotatebox[2]{#2}%
  \newcommand*\fsize{\dimexpr\f@size pt\relax}%
  \newcommand*\lineheight[1]{\fontsize{\fsize}{#1\fsize}\selectfont}%
  \ifx\svgwidth\undefined%
    \setlength{\unitlength}{192.00001153bp}%
    \ifx\svgscale\undefined%
      \relax%
    \else%
      \setlength{\unitlength}{\unitlength * \real{\svgscale}}%
    \fi%
  \else%
    \setlength{\unitlength}{\svgwidth}%
  \fi%
  \global\let\svgwidth\undefined%
  \global\let\svgscale\undefined%
  \makeatother%
  \begin{picture}(1,0.1)%
    \lineheight{1}%
    \setlength\tabcolsep{0pt}%
    \put(0,0){\includegraphics[width=\unitlength,page=1]{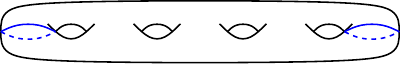}}%
    \put(1.05,.06){\color[rgb]{0,0,1}\makebox(0,0)[lt]{\lineheight{1.25}\smash{\begin{tabular}[t]{l}$\alpha$\end{tabular}}}}%
  \end{picture}%
\endgroup%

\begingroup%
  \makeatletter%
  \providecommand\color[2][]{%
    \errmessage{(Inkscape) Color is used for the text in Inkscape, but the package 'color.sty' is not loaded}%
    \renewcommand\color[2][]{}%
  }%
  \providecommand\transparent[1]{%
    \errmessage{(Inkscape) Transparency is used (non-zero) for the text in Inkscape, but the package 'transparent.sty' is not loaded}%
    \renewcommand\transparent[1]{}%
  }%
  \providecommand\rotatebox[2]{#2}%
  \newcommand*\fsize{\dimexpr\f@size pt\relax}%
  \newcommand*\lineheight[1]{\fontsize{\fsize}{#1\fsize}\selectfont}%
  \ifx\svgwidth\undefined%
    \setlength{\unitlength}{192.00001153bp}%
    \ifx\svgscale\undefined%
      \relax%
    \else%
      \setlength{\unitlength}{\unitlength * \real{\svgscale}}%
    \fi%
  \else%
    \setlength{\unitlength}{\svgwidth}%
  \fi%
  \global\let\svgwidth\undefined%
  \global\let\svgscale\undefined%
  \makeatother%
  \begin{picture}(1,0.23237929)%
    \lineheight{1}%
    \setlength\tabcolsep{0pt}%
    \put(0,0){\includegraphics[width=\unitlength,page=1]{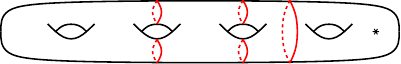}}%
    \put(1.05,.06){\color[rgb]{1,0,0}\makebox(0,0)[lt]{\lineheight{1.25}\smash{\begin{tabular}[t]{l}$\mu$\end{tabular}}}}%
  \end{picture}%
\endgroup%

\begingroup%
  \makeatletter%
  \providecommand\color[2][]{%
    \errmessage{(Inkscape) Color is used for the text in Inkscape, but the package 'color.sty' is not loaded}%
    \renewcommand\color[2][]{}%
  }%
  \providecommand\transparent[1]{%
    \errmessage{(Inkscape) Transparency is used (non-zero) for the text in Inkscape, but the package 'transparent.sty' is not loaded}%
    \renewcommand\transparent[1]{}%
  }%
  \providecommand\rotatebox[2]{#2}%
  \newcommand*\fsize{\dimexpr\f@size pt\relax}%
  \newcommand*\lineheight[1]{\fontsize{\fsize}{#1\fsize}\selectfont}%
  \ifx\svgwidth\undefined%
    \setlength{\unitlength}{192.00001153bp}%
    \ifx\svgscale\undefined%
      \relax%
    \else%
      \setlength{\unitlength}{\unitlength * \real{\svgscale}}%
    \fi%
  \else%
    \setlength{\unitlength}{\svgwidth}%
  \fi%
  \global\let\svgwidth\undefined%
  \global\let\svgscale\undefined%
  \makeatother%
  \begin{picture}(1,0.23237929)%
    \lineheight{1}%
    \setlength\tabcolsep{0pt}%
    \put(0,0){\includegraphics[width=\unitlength,page=1]{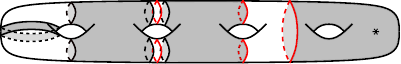}}%
    \put(1.05,.06){\makebox(0,0)[lt]{\lineheight{1.25}\smash{\begin{tabular}[t]{l}$U_\mu$\end{tabular}}}}%
  \end{picture}%
\endgroup%

	\caption{An example of the subsurface $U_\mu$. The top picture shows the fixed multicurve $\alpha$ on the closed surface $S$. The second picture is of the multicurve $\mu \subseteq X_\alpha$ on the surface $\Sdot$ (where $\ast$ denotes the marked point $z$). The bottom picture is the subsurface $U_\mu$. The core curve of the annulus in the bottom picture is $\alpha_\mu$.}
	\label{fig:U_mu}
\end{figure}
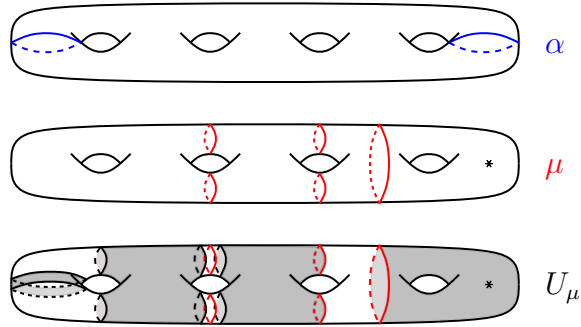

The subsurface $U_\mu$ is automatically a subsurface of  $\Sdot -\mu$, but $U_\mu$  might be strictly smaller than the collection of non-pants subsurfaces of $\Sdot -\mu$. To describe $U_\mu$, let $M_\mu$ be the component of $\Sdot -\mu$ that contains the marked point of $\Sdot$ then let $\alpha_\mu$ be the set of curves $\{ a \in \Pi^{-1}(\alpha) : a \subseteq (\Sdot - \mu) -M_\mu\}$. Since $(\Sdot - \mu) -M_\mu$ does not contain the marked point, $\alpha_\mu$ is a multicurve on $\Sdot - \mu$.  The next lemma asserts that $U_\mu$ is the union of the non-pants components of $\Sdot -(\mu \cup \alpha_\mu)$ and the  annuli whose core curves are in $\alpha_\mu$. Figure \ref{fig:U_mu} gives an example of $U_\mu$ that illustrates the lemma.

\begin{lem}\label{lem:description_of_U}
	Let $\mu$ be a simplex of $X_\alpha$. In the notation of the proceeding paragraph, $U_\mu$ is the disjoint union of the non-pants components of $\Sdot -(\mu \cup \alpha_\mu)$ plus the annuli whose core curves are in $\alpha_\mu$. Moreover, $\partial U_\mu \subseteq X_\alpha$ since every curve of $\partial U_\mu$ is either a curve of $\mu$ or a curve of $\alpha_\mu$.
\end{lem}

\begin{proof}	
	Let $V_\mu$ be the disjoint union of the non-pants components of $\Sdot -(\mu \cup \alpha_\mu)$ plus the annuli whose core curves are in $\alpha_\mu$. The curves in $\alpha_\mu \cup \lk(\mu \cup \alpha_\mu)$ fill $V_\mu$.  Since $U_\mu$ is  filled by the curves in $\lk(\mu)$, it suffices to prove $\lk(\mu) = \alpha_\mu \cup \lk(\mu \cup \alpha_\mu)$. Since $\alpha_\mu \subseteq \lk(\mu)$  we have $\alpha_\mu \cup \lk(\mu \cup \alpha_\mu) \subseteq \lk(\mu)$. For the other direction, let $c \in \lk(\mu)$ and let $K$ be the component of $\Sdot-\mu$ that contains $c$. If $K = M_\mu$, then $c$ is disjoint from each curve in $\alpha_\mu$ making $c \in \alpha_\mu \cup \lk(\mu \cup \alpha_\mu)$. If $K \neq M_\mu$, then  $\Pi$ restricted to $K$ is a homeomorphism. Hence every curve of $X_\alpha$ that is on $K$ must be either an element of $\alpha_\mu \cap K$ or disjoint from each curve in $\alpha_\mu \cap K$. Thus, either $c \in \alpha_\mu$ or $c \in \lk(\mu \cup \alpha_\mu)$.
\end{proof}

Using the subsurfaces $U_\mu$, we can relate the combinatorics of the links in $X_\alpha$ to the topology of subsurfaces on $\Sdot$.

\begin{lem}[Topology to cHHS dictionary]\label{lem:top_CHHS_dictionary}
	Let $\mu$ and $\mu'$ be non-maximal simplices of $X_\alpha$.
	\begin{enumerate}
		\item \label{dictionary:link} $\lk(\mu)$ is spanned by the set of all  curves $\{ x \in X^0_\alpha: x \subseteq U_\mu\}$. In particular, $\lk(\mu)$ is a join whenever  $U_\mu$ is disconnected.
		\item \label{dictionary:nesting} $\lk(\mu) \subseteq \lk(\mu')$ if and only if $U_\mu \subseteq U_{\mu'}$.
		\item \label{dictionary:equivalenet} $\lk(\mu) = \lk(\mu')$  if and only if $U_\mu = U_{\mu'}$.
		\item \label{dictionary:saturation} $\sat(\mu)$ is the set of all curves of $X_\alpha$ on $\Sdot - U_\mu$ plus the set of curves in $\partial U_\mu$ that are not  core curves of any annular components of $U_\mu$. 
	\end{enumerate}
\end{lem}

\begin{proof} 
	\begin{enumerate}
		\item  The first item is a direct consequence of the fact that edges in $X_\alpha$ correspond to the disjointness of curves plus the definition of $U_\mu$ as the fill of $\lk(\mu)$.
		\item Since $U_\mu$ and $U_{\mu'}$ are the subsurfaces filled by the curves in $\lk(\mu)$ and $\lk(\mu')$ respectively, $\lk(\mu) \subseteq \lk(\mu')$ implies $U_\mu \subseteq U_{\mu'}$. 
		
		On the other hand, if $U_\mu \subseteq U_{\mu'}$ and $c$ is a curve in $\lk(\mu)$, then  $c \subseteq U_\mu \subseteq U_{\mu'}$. This makes  $c \in \lk(\mu')$ by the first item. 	
		
		\item This item is an immediate consequence of the second item.
		\item Let $\tau$ be the subset of curves in $\partial U_\mu$ that are not core curves of any annular components of $U_\mu$. By the first item, if a curve $c$ is in $\sat(\mu)$, then $c$ must be disjoint from $U_\mu$.  Hence either $c \subseteq \Sdot-U_\mu$ or $c \in \tau$. To see that $\sat(\mu)$ contains every  curve $c \in X_\alpha^0$ that is a curve on $\Sdot- U_\mu$, observe that Lemma \ref{lem:Extending_To_Pants_decomposition} provides a pants decomposition $\rho$ of $\Sdot-U_\mu$ that contains $c$ and uses only curves in $X_\alpha$. Thus, $\rho \bowtie \tau$ is a simplex of $X_\alpha$ with $\lk(\rho \bowtie \tau) = \lk(\mu)$ because $U_{\rho \bowtie \tau} = U_\mu$.  Hence $c \in \sat(\mu)$. \qedhere
	\end{enumerate}
\end{proof}

Armed with Lemma \ref{lem:top_CHHS_dictionary}, we now verify the combinatorial criteria for $(X_\alpha,W_\alpha)$ to be a combinatorial HHS.

\begin{prop}\label{prop:comb_nesting_condition}
	Suppose $\mu$ and $\mu'$ are non-maximal simplices of $X_\alpha$ so that there exists a non-maximal simplex $\nu$ of $X_\alpha$ with $\lk(\nu) \subseteq \lk(\mu) \cap \lk(\mu')$ and $ \lk(\nu)$ not a join or a single vertex. There exists a (possibly empty) simplex $\rho \subseteq \lk(\mu')$ so that $\lk(\mu' \bowtie \rho) \subseteq \lk(\mu)$ and if $\nu$ is any simplex as in the proceeding sentence, then $\lk(\nu) \subseteq \lk(\mu' \bowtie \rho)$.	
\end{prop}

\begin{proof}

		Assume $\mu$ and $\mu'$ are simplices of $X_\alpha$ as described in the proposition. Let $V_0$ be the subsurface of $\Sdot$ filled by the curves in $\lk(\mu) \cap \lk(\mu')$, then let $V$ be the disjoint union of the non-annular components of $V_0$.  Note,  if $V = \emptyset$, then $\lk(\mu) \cap \lk(\mu')$ is either empty or is a simplex of $X_\alpha$. However, the assumption that there is a non-maximal simplex $\nu$ of $X_\alpha$ so that $\lk(\nu)$ is not a join or a single vertex and $\lk(\nu) \subseteq \lk(\mu) \cap \lk(\mu')$ prevents 	$\lk(\mu) \cap \lk(\mu')$ from being empty or a simplex. Thus, we know that $V \neq \emptyset$ and contains an infinite number of curves of $X_\alpha$.
		
	Suppose there is a simplex $\rho \subseteq \lk(\mu')$ with the property that $U_{\mu' \bowtie \rho} = V$  and let   $\nu$ be a non-maximal simplex of $X_\alpha$ so that 
		\begin{enumerate}
			\item $\lk(\nu) \subseteq \lk(\mu) \cap \lk(\mu')$;
			\item $\lk(\nu)$ is not a join or a single vertex.
		\end{enumerate}
		\noindent Then $U_\nu \subseteq V_0$ because $\lk(\nu) \subseteq \lk(\mu) \cap \lk(\mu')$. Since $\lk(\nu)$ is not a join, no curve in $\lk(\nu)$ can be a core curve for an annular component of $V_0$. Hence $U_\nu \subseteq V = U_{\mu' \bowtie \rho}$ which means $\lk(\nu) \subseteq \lk(\mu' \bowtie \rho)$ by Lemma \ref{lem:top_CHHS_dictionary}. Since $U_{\mu' \bowtie \rho} = V \subseteq U_{\mu}$ implies $\lk(\mu'\bowtie \rho ) \subseteq \lk(\mu)$, the simplex $\rho$ would satisfy the proposition.
		
		The proof will therefore be complete if we can find a simplex $\rho \subseteq \lk(\mu')$ so that $U_{\mu' \bowtie \rho} = V$. We will find  $\rho$ by applying Lemma \ref{lem:Extending_To_Pants_decomposition} to the subsurface $U_{\mu'} - V$. In order to do so,  we need to verify  that the boundary curves of $V$ are also curves in $X_\alpha$.

		\begin{claim}\label{claim:boundary_V_is_in_X}
			$\partial V$ is a simplex of $X_\alpha$.
		\end{claim}
	
		\begin{subproof}
			For the purposes of contradiction, assume  $\partial V \not \subseteq X^0_\alpha$. This means there exists $a \in \Pi^{-1}(\alpha)$ so that $\Pi(a)$ and $\Pi(\partial V)$ intersect. Further, by performing surgery across the marked point as in Figure \ref{fig:surger_across_puncture}, we can pick $a$ so that $a$ and $\partial V$ do not form any bigons around the marked point.  Pick a representative of $a$ that minimizes the intersection number with $\partial V$ and let $a'$ be the set of arcs produced by intersecting that representative with $V$.
			
			First assume  that $V$ does not contain the marked point $z$. In this case, the forgetful map $\Pi$ restricted to $V$ is a homeomorphism. Thus no curve  of $\lk(\mu) \cap \lk(\mu')$ on $V$ can intersect $a'$. Hence $V - a'$ is a proper subsurface of $V$ that contains all the curves of $\lk(\mu) \cap \lk(\mu')$ that are not the core curves of an annular component of $V_0$. However this is impossible as $V$ is the subsurface filled by these curves. Hence, $\partial V$ must be a simplex of $X_\alpha$.
			
			Now assume that $V$ does contain the marked point $z$.	Because $a$ is not contained in $V_0$, it must be the case that $a$ intersects either $\mu$ or $\mu'$. Without loss of generality, assume $a$ intersects $\mu$. Since $\Pi(a)$  does not intersect $\Pi(\mu)$, $a$ and $\mu$ must form an inner most bigon $B$ around the marked point. Because $V$ contains the marked point, the bigon $B$ must intersect the subsurface $V$. Since $a$ and $\partial V$ do not intersect in any bigon around the marked point, the only way for $B$ to intersect $V$ is for $\mu$ to also intersects $V$. However, this is impossible because $V$ is filled by curves from $\lk(\mu)$. Hence $\partial V$ must be a simplex of $X_\alpha$
		\end{subproof}

		Since both $\partial V$ and $\partial U_{\mu'}$ are both simplices  of $X_\alpha$, the boundary of $U_{\mu'} - V$ is also a simplex of $X_\alpha$. Hence,  Lemma \ref{lem:Extending_To_Pants_decomposition} provides a pants decomposition $\tau$ of $U_{\mu'} - V$ made up of curves of $X_\alpha$. 
		  Define $\rho$ to be the simplex spanned by $\tau \cup  \partial V$. 
		
		\begin{claim}\label{claim:U_rho=V}
			$U_{\mu' \bowtie \rho} = V$.
		\end{claim}
		
		\begin{subproof}
		 Every curve on $V$ is disjoint from $\mu' \cup \rho$ by construction. Since $V$ is filled by the curves of $X_\alpha$ that it contains and $U_{\mu' \bowtie \rho}$ is the subsurface filled by $\lk(\mu' \bowtie \rho)$, this means $V \subseteq U_{\mu' \bowtie \rho}$. 
		 
		 For the other direction, consider $c \in \lk(\mu' \bowtie \rho)$. Such a curve must be contained in a  component of $U_{\mu'}$ and must be disjoint from $\rho$. Since $\rho$ contains both a pants decomposition of  $U_{\mu'} -V$ and $\partial V$, the only way for that to happen is for $c \subseteq V$.  This implies $U_{\mu' \bowtie \rho} \subseteq V$ since  $U_{\mu' \bowtie \rho}$ is filled by the curves in $\lk(\mu' \bowtie \rho)$.		
	\end{subproof}
 As described before Claim \ref{claim:boundary_V_is_in_X},	Claim \ref{claim:U_rho=V} completes the proof of Proposition \ref{prop:comb_nesting_condition}.
\end{proof}

\begin{prop}\label{prop:C=C_0_condition}
	Let $\mu$ be a non-maximal simplex of $X_\alpha$  and let $x,y$ be curves of $X_\alpha$ contained in $\lk(\mu)$. If $x$ and $y$ are joined by a $W_\alpha$-edge but not an $X_\alpha$-edge in $X_\alpha^{+W_\alpha}$, then there exists a multicurve $\nu$ so that $\mu \subseteq \nu$ and $\nu \cup x$ and $\nu \cup y$ are adjacent vertices of $W_\alpha$.
\end{prop}

\begin{proof}
	If $x,y \in X_\alpha ^{+W_\alpha}$  are vertices that are joined be a $W_\alpha$-edge but not an $X_\alpha$-edge, then the curve $x$ and the curve $y$ fill a complexity 1 subsurface $V$ of $\Sdot$. Further, $\partial V \subseteq X_\alpha$ because the curves of $\partial V$ are part of the maximal simplices of $X_\alpha$ that contain $x$ and $y$ and are joined by an edge in $W_\alpha$.
	
	 Since $x,y \in \lk(\mu)$, every curve in $\mu$ is either a boundary curve of $V$ or contained in $\Sdot -V$.   By Lemma \ref{lem:Extending_To_Pants_decomposition}, we can extend $\mu \cap (\Sdot -V)$ to a pants decomposition $\rho$ of $\Sdot -V$ that contains only curves of $X_\alpha$. Thus, $\nu = \rho \cup \partial V$ is the desired multicurve.
\end{proof}

\subsection{Geometric conditions for $(X_\alpha,W_\alpha)$}\label{subsec:geom_conditions}
We  begin the verification of the geometric conditions of a combinatorial HHS by checking that $X_\alpha^{+W_\alpha}$ is uniformly quasi-isometric to $X_\alpha$. By Proposition \ref{prop:X_is_tree_adj}, this makes $X_\alpha^{W_\alpha}$ a quasi-tree.

\begin{cor}\label{cor:augment_X_quasi-tree}
	$X_\alpha^{+W_\alpha}$ is uniformly $E_\alpha$-equivariantly quasi-isometric to $X_\alpha$, and hence is a uniform quasi-tree.
\end{cor}

\begin{proof}
	Because $X_\alpha^{+W_\alpha}$ is a copy of $X_\alpha$ with additional edges, it suffices to verify that whenever  two vertices $x,y$ of $X_\alpha^{+W_\alpha}$ are joined by a $W_\alpha$-edge, then the distance between $x$ and $y$ in $X_\alpha$ is uniformly bounded. If $x$ and $y$ are joined by a $W_\alpha$-edge, then there exists a multicurve $\nu$ so that  $\nu \cup x$ and $\nu \cup y$ are adjacent edges of $W_\alpha$. As  the vertices of $W_\alpha$ are maximal simplices in $X_\alpha$, both $x$ and $y$ are joined by an $X_\alpha$-edge to each curve of $\nu$. Hence $d_{X_{\alpha}}(x,y) \leq 2$.
\end{proof}

For the remainder of this section, we analyze the links in $X_\alpha$ to verify the second geometric condition for $(X_\alpha,W_\alpha)$ to be a combinatorial HHS (Item \ref{CHHS:geom_link_condition} of Definition \ref{defn:CHHS}). We begin by recalling the notation we set out in Section \ref{sec:cHHS}. Note,  all links of simplices  in this section are taken in the the complex $X_\alpha$.

Recall, that if $Z$ is a subset of vertices of a graph $X$, then $X - Z$ is the subgraph spanned by the vertices in $X$ that are not in $Z$. For a simplex $\mu$ in $X_\alpha$ the space $Y_\mu$ is defined as the space $X_\alpha^{+W_\alpha} - \sat(\mu)$ and  $H(\mu)$ is the subset of $Y_\mu$ spanned by $\lk(\mu)$. We endow both $Y_\mu$ and $H(\mu)$ with their intrinsic path metrics and not the subspace metric from $X_\alpha$. Over the next two subsections, we will prove that $H(\mu)$ is uniformly hyperbolic and that the inclusion of $H(\mu)$ into $Y_\mu$ is uniformly a quasi-isometric embedding.

The simplices of $X_\alpha$ come in two types and we split our arguments accordingly.

\begin{defn}\label{defn:subsurface_and_tree_type}
	Let $\mu$ be a non-maximal simplex of $X_\alpha$ and $F_1,\dots, F_m$ be the fibers over $\alpha_1,\dots, \alpha_m$.
	
	\begin{enumerate}
		\item (Subsurface type) We say $\mu$ is of \emph{subsurface type} if $\mu$ contains a vertex of each fiber $F_i$. In this case, $\lk(\mu)$ is precisely $\C(\Sdot - \mu)$, because any curve on $\Sdot$ that is disjoint from $\mu$ will be disjoint from $\alpha$ in the image of $\Pi$. By Lemma \ref{lem:description_of_U}, $U_\mu$ will have no annular components, ensuring that $\C(\Sdot -\mu) = \C(U_\mu)$ in this case.
		\item (Tree guided type) We say $\mu$ is of \emph{tree guided type} if there exists a fiber $F_i$ so that $\mu$ does not contain a vertex of $F_i$. In this case, $\lk(\mu)$ is a proper subset of $\C(\Sdot - \mu)$.			
	\end{enumerate} 
\end{defn}

\subsubsection{Subsurface type links}
We handle the case where $\mu$ is of subsurface type first.  Because of the definition of $W_\alpha$ and the construction of $X_\alpha^{+W_\alpha}$,  Proposition \ref{prop:C=C_0_condition} implies that $H(\mu)$ is a copy of $\lk(\mu)$ with  additional edges between curves that intersect at most four times in a connected, complexity 1 subsurface of $\Sdot$.  Since $\lk(\mu) = \C(\Sdot -\mu)$, this makes $H(\mu)$ equal to $\C'(\Sdot - \mu)$, the modified curve graph from Definition \ref{defn:C'(S)}. We can therefore use the subsurface projection map to prove $H(\mu)$ is quasi-metrically embedded in $Y_\mu$.

\begin{prop}\label{prop:H_QI_to_C}
	If $\mu$ is a subsurface type simplex of $X_\alpha$, then   $H(\mu)=\C'(\Sdot -\mu)$. In particular, $H(\mu)$ is uniformly hyperbolic.
\end{prop}

\begin{proof}
	The vertices of $H(\mu)$ are all of the  curves of $\Sdot$ that are contained in $\Sdot - \mu$. By Proposition \ref{prop:C=C_0_condition}, there is an edge between two vertices of $H(\mu)$ if the curves are either disjoint or if they intersect at most four times  inside of a complexity 1 subsurface of $\Sdot -\mu$. However, this is precisely the definition of  $\C'(\Sdot -\mu)$.
\end{proof}

\begin{prop}\label{prop:subsurface_type_QI_embed}
	If $\mu$ is a subsurface type simplex of $X_\alpha$, then the inclusion of $H(\mu)$ into $Y_\mu$ is a uniform quasi-isometric embedding.
\end{prop}

\begin{proof}
	Because $Y_\mu = X_\alpha^{+W_\alpha} - \sat(\mu)$, every curve $y \in Y^0_\mu$ must intersect the subsurface $U_\mu \subseteq \Sdot -\mu$ by part 4 of Lemma \ref{lem:top_CHHS_dictionary}. Thus, subsurface projection gives a coarse  map $\pi_{\Sdot -\mu} \colon Y_\mu \to \C'(\Sdot -\mu) = H(\mu)$. This map is uniformly coarsely Lipschitz by Lemmas  \ref{lem:subsurface_projection_coarsely_Lipschtiz} and \ref{lem:checking_coarsely_lipschitz}. Since $\pi_{\Sdot -\mu}$ is also the identity on the vertices of $H(\mu)$, Lemma \ref{lem:coarsel_lipschitz} ensures that the inclusion of $H(\mu)$ into $Y_\mu$ is a uniform quasi-isometric embedding.  		
\end{proof}

\subsubsection{Tree guided type links}
We now handle the case when $\mu$ is not of subsurface type. In this case, $\lk(\mu)$ is not all of $\C(\Sdot - \mu)$ as there are curves on $\Sdot$ that are disjoint from $\mu$ but not in $X_\alpha$,  because their image under $\Pi$ intersects $\alpha$. This means we cannot rely on the hyperbolicity of $\C'(\Sdot-\mu)$ or use the subsurface projection map. Instead, we prove that  $X_\alpha -\sat(\mu)$ is $F_i$-guided for each fiber $F_i$ where $\mu \cap F_i = \emptyset$. This implies  $Y_\mu$ is a quasi-tree. We then show that $H(\mu)$ is quasi-isometrically embedded by using a combination of the subsurface projection maps and the tree guided structure of $X_\alpha - \sat(\mu)$.

\begin{prop}\label{prop:Y_mu_is_quasi-tree}
	Let $\mu$ be a tree guided type simplex of $X_\alpha$. For each $F_i$, if $\mu \cap F_i = \emptyset$, then $\sat(\mu) \cap F_i = \emptyset$ and $X_\alpha - \sat(\mu)$ is $F_i$-guided. Further, $Y_\mu$ is uniformly quasi-isometric to the uniform  quasi-tree $X_\alpha - \sat(\mu)$.
\end{prop}

\begin{proof}
	Let $F_i$ be a fiber so that $\mu \cap F_i = \emptyset$. We first show that $\sat(\mu) \cap F_i$ must also be empty. 
	Let $c$ be a curve in $\sat(\mu)$. By part 4 of Lemma \ref{lem:top_CHHS_dictionary}, $c$ is either a curve on $\Sdot-U_\mu$ or a curve of $\partial U_\mu$ that is not a core curve for an annular component of $U_\mu$. This means $\Pi(c)$ can not be contained in $\Pi(U_\mu)$.  
		However, we  have $\alpha_i \subseteq \Pi(U_\mu)$ because $\mu \cap F_i = \emptyset$. Thus, $c \not\in F_i$ when $c \in \sat(\mu)$. 
	
	Since $\sat(\mu) \cap F_i = \emptyset$, Lemma \ref{lem:T-guide_after_removal} says $X_\alpha - \sat(\mu)$ is also $F_i$-guided. Lemma \ref{lem:tree_adj_is_quasi-tree} then implies $X_\alpha - \sat(\mu)$ is a uniform quasi-tree. To see that $Y_\mu$ is uniformly quasi-isometric to $X_\alpha - \sat(\mu)$, we will show that for every $x,y\in Y^0_\mu$ that are not joined by an edge in $X_\alpha -\sat(\mu)$ but are joined by an edge in $Y_\mu$ the distance between $x$ and $y$ in $X_\alpha - \sat(\mu)$ is at most 2.
	
	 The definition of edges in $W_\alpha$ ensures that there exists a multicurve $\nu$ so that $\nu \cup x$ and $\nu \cup y$ are adjacent vertices of $W_\alpha$. Since $\Pi(\nu \cup x)$ is a pant decomposition of $S$ that contains $\alpha$, there exists $c \in \nu$ that is a curve in $F_i$. Since $\sat(\mu) \cap F_i =\emptyset$, the curve $c$ is a vertex of $X_\alpha-\sat(\mu)$. Hence $x$ and $y$ are 2 apart in $X_\alpha -\sat(\mu)$ as they are both disjoint from $c$. 
\end{proof}

\begin{prop}\label{prop:QI_embedding_tree_quided_link}
	Let $\mu$ be a tree guided type simplex of $X_\alpha$. The inclusion of $H(\mu)$  into the quasi-tree $Y_\mu$ is a $(15,0)$-quasi-isometric embedding, making $H(\mu)$ a  uniform quasi-tree. 
\end{prop}

\begin{proof}
	By the first item of Lemma \ref{lem:top_CHHS_dictionary}, $\lk(\mu)$ is a join whenever $U_\mu$ is disconnected. Since $\diam(H(\mu)) \leq 2$ whenever $\lk(\mu)$ is a join, we can assume $U_\mu$ is connected. Further, we can assume that $U_\mu$ is not a single annulus as $\lk(\mu)$ would be a single vertex in that case. Since $\mu$ is not of subsurface type, Lemma \ref{lem:description_of_U} says $U_\mu$ must contain the marked point $z$ as well as an infinite number of curve in $\Pi^{-1}(\alpha)$.

	We will construct a coarsely Lipschitz map $\psi_\mu \colon Y_\mu \to H(\mu)$.
	Since each curve in $Y_\mu^0$ must intersect $U_\mu$,  $\pi_{U_\mu}(y)$ is non-empty for each $y \in Y_\mu^0$. However, because $\mu \cap F_i = \emptyset$ it might be the case that none of the curve of $\pi_{U_\mu}(y)$ are  vertices of $X_\alpha$, that is,  the curves in  $\Pi(\pi_{U_\mu}(y))$ could intersect $\alpha$ on $S$.  This prevents us from using the subsurface projection map alone to construct $\psi_\mu$; instead we will use a combination of the subsurface projection map and the tree-guided structure of $X_\alpha - \sat(\mu)$.

	If the fiber $F_i$ has $\mu \cap F_i = \emptyset$, then Proposition \ref{prop:Y_mu_is_quasi-tree} implies that $\sat(\mu) \cap F_i = \emptyset$ and Lemma \ref{lem:links_are_subtrees_of_fiber} showed that $\lk(\mu) \cap F_i \subseteq Y_\mu$ is connected (in fact, it is the subtree of curves in $F_i$ that are contained in $U_\mu$). Because $\lk(\mu) \cap F_i$ is a connected subset of the tree $F_i$,  the closest point projection of $F_i$ onto $\lk(\mu) \cap F_i$ is well defined  and Lipschitz. We will use this map to build a coarsely Lipschitz map from $Y_\mu$ to $H(\mu)$. The first step is to project each $y \in  Y_\mu - H(\mu)$ with  $\lk(y) \cap \lk(\mu) \cap F_i = \emptyset$ onto $\lk(\mu) \cap F_i$.

	Let $y \in Y^0_\mu - H^0(\mu)$ and suppose $y$ is not disjoint from any curve in $\lk(\mu) \cap F_i$ where $\mu \cap F_i = \emptyset$, i.e., $\lk(y) \cap \lk(\mu) \cap F_i = \emptyset$.  Let $y' = y$ if $y \in F_i$ and $y'$ be any element of $\lk(y) \cap F_i$ if $y \not \in F_i$. By Lemma \ref{lem:links_are_subtrees_of_fiber}, $\lk(y) \cap F_i$ is a connected subset of $F_i$ in the latter case. Define $\psi_i(y)$ to be the closest point projection in $F_i$ of $y'$ onto $\lk(\mu) \cap F_i$. This is does not depend on the choice of $y'$ because $\lk(y) \cap \lk(\mu) \cap F_i = \emptyset$ means  $\lk(y) \cap F_i$ is a subtree of $F_i$ that is disjoint from the subtree $\lk(\mu) \cap F_i$ when $y \not \in F_i$.
	
	We now verify that for each $y \in Y^0_\mu - H^0(\mu)$ the set $\{ \psi_i(y) :\lk(y) \cap \lk(\mu) \cap F_i = \emptyset\}$ is a diameter 1 subset of $H(\mu)$ when it is non-empty.
	
	\begin{claim}\label{claim:psi_is_independent_of_fiber}
		Suppose  $y \in Y^0_\mu - H^0(\mu)$ so that $\lk(y) \cap \lk(\mu) \cap F_i = \emptyset$ and  $\lk(y) \cap \lk(\mu) \cap F_j = \emptyset$ where $i \neq j$ and $\mu \cap F_i =\emptyset$ and $\mu \cap F_j = \emptyset$. Then $\psi_i(y)$ is connected by an $X_\alpha$-edge to $\psi_j(y)$.
	\end{claim}
	
	\begin{subproof}
		Because $F_i$ and $F_j$ are both contained in $Y_\mu$ and because $\Pi(y)$ will be disjoint from or equal to each of $\alpha_i$ and $\alpha_j$ on $S$, there exists $a_i \in F_i$ and $a_j \in F_j$ so that $a_i$ and $a_j$ are disjoint and $i(a_i,y) =0$ and $i(a_j ,y) =0$; see Figure \ref{figure:Tree_a}. 
		By the definition of $\psi_i$ and $\psi_j$, we have $\psi_i(y) = \psi_i(a_i)$ and $\psi_j(y) = \psi_j(a_j)$. 
		
		\begin{figure}[!h]
			\centering
			\begin{subfigure}[b]{0.31\textwidth}
				\centering
				\def\svgscale{.6}
				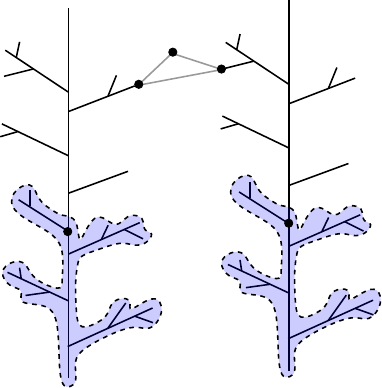
				\subcaption{}
				\label{figure:Tree_a}
			\end{subfigure}
			\quad
			\begin{subfigure}[b]{0.31\textwidth}
				\centering
				\def\svgscale{.6}
				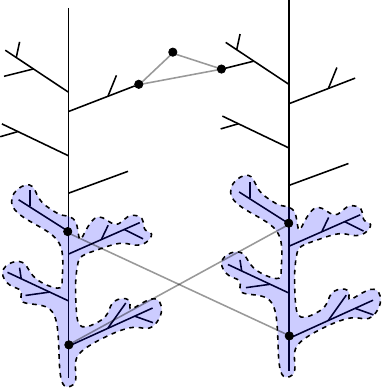
				\subcaption{}
				\label{figure:Tree_b}
			\end{subfigure}
			\quad 
			\begin{subfigure}[b]{0.31\textwidth}
				\centering
				\def\svgscale{.6}
				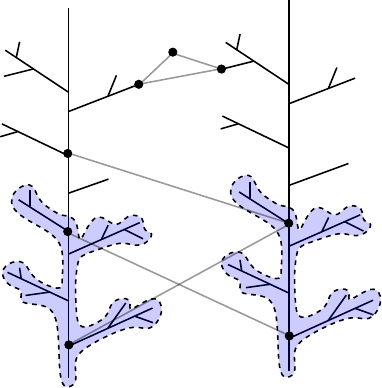
				\subcaption{}
				\label{figure:Tree_c}
			\end{subfigure}
			
			\caption{Schematic for the proof of Claim \ref{claim:psi_is_independent_of_fiber}}\label{figure:tree}
		\end{figure}

		Since $\mu \cap F_i$ and $\mu \cap F_j$ are both empty, there must exist a curve $b_j \in F_j$ so that $b_j$ is disjoint from both $\mu$ and $\psi_i(a_i)$; see Figure \ref{figure:Tree_b}. 
		Similarly, there must exists $b_i \in F_i$   so that $b_i$ is disjoint from both $\mu$ and $\psi_j(a_j)$. Let $\gamma_i$ be the geodesic in the tree $F_i$ that connects $a_i$ to $\psi_i(a_i)$.  
		Because $\psi_i(a_i)$ is connected by an $X_\alpha$-edge to $b_j$ and $X_\alpha-\sat(\mu)$ is $F_j$-guided, the union of the $F_j$-sequence (Definition \ref{defn:T_sequence}) for $\gamma_i$ is a connected subset of $F_j$ that contains both $a_j$ and $b_j$.
		Hence, the $F_j$-sequence of $\gamma_i$ contains $\psi_j(a_j)$. There is therefore some vertex $c_i$ of $\gamma_i$ that is joined by an $X_\alpha$-edge to $\psi_j(a_j)$; see Figure \ref{figure:Tree_c}. Since both $c_i$  and $b_i$ are in the subtree $\lk(\psi_j(a_j)) \cap F_i$, the $F_i$-geodesic connecting $c_i$  and $b_i$ is also contained in $\lk(\psi_j(a_j)) \cap F_i$ by Lemma \ref{lem:links_are_subtrees_of_fiber}. 
		However, because $b_i$ is in $\lk(\mu) \cap F_i$, the  $F_i$-geodesic connecting $c_i$  and $b_i$ must pass through $\psi_i(a_i)$. This implies $\psi_i(a_i) \in \lk(\psi_j(a_j)) \cap F_i$.
		\end{subproof}

	We now use the $\psi_i$'s and subsurface projection to define a coarse map $\psi_{\mu} \colon Y^0_\mu \to H^0(\mu)$: 
		\begin{itemize}
		\item For $y \in H^0(\mu)$, define $\psi_\mu(y) = y$. In this case, $y$ is a curve on $U_\mu$, so $\psi_\mu(y) = \pi_{U_\mu}(y)$.
		\item For $y \in Y^0_\mu - H^0(\mu)$ we have two subcases:
		\begin{enumerate}
			\item   For all the fibers $F_i$ where $\mu \cap F_i = \emptyset$, we have $\lk(y) \cap \lk(\mu) \cap F_i \neq \emptyset$. Thus, for each fiber with $\mu \cap F_i = \emptyset$, there is $ a_i \in \lk(\mu) \cap F_i$ that is disjoint from  $y$. Since $a_i \subseteq U_\mu$ for every such $a_i$,  $\pi_{U_\mu}(y)$ is disjoint from each  $ a_i$  as well.  Thus, $\pi_{U_\mu}(y)$ is a collection of curves in $X_\alpha$ that is contained in $\lk(\mu)$ and we  define $\psi_\mu(y) = \pi_{U_\mu}(y)$. This is a uniformly bounded subset by Lemma \ref{lem:subsurface_projection_coarsely_Lipschtiz}. 
			\item There exists  a fiber $F_i$ so that $\lk(y) \cap \lk(\mu) \cap F_i = \emptyset$, i.e., $y$ is not disjoint from any curve in $\lk(\mu) \cap F_i$. In this case, we define \[ \psi_\mu(y) =  \{ \psi_i(y) : \lk(y) \cap \lk(\mu) \cap F_i = \emptyset \}.\] By Claim \ref{claim:psi_is_independent_of_fiber}, $\psi_\mu(y)$ is a diameter 1 subset of $H(\mu)$.
		\end{enumerate}
	\end{itemize}
	This definition means we can break $\psi_\mu(y)$ into two mutually exclusive case: one where $\psi_\mu(y) = \pi_{U_\mu}(y)$ and one where $\lk(y) \cap \lk(\mu) \cap F_i = \emptyset$ for some fiber $F_i$. In both case, $\psi_\mu(y)$ is uniformly bounded by either Lemma \ref{lem:subsurface_projection_coarsely_Lipschtiz} or Claim \ref{claim:psi_is_independent_of_fiber}.
	
	Let $y_1,y_2$ be vertices of $Y_\mu$ that are connected by an edge of $Y_\mu$. We will show that the $H(\mu)$-distance between $\psi_\mu(y_1)$ and $\psi_\mu(y_2)$ is uniformly bounded. This implies $\psi_\mu$ is coarsely Lipschitz (Lemma \ref{lem:checking_coarsely_lipschitz}), which  finishes the proof by Lemma \ref{lem:coarsel_lipschitz}. We first handle the case when $y_1$ and $y_2$ are actually disjoint curves, i.e., joined by an $X_\alpha$-edge.

	\begin{claim}\label{claim:disjoint_case}
		If $y_1,y_2 \in Y_\mu^0$ are disjoint curves, then $d_{H(\mu)}(\psi_\mu(y_1),\psi_\mu(y_2)) \leq 2$.
	\end{claim}
	
	\begin{subproof}
		Let  $y_1,y_2 \in Y_\mu^0$ be disjoint curves.
		
	Assume first that $\psi_\mu(y_1) = \pi_{U_\mu}(y_1)$ and $\psi _\mu(y_2)= \pi_{U_\mu}(y_2)$. Without loss of generality, we can assume $y_1 \not \in H(\mu)$. If $y_2 \in H(\mu)$, then $y_2 = \psi_\mu(y_2) = \pi_{U_\mu}(y_2)$ does not intersect $\psi_\mu(y_1) = \pi_{U_\mu}(y_1)$ by Lemma \ref{lem:subsurface_projection_coarsely_Lipschtiz}. 
	If instead $y_2 \not \in H(\mu)$, then $ \lk(\mu) \cap F_i$ contains a curve of both $\lk(y_1) \cap F_i$ and $\lk(y_2) \cap F_i$ for each fiber with $F_i \cap \mu = \emptyset$. Further, $\lk(y_1) \cap \lk(y_2) \cap F_i$ is non-empty and connected because $y_1$ and $y_2$ are joined by an edge of $X_\alpha - \sat(\mu)$ and $X_\alpha -\sat(\mu)$ is $F_i$-guided. 
	Therefore, the Helly property of trees implies that  $\lk(y_1) \cap \lk(y_2)\cap \lk(\mu) \cap F_i \neq \emptyset$. 
	Thus, there is a curve $c \in X_\alpha - \sat(\mu)$ that is contained in $U_\mu$ and disjoint from both $y_1$ and $y_2$. This implies $d_{H(\mu)}(\psi_\mu(y_1),\psi_\mu(y_2)) \leq 2$ because both $\psi_\mu(y_1) = \pi_{U_\mu}(y_1)$ and $\psi _\mu(y_2)= \pi_{U_\mu}(y_2)$ will  not intersect $c$.

	Next assume $\lk(y_1) \cap \lk(\mu) \cap F_i = \emptyset$ and $\lk(y_2) \cap \lk(\mu) \cap F_i = \emptyset$ for  single fiber $F_i$. There must exist $a \in F_i$ so that $i(a,y_1)=0$ and $i(a,y_2) =0$. Thus, $d_{H(\mu)}(\psi_\mu(y_1),\psi_\mu(y_2)) =0$  because $\psi_i(y_1) = \psi_i(a) = \psi_i(y_2)$.
	
	Now assume $\lk(y_1) \cap \lk(\mu) \cap F_i = \emptyset$ and $\lk(y_2) \cap \lk(\mu) \cap F_j = \emptyset$ for distinct fibers $F_i$ and $F_j$ with $\mu \cap F_i$ and $\mu \cap F_j$ both empty. There exist $a_i \in F_i$ and $a_j \in F_j$ that are both disjoint from  or equal to each of $y_1$ and $y_2$ (however $a_i$ and $a_j$ might not be disjoint). 
	We can now employ a similar argument to  the proof of Claim \ref{claim:psi_is_independent_of_fiber}. There must be curves $b_i\in F_i$ and $b_j \in F_j$ that are disjoint from $\mu$ and respectively disjoint from $\psi_j(y_2) = \psi_j(a_j)$ and $\psi_i(y_1) = \psi_i(a_i)$. Let $\gamma_i$ be the concatenation of a $F_i$-geodesic connecting $\psi_i(a_i)$ to $a_i$ with the $X_\alpha$-edge from $a_i$ to $y_1$. 
	By examining the $F_j$-sequence for $\gamma_i$, we find a vertex $c_i$ on the $F_i$-geodesic between $a_i$ and $\psi_i(a_i)$ that is disjoint from $\psi_j(a_j)$. As in the end of Claim \ref{claim:psi_is_independent_of_fiber}, this forces $\psi_j(a_j) = \psi_j(y_2)$ to be joined by an $X_\alpha$-edge to $\psi_i(a_i) = \psi_i(y_1)$ because $\lk(\psi_j(a_j)) \cap F_i$ is connected. Since $\psi_i(y_1) \subseteq \psi_\mu(y_1)$ and $\psi_j(y_2) \subseteq \psi_\mu(y_2)$, this shows that the $H(\mu)$-distance between $\psi_\mu(y_1)$ and $\psi_\mu(y_2)$ is at most $1$.
	
	Finally, assume $\lk(y_1) \cap \lk(\mu) \cap F_i = \emptyset$, but $\psi_{\mu}(y_2)=\pi_{U_\mu}(y_2)$. In this case, $y_2$ must be disjoint from a curve in $\lk(\mu) \cap F_i$, i.e., $\lk(y_2) \cap \lk(\mu) \cap F_i \neq \emptyset$. In particular, $\psi_{\mu}(y_2) = \pi_{U_\mu}(y_2)$  does not intersect an element of $\lk(y_2) \cap \lk(\mu) \cap F_i$. Thus, we can prove $d_{H(\mu)}(\psi_\mu(y_1), \psi_\mu(y_2) ) \leq 1$ by proving that  $\psi_i(y_1) \subseteq \psi_\mu(y_1)$ is an element of $\lk(y_2) \cap \lk(\mu) \cap F_i$.
	
	Since $y_1$ and $y_2$ are disjoint, there exists $a \in F_i$ so that $i(a,y_1) =0$ and $i(a,y_2) =0$. Since $a \in \lk(y_2) \cap F_i$ and $\lk(y_2) \cap F_i$ is a subtree of $F_i$ (Lemma \ref{lem:links_are_subtrees_of_fiber}), the geodesic in $F_i$ from $a$ to any point in $\lk(y_2) \cap \lk(\mu) \cap F_i$ is contained in $\lk(y_2) \cap F_i$. By construction $\psi_i(y_1)$ is the endpoint of the unique geodesic in $F_i$ connecting $a$ to $\lk(\mu) \cap F_i$. Since this path is a subpath of any geodesic in $F_i$ that connects $a$ to a point in $\lk(y_2) \cap \lk(\mu) \cap F_i$, we have that $ \psi_i(y_1) \in \lk(y_2) \cap \lk(\mu) \cap F_i$ as desired.
\end{subproof}

	To complete the proof that $\psi_\mu$ is coarsely Lipschitz, suppose $y_1,y_2 \in Y_\mu^0$ are $Y_\mu$-adjacent but not $X_\alpha$-adjacent. There then exists some multicurve $\nu$ so that  $\nu \cup y_1$ and $\nu\cup y_2$ are adjacent pants decompositions in $W_\alpha$. Since $\nu$ is a single curve away from being a maximal simplex of $X_\alpha$, $\nu$ must contain a curve from each  fiber $F_i$. Because $\sat(\mu) \cap F_i =\emptyset$ for some fiber $F_i$, $Y_\mu = X_\alpha^{+W_\alpha} - \sat(\mu)$ contains all of this fiber $F_i$. Hence, the multicurve $\nu$ must contain a curve $c$ that is a vertex of $Y_\mu$. Thus there is $c \in Y_\mu$ with $c$ disjoint from both $y_1$ and $y_2$.  Therefore, Claim \ref{claim:disjoint_case} implies $$d_{H(\mu)}(\psi_\mu(y_1),\psi_\mu(y_2)) \leq d_{H(\mu)}(\psi_\mu(y_1),\psi_\mu(c))+ d_{H(\mu)}(\psi_\mu(c),\psi_\mu(y_2)) + 1 \leq 5.$$	Lemma \ref{lem:checking_coarsely_lipschitz} now implies $\psi_\mu$ is $(15,0)$-coarsely Lipschitz.
\end{proof}

\begin{remark}[$\psi_\mu$ versus closest point projection]\label{rem:cpp_vs_psi}
	Because $H(\mu)$ is quasi-isometrically embedded in $Y_\mu$ and $Y_\mu$ is hyperbolic  there is a coarse closest point projection $\mf{p}_\mu \colon Y_\mu \to H(\mu)$ defined by sending $y$ to the uniformly bounded set  $\{x \in H(\mu): d_{Y_\mu}(x,y) = d_{Y_\mu}(x,H(\mu)) \}$.
	This closest point projection is within uniformly bounded distance of the map $\psi_\mu$ defined in the proof of the above theorem.
	 To see this, let $y \in Y^0_\mu$. If $\psi_\mu(y) = \pi_{U_\mu}(y)$, then $y$ and $\psi_\mu(y)$ are at most 2 apart by Lemma \ref{lem:subsurface_projection_coarsely_Lipschtiz}. If instead $\lk(y) \cap \lk(\mu) \cap F_i =\emptyset$ for some fiber with $\mu \cap F_i = \emptyset$, then let $\gamma$ be the $Y_\mu$-geodesic connecting $y$ to $\mf{p}_\mu(y)$. Let $z$ be the endpoint of $\gamma$ in $\mf{p}_\mu(y)$. There must exist a curve $b\in F_i \cap \lk(\mu)$ so that either $b=z$ or $z$ and $b$ are disjoint. 
	 Thus, as in the proof of Claim \ref{claim:psi_is_independent_of_fiber}, examining the $F_i$-sequence of $\gamma$ produces vertex $c$ of $\gamma$ that is distance $1$ from $\psi_i(y)$. Since $\psi_i(y) \in H(\mu)$ and $z$ minimizes the distance from $y$ to $H(\mu)$, $c$ and $z$ must be at most 1 apart. Thus, $\psi_\mu(y)$ within distance 2 from $\mf{p}_\mu(y)$.
\end{remark}

\subsection{$(X_\alpha,W_\alpha)$ is a combinatorial HHS.}

We now combine the result of the previous two sections to conclude that $(X_\alpha,W_\alpha)$ is a combinatorial HHS with a cobounded action by $E_\alpha$

\begin{thm}
	The graph $W_\alpha$ is connected and a hierarchically hyperbolic space because the  pair $(X_\alpha,W_\alpha)$ is a combinatorial HHS. Further, the action of $E_\alpha$ on $\Sdot$ induces a cobounded action of $E_\alpha$ on $W_\alpha$.
\end{thm}

\begin{proof}
	We  first check the five parts of Definition \ref{defn:CHHS} for some $\delta \geq 2$ ultimately depending only on $S$. Theorem \ref{thm:CHHS_are_HHS} will then imply that $W$ is connected and a hierarchically hyperbolic space.
	\begin{enumerate}
		\item If $\mu_1,\dots, \mu_n$ is a sequence of non-maximal simplices of $X_\alpha$ so that $$\lk(\mu_1) \subsetneq \lk(\mu_2) \subsetneq \dots \subsetneq \lk(\mu_n),$$ then Lemma \ref{lem:top_CHHS_dictionary} implies $n$ is bounded in terms of the complexity of $S$.
		\item Corollary \ref{cor:augment_X_quasi-tree} proves that $X_\alpha^{+W_\alpha}$ is quasi-isometric to a tree and hence hyperbolic.
		\item The links of simplices of $X_\alpha$ split into subsurface and tree-guided types; see Definition \ref{defn:subsurface_and_tree_type}. The  proof of uniform hyperbolicity and  uniform quasi-isometric embedding for subsurface type is shown in Propositions \ref{prop:H_QI_to_C} and  \ref{prop:subsurface_type_QI_embed}, while the proof for the tree-guided type is completed in Proposition \ref{prop:QI_embedding_tree_quided_link}. 
		\item Because $\diam(H(\mu)) >2$ implies $\lk(\mu)$ is not a join or a single vertex, this follows from Proposition \ref{prop:comb_nesting_condition}.
		\item This condition is ensured by Proposition \ref{prop:C=C_0_condition}.
	\end{enumerate}
	
	To see that the action of $E_\alpha$ on $W_\alpha$ is cobounded, start with the fact that $\Pi$ will map each vertex of $W_\alpha$ to a pants decomposition of $S$ containing the curves in $\alpha$. Recall, $E_\alpha$ fits into the short exact sequence \[1 \to \Push(S,z) \to E_\alpha \to G_\alpha \to 1\] and $\Pi$ is equivariant with respect to the quotient map $E_\alpha \to G_\alpha$. For each pants decomposition $\rho$ of $S$ with $\alpha \subseteq \rho$, there only exists finitely many $\Push(S;z)$-orbits of elements of $\{\mu \in W^0_\alpha : \Pi(\mu) = \rho\}$. Moreover,  there exist only finitely many $G_\alpha = \stab_{\mcg(S)}(\alpha)$-orbits of pants decompositions of $S-\alpha$. Hence,  there are only finitely many $E_\alpha$-orbits of vertices of $W_\alpha$, making the $E_\alpha$-action cobounded. 
\end{proof}

\section{Adding in annuli}\label{sec:annuli}

We continue with the assumptions of Section \ref{sec:non-annular_case}. That is, $S \cong S_{g,0}^0$ with $g \geq 2$, and $\alpha = \alpha_1 \cup \dots \cup \alpha_m$ is a fixed multicurve on $S$. For a fixed $z \in S$,  the group $E_\alpha$ is then the full preimage in $\mcg(S;z)$ of the stabilizer of $\alpha$ in $\mcg(S)$. The spaces $X_\alpha$ and $W_\alpha$ are as defined in Definitions \ref{defn:X} and \ref{defn:W}. 

The extension group $E_\alpha$ does not act metrically properly on the space $W_\alpha$. This is because the action of $E_\alpha$ on a maximal simplex in $X_\alpha$ has large stabilizer, namely the Dehn twists around the curves in the simplex. To fix this, we need to ``blow-up'' $X_\alpha$ by adding  marking arcs for each of the curves  in $X_\alpha$  that account for the action of Dehn twists. We denote this blow up $B(X_\alpha)$. The corresponding $B(X_\alpha)$-graph will be denoted $B(W_\alpha)$. As in the previous section, we will as a quantity is \emph{uniform} if it does not depend on the multicurve $\alpha$, the surface $S$, or specific simplices/vertices of $B(X_\alpha)$ and $B(W_\alpha)$.

\subsection{Defining the blow-ups of $X_\alpha$ and $W_\alpha$}

For a curve $c$ in $\Sdot$, let $A_c$ be the annulus with core curve $c$. For each $c \in X_\alpha^0$, let $B(c)$ denote the set of vertices of the annular complex $\C(A_c)$. We call the elements of $B(c)$ the \emph{marking arcs} for $c$.	Let $B(X_\alpha)$ be the space obtained from $X_\alpha$ by adding  in the following vertices and edges:
\begin{enumerate}
	\item For each curve $c \in X_\alpha$, add the elements of $B(c)$ as additional vertices.
	\item For each curve $c \in X_\alpha$, add an edge between $c$ and each element of $B(c)$.
	\item If $c$ and $c'$ are disjoint curves in $X_\alpha$, then for each $v \in B(c)$ and $v' \in B(c') \cup \{c'\}$  add an edge connecting $v$ and $v'$.
\end{enumerate}

\noindent For a simplex $\mu$  of $B(X_\alpha)$ we define three sets of vertices of $B(X_\alpha)$:
\begin{itemize}
	\item $\base(\mu)$ is the set  $\{ c \in X^0_\alpha : c \in \mu \}$.
	\item $\supp(\mu)$ is the set $\base(\mu) \cup \{c \in X^0_\alpha : B(c) \cap \mu \neq \emptyset\}$.
	\item $\markcurve(\mu)$ is the set of vertices of $\mu -\base(\mu)$.
\end{itemize}
We call the curves in $\base(\mu)$ and $\supp(\mu)$ the \emph{base} and \emph{support} curves of $\mu$ respectively. The vertices of $\markcurve(\mu)$ are called the \emph{marking arcs} of $\mu$. Note, each curve of $\supp(\mu) - \base(\mu)$ corresponds to a marking arc of $\mu$ that does not have a base curve in $\mu$.  We say that $c \in \base(\mu)$ is  an \emph{unmarked base curve of $\mu$} if $B(c) \cap \markcurve(\mu) = \emptyset$. We denote the set of unmarked base curves of $\mu$ by $\unmark(\mu)$. The maximal simplices of $B(X_\alpha)$ are precisely the (not necessarily clean) markings of $\Sdot$ whose base curves are maximal simplices in $X_\alpha$.

We now define our $B(X_\alpha)$-graph which we denote $B(W_\alpha)$. Recall, for $t,t' \in B(c) = \C^0(A_c)$,  $d_{A_c}(t,t')$ denotes the distance between $t$ and $t'$ in the annular complex $\C(A_c)$.  The graph $B(W_\alpha)$  has vertex set the maximal simplices of $B(X_\alpha)$ with edges between two vertices $\mu,\nu \in B(W_\alpha)$ corresponding to \emph{twist} and \emph{flip} moves:
\begin{itemize}
	\item (twist) $\base(\mu) = \base(\nu)$ and there exists exactly one curve  $c \in \base(\mu)$ so that $\markcurve(\mu) \cap B(c) \neq \markcurve(\nu) \cap B(c)$. Further, for this curve $c$, if $t \in \markcurve(\mu) \cap B(c)$ and $t'\in \markcurve(\nu) \cap B(c)$,  then $d_{A_c}(t,t')=1$.
	
	\item (flip) There exists curves $c \in \base(\mu)$ and $d \in \base(\nu)$ so that
	\begin{itemize}
		\item $(\base(\mu) - c) \cup d = \base(\nu)$ and $i(c,d) \leq 4$, i.e.,  $\base(\mu)$ and $\base(\nu)$ differ by an edge in $W_\alpha$;
		\item if $t_c \in \markcurve(\mu)$ and $t_d \in \markcurve(\nu)$ are the marking arcs for $c$ and $d$, then \[d_{A_c}(t_c,d) \leq 1 \text{ and } d_{A_d}(t_d,c) \leq 1;\]
		\item for all curves $b \in \base(\mu) -c = \base(\nu) - d$, the marking arc of $b$ in $\mu$ is the same as the marking arc of $b$ in $\nu$.
	\end{itemize}
\end{itemize}

The graphs $B(X_\alpha)$ and $B(W_\alpha)$ are related to $X_\alpha$ and $W_\alpha$ by simplicial maps induced by sending every simplex $\mu \subseteq B(X_\alpha)$ to $\supp(\mu)$. This will allow us to reduce properties of $(B(X_\alpha),B(W_\alpha))$ to properties about $(X_\alpha,W_\alpha)$.

\begin{lem}\label{lem:forget_marking}
	The following  maps on vertices extend to simplicial maps between the given complexes.
	\begin{itemize}
		\item $B^0(X_\alpha) \to X_\alpha$	by $x \to \supp(x)$.
		\item $B^0(X_\alpha)^{+B(W_\alpha)} \to X_\alpha^{+W_\alpha}$ by $x \to \supp(x)$.
		\item $B^0(W_\alpha) \to W_\alpha$ by $\mu \to \supp(\mu)$.
	\end{itemize}	
\end{lem}

\begin{proof}
	In all three cases, it suffices to show that two vertices joined by an edge in the domain graph have supports that are either equal or joined by an edge in the codomain graph.
	
	Let $x,y$ be vertices of $B(X_\alpha)$ that are joined by an edge. Thus either $\supp(x) = \supp(y)$ or $\supp(x)$ is disjoint from $\supp(y)$. Hence $\supp(x)$ and $\supp(y)$ are either equal or joined by an edge in $X_\alpha$.
	
	Let $x,y$ be vertices of $B(X_\alpha)^{+B(W_\alpha)}$ that are joined by an edge. By the previous paragraph, we can assume $x$ and $y$ are joined by a $B(W_\alpha)$-edge. If $x$ and $y$ are joined by an edge coming from a twist move in $B(W_\alpha)$, then $\supp(x) = \supp(y)$. If instead $x$ and $y$ are joined be an edge coming from a flip move, then $\supp(x)$ and $\supp(y)$ fill a connected complexity 1 subsurface of $\Sdot$ and intersect at most 4 times. Hence, $\supp(x)$ and $\supp(y)$ will be joined by a $W_\alpha$-edge of $X_\alpha ^{W_\alpha}$.
		
	Let $\mu,\nu$ be vertices of $B(W_\alpha)$ that are joined by an edge. If $\mu$ and $\nu$ are joined by a twist move, then $\supp(\mu) = \supp(\nu)$. If $\mu$ and $\nu$ are joined by a flip move, then $\supp(\mu)$ and $\supp(\nu)$ are joined by a flip move in $W_\alpha$	by definition.
\end{proof}

\subsection{The combinatorial conditions with annuli}
We now verify that $(B(X_\alpha),B(W_\alpha))$ satisfies the combinatorial condition of being a combinatorial HHS. We begin by describing the links of simplices of $B(X_\alpha)$. Henceforth, we use $\blk(\cdot)$ to denote the link of simplices in $B(X_\alpha)$ and $\xlk(\cdot)$ to denote the link of simplices in $X_\alpha$.

\begin{lem}\label{lem:B-links}
	For any non-maximal simplex $\mu$ of $B(X_\alpha)$, the link of $\mu$ in $B(X_\alpha)$ is the subgraph spanned by the union of the following vertices:
	\begin{enumerate}
		\item The vertices of $B(c)$ for each $c \in \unmark(\mu)$.
		\item The set of vertices $\{x \in B^0(X_\alpha) : \supp(x) \in \xlk(\supp(\mu))\}$.
		\item The vertices of $\supp(\mu) - \base(\mu)$. 
	\end{enumerate}
	Moreover, if $\supp(\mu) \neq \base(\mu)$, then $\blk(\mu)$ is either a single vertex (which will be $\supp(\mu) - \base(\mu)$) or  is a join of some graph with $\supp(\mu) -\base(\mu)$. 
\end{lem}

\begin{proof}
	
	If $c \in \unmark(\mu)$, then  all the vertices of $B(c)$ are in $\blk(\mu)$ by the definition of the edges in $B(X_\alpha)$.
	
	If $x \in B^0(X_\alpha)$ with $\supp(x) \in \xlk(\supp(\mu))$, then $\supp(x)$ is disjoint from  all the curves in $\supp(\mu)$. By the edge relations in $B(X_\alpha)$, this makes $x \in \blk(\mu)$.
	
	Finally, suppose $c$ is a vertex of  $\supp(\mu) - \base(\mu)$. Since $\mu$ is a simplex,  $c$ must be disjoint from $\base(\mu)$. However, this means $c$ must be joined by an edge of every element of $\mu$ by the definition of the $B(X_\alpha)$ edges. 
	
	To verify these are the only elements of $\blk(\mu)$ let $x$ be a vertex of $\blk(\mu)$. If $x \in X_\alpha^0$, then either $x \in \xlk(\supp(\mu))$ or $x \in \supp(\mu) -\base(\mu)$. If instead $x$ is in $B(c)$ for some curve $c \in X^0_\alpha$, then either $c \in \unmark(\mu)$	 or $c \in \xlk(\supp(\mu))$.
	
	For the moreover statement, let $c \in \supp(\mu) - \base(\mu)$ and let  $t_c$ be the element of $\markcurve(\mu) \cap B(c)$. Any vertex $x \in \blk(\mu)$ will be joined by an edge to $t_c$ by definition. However, this makes $x$ and $c$ joined by an edge from the  the definition of  edges in $B(X_\alpha)$.
\end{proof}

We now describe the nesting of links of $B(X_\alpha)$ by using the topology of subsurfaces of $\Sdot$. This is an extension of Lemma \ref{lem:top_CHHS_dictionary} from $X_\alpha$  to $B(X_\alpha)$ and is our main tool in this section. Recall that Definition \ref{defn:U_mu} defined $U_{\supp(\mu)}$ as subsurface filled by the curves in $\xlk(\supp(\mu))$.  We remark that it  is possible that $U_{\supp(\mu)}$ is empty for a non-maximal simplex $\mu \subseteq B(X_\alpha)$. This occurs when $\supp(\mu)$ is a pants decomposition of $\Sdot$, but not every base curve of $\mu$ has a marking arc.

\begin{lem}\label{lem:B-nesting}
	Suppose $\mu$ and $\mu'$ are non-maximal simplices of $B(X_\alpha)$ with $\supp(\mu) = \base(\mu)$ and $\supp(\mu') = \base(\mu')$.
	We have $\blk(\mu) \subseteq \blk(\mu')$ if and only if
	\begin{itemize}
		\item $U_{\supp(\mu)} \subseteq U_{\supp(\mu')}$; and
		\item  for each $c \in \unmark(\mu)$, either $c \in \unmark(\mu')$ or $c \subseteq U_{\supp(\mu')}$.
	\end{itemize}
\end{lem}

\begin{remark}\label{rem:nested_subsurfaces}
	The above conditions for the nesting of links can be stated purely in terms of subsurfaces by using annuli. For a simplex $\mu$ where $\supp(\mu) =\base(\mu)$ define $C_\mu$ to be the disjoint union of the annuli $\{A_c: c \in \unmark(\mu)\}$. The conclusion of Lemma \ref{lem:B-nesting} can then be phrased as $\blk(\mu) \subseteq \blk(\mu')$ if and only if $U_{\supp(\mu)} \subseteq U_{\supp(\mu')}$ and $C_\mu \subseteq U_{\supp(\mu')}  \cup C_{\mu'}$.
\end{remark}

\begin{proof}[Proof of Lemma \ref{lem:B-nesting}]
	Before starting the proof, observe that  by applying  Lemma \ref{lem:B-links}, the assumption that $\base(\mu) = \supp(\mu)$  implies that a vertex $x \in \blk(\mu)$  if  and only if  $\supp(x) \in \xlk(\supp(\mu))$ or $x \in B(c)$ for some  $c \in \unmark(\mu)$.
	
	Suppose   $\blk(\mu) \subseteq  \blk(\mu')$. For the first bullet, we prove $\xlk(\supp(\mu)) \subseteq \xlk(\supp(\mu'))$, which implies $U_{\supp(\mu)} \subseteq U_{\supp(\mu')}$ by  Lemma \ref{lem:top_CHHS_dictionary}. If $x \in \xlk(\supp(\mu))$, then $x \in \blk(\mu)$ and $x\in X_\alpha$. Since $ \blk(\mu) \subseteq \blk(\mu')$, $x$ must be an element of $\blk(\mu')$ that is also an element of $X_\alpha$. By the above observation, this only occurs if $x = \supp(x) \in \xlk(\supp(\mu'))$ as well.   This shows $\xlk(\supp(\mu)) \subseteq \xlk(\supp(\mu'))$ as desired.
	
	For the second bullet, let $c$ be an unmarked base curve of $\mu$.  By Lemma \ref{lem:B-links}, this implies $B(c) \subseteq \blk(\mu)$, which means $B(c)$ is also contained in $\blk(\mu')$. Thus either $c \in \unmark(\mu')$ or $c \in \xlk(\supp(\mu'))$. In the latter cases, we have that $c \subseteq U_{\supp(\mu')}$ by Definition \ref{defn:U_mu}.
	
	Now suppose  $U_{\supp(\mu)} \subseteq U_{\supp(\mu')}$ and  for all $c \in \unmark(\mu)$, either $c \subseteq U_{\supp(\mu')}$ or $c \in\unmark(\mu')$. Let $x \in \blk(\mu)$. By the observation at the beginning, we know that either $\supp(x) \in \xlk(\supp(\mu))$ or $x$ is a marking arc for some unmarked base curve $c$ of $\mu$. In the first case, Lemma \ref{lem:top_CHHS_dictionary} implies $\supp(x) \in \xlk(\supp(\mu'))$ because $U_{\supp(\mu)} \subseteq U_{\supp(\mu')}$. Thus, $x\in\blk(\mu')$ by Lemma \ref{lem:B-links}. In the second case, we either have $c \in \unmark(\mu')$ or $c \subseteq U_{\supp(\mu')}$  by hypothesis. In both cases $x \in \blk(\mu')$ by Lemma \ref{lem:B-links}. Thus we have $\blk(\mu) \subseteq \blk(\mu')$
\end{proof}

Using Lemma \ref{lem:B-nesting}, we show that $(B(X_\alpha),B(W_\alpha))$ satisfies the combinatorial conditions of being a combinatorial HHS. The essence of these arguments are the same as Proposition \ref{prop:comb_nesting_condition} and Proposition \ref{prop:C=C_0_condition} in the non-annular case.

\begin{prop}\label{prop:annular_nesting_condition}		
	Suppose $\mu$ and $\mu'$ are non-maximal simplices of $B(X_\alpha)$ so that there exists a non-maximal simplex $\nu$ of $X_\alpha$ with $\blk(\nu) \subseteq \blk(\mu) \cap \blk(\mu')$ and $ \blk(\nu)$ not a join or a single vertex. There exists a (possibly empty) simplex $\rho$ of $\blk(\mu')$ so that $\blk(\mu' \bowtie \rho) \subseteq \blk(\mu)$ and if $\nu$ is any simplex as in the proceeding sentence, then $\blk(\nu) \subseteq \blk(\mu' \bowtie \rho)$.
\end{prop}

\begin{proof}
		Let $\eta = \markcurve(\mu) \bowtie \supp(\mu)$ and $\eta' = \markcurve(\mu') \bowtie \supp(\mu')$. Lemma \ref{lem:B-links} implies that $\blk(\eta)$ is  $\blk(\mu)$ minus the vertices of $\supp(\mu) - \base(\mu)$ and that the same holds for $\eta'$ and $\mu'$. Thus, $\blk(\mu) = \blk(\eta) \bowtie (\supp(\mu) - \base(\mu))$ and the same holds for $\eta'$ and $\mu'$. We first verify that it suffices to prove the proposition is true for $\eta$ and $\eta'$, instead of $\mu$ and $\mu'$. The advantage of working with $\eta$ and $\eta'$ over $\mu$ and $\mu'$ is that Lemma  \ref{lem:B-nesting} applies to $\eta$ and $\eta'$ because $\supp(\eta) = \base(\eta)$ and $\supp(\eta') = \base(\eta')$. 
	
	\begin{claim}
		If $\rho$ is a simplex that satisfies Proposition \ref{prop:annular_nesting_condition}  for $\eta$ and $\eta'$, then  the join of $\rho$ and $\supp(\mu') - \base(\mu')$ is a simplex that satisfies the proposition for $\mu$ and $\mu'$.
	\end{claim}
	
	\begin{subproof}	
		Let $\rho$ be the simplex that satisfies Proposition \ref{prop:annular_nesting_condition} for $\eta$ and $\eta'$ and let $\rho' = \rho \bowtie (\supp(\mu')) -\base(\mu'))$. Since $\rho \subseteq \blk(\eta')$ and $ \blk(\mu') = \blk(\eta') \bowtie (\supp(\mu')) -\base(\mu'))$, $\rho'$ is in fact a simplex of $B(X_\alpha)$. Moreover, $\eta' \cup \rho = \mu' \cup \rho'$, hence $\blk(\eta' \bowtie \rho) = \blk(\mu' \bowtie \rho')$. As $\rho$ satisfies Proposition \ref{prop:annular_nesting_condition} for $\eta$ and $\eta'$, we have $\blk(\eta' \bowtie \rho) = \blk(\mu' \bowtie \rho') \subseteq \blk(\eta)$. Since $\blk(\eta) \subset \blk(\mu)$, this gives $\blk(\mu' \bowtie \rho') \subseteq \blk(\mu)$ as desired.
			
		We now show the second requirement of Proposition \ref{prop:annular_nesting_condition} for $\rho'$ and $\mu',\mu$. If  there exists a non-maximal simplex $\nu \subseteq B(X_\alpha)$, so that $\blk(\nu) \subseteq \blk(\mu)$   and $\blk(\nu) \not \subseteq \blk(\eta)$, then $\supp(\mu) - \base(\mu) \neq \emptyset$ and $\blk(\nu)$ must contain vertices of $\supp(\mu) - \base(\mu)$. In this case, Lemma \ref{lem:B-links} says $\blk(\mu)$  is either a single vertex contained in $\supp(\mu) - \base(\mu)$ or the join of some graph with the vertices in $\supp(\mu) - \base(\mu)$. Hence, $\blk(\nu)$ would need to be either a single vertex or a join as well. Since the same reasoning applies to $\blk(\nu) \subseteq \blk(\mu')$, we have that whenever $\blk(\nu) \subseteq \blk(\mu) \cap \blk(\mu')$ and  $\blk(\nu)$ is not a join or a single vertex, then $\blk(\nu)$ must be contained in $\blk(\eta) \cap \blk(\eta')$. 
		Since $\rho$ satisfies Proposition \ref{prop:annular_nesting_condition} for $\eta$ and $\eta'$, then any such simplex $\nu$ will have $\blk(\nu) \subseteq \blk(\eta' \bowtie \rho) =\blk(\mu' \bowtie \rho')$.  
	\end{subproof}
	
	We now construct the required simplex $\rho$ for $\eta$ and $\eta'$.  We will define two simplices, $\rho_U$ and $\rho_A$, of $B(X_\alpha)$ so that the desired $\rho$  will be $\rho_U \bowtie \rho_A$. 
	
		To define $\rho_A$, let $\mc{A}$ to be the set of  curves $c \in \unmark(\eta')$ so that  $c \not \in \unmark(\eta)$  and $c \not\subseteq U_{\supp(\eta)}$. If $\mc{A} = \emptyset$, then define $\rho_A = \emptyset$. Otherwise, for each $c \in \mc{A}$, pick any vertex $t_c \in B(c)$, then define $\rho_A = \{ t_c : c \in   \mc{A}\}$.
		
		The definition of $\rho_U$ will be similar to the definition of $\rho$ in the proof of Proposition \ref{prop:comb_nesting_condition} for the non-annular case, with extra care to keep track of the set $B(c)$ that are contained in $\blk(\eta) \cap \blk(\eta')$. Let $V_0$ be the subsurface of $\Sdot$ filled by the curves in $\xlk(\supp(\eta)) \cap \xlk(\supp(\eta'))$, then let $V$ be the disjoint union of the non-annular components of $V_0$. Let $\sigma$ be the set of curves in $\unmark(\eta)$ that are contained in $U_{\supp(\eta')}$. 
		Since $\sigma$ contains base curves of $\eta$, none of the curves of $\sigma$ can be contained in $V_0$, thus $\sigma \subseteq U_{\supp(\eta')} - V_0$.  The proof of Claim \ref{claim:boundary_V_is_in_X} implies $\partial V$ is a simplex of $X_\alpha$, hence Lemma \ref{lem:Extending_To_Pants_decomposition} produces a pants decomposition $\tau$ of $U_{\supp(\eta')} - V_0$ so that $\sigma \subseteq \tau$ and $\tau \subseteq X_\alpha$. For each $c \in \tau \cup \partial V_0$, if $c \not \in \unmark(\eta)$ and $c \not \subseteq U_{\supp(\eta)}$, select a marking arc $t_c \in B(c)$. Define $$\rho_U = \tau \cup \partial V_0 \cup \{ t_c : c \in  \tau \cup \partial V_0 \text{ with } c \not \in \unmark(\eta) \text{ and } c\not\subseteq U_{\supp(\eta)} \}.$$
		
		We now verify that $\rho = \rho_A \bowtie \rho_U$ has the three desired properties given in the proposition. Note, $\supp(\eta' \bowtie \rho) = \base(\eta' \bowtie \rho) = \base(\eta') \cup \base(\rho_U)$ by construction.
		\begin{itemize}
			\item  $\rho \subseteq \blk(\eta')$: This is immediate from the construction of $\rho$, as every vertex  $x \in \rho$ is either a marking arc for an unmarked base curve of $\eta'$ or has $\supp(x) \in \xlk(\supp(\eta'))$.
			
			\item $\blk(\eta' \bowtie \rho) \subseteq \blk(\eta)$: Since $\base(\rho_U) = \tau \cup \partial V$, we have $U_{\supp(\eta' \bowtie \rho)} = U_{\supp(\eta') \bowtie (\tau \cup \partial V)}$. Thus, a repeat of the proof of Claim \ref{claim:U_rho=V} proves that $U_{\supp(\eta' \bowtie \rho)} = U_{\supp(\eta') \bowtie (\tau \cup \partial V)} = V$. Since $V$ is filled by a subset of curves in $\xlk(\supp(\eta))$, we have  $U_{\supp(\eta' \bowtie \rho)} = V \subseteq U_{\supp(\eta)}$. This verifies the first condition in Lemma \ref{lem:B-nesting}.
			 
			 For the second condition in Lemma \ref{lem:B-nesting}, let $c$ be an unmarked base curve of $\eta' \bowtie \rho$. If $c \in \base(\eta')$, then the choice of $\rho_A$ ensures that $c$  is either an unmarked base curve of $\eta$ or is contained on $U_{\supp(\eta)}$. If instead $c \in \base(\rho_U)$, then the choice of marking arcs in $\rho_U$ again ensures that $c$ is  either an unmarked base curve of $\eta$ or is contained in $U_{\supp(\eta)}$. Hence, $\blk(\eta' \bowtie \rho) \subseteq \blk(\eta)$.
			 
			\item If $\nu$ is a non-maximal simplex of $B(X_\alpha)$ with $\blk(\nu) \subseteq \blk(\eta) \cap \blk(\eta')$, then either $\blk(\nu) \subseteq \blk(\eta' \bowtie \rho)$ or $\blk(\nu)$ is a single vertex or a join: 
			Assume $\blk(\nu)$ is not a single vertex or a join  and $\blk(\nu) \subseteq \blk(\eta) \cap \blk(\eta')$. Because $\supp(\nu) \neq \base(\nu)$ would imply that $\blk(\nu)$ is a join or a single vertex by Lemma \ref{lem:B-links}, we know $\supp(\nu) = \base(\nu)$  and can apply Lemma \ref{lem:B-nesting} to  $\nu$.

			If $\nu$ contains an unmarked base curve $c$, then every vertex in $\blk(\nu) - B(c)$ will be joined by an edge to every vertex of $B(c)$. Since $\blk(\nu)$ is not a join, this means either $\blk(\nu) = B(c)$ or $\nu$ does not have any unmarked base curves.  
			
			If $\blk(\nu) = B(c)$, then $U_{\supp(\nu)} = \emptyset$ and $c$ is the only unmarked base curve of $\nu$.  Since $\blk(\nu) \subseteq \blk(\eta) \cap \blk(\eta')$,  we have four possibilities coming from Lemma \ref{lem:B-nesting}. In each case we verify that $c$ is either an unmarked base curve of $\eta' \bowtie \rho$ or is contained in $U_{\supp(\eta' \bowtie \rho)} = V$.
			\begin{itemize}
				\item If $c \in \unmark(\eta) \cap \unmark(\eta')$ or  $c \subseteq U_{\supp(\eta)}$ and $c \in \unmark(\eta')$, then $c$ is an unmarked base curve of $\eta' \bowtie \rho$, because of the choice of the marking arcs in $\rho_A$.
				\item If $c \in \unmark(\eta)$ and $c \subseteq U_{\supp(\eta')}$, then $c$ is a curve in the multicurve $\sigma$ defined while constructing $\rho_U$. Hence $c \in \base(\rho_U)$ and $c$ is unmarked by the choice of marking arcs in $\rho_U$.
				\item  	If $c \subseteq U_{\supp(\eta)}$ and $c\subseteq U_{\supp(\eta')}$, then $c \subseteq V_0$ because $c \in \xlk(\supp(\eta)) \cap \xlk(\supp(\eta'))$. If $c$ is a core curve of an annular component of $V_0$, then $c \in \partial V_0$. Thus $c \in \base(\rho_U)$ and $c$ will be unmarked in $\eta' \bowtie \rho$ because of the choice of marking arcs in $\rho_U$. If instead $c \subseteq V$, then $c \subseteq U_{\eta' \bowtie \rho}$.
			\end{itemize}
			 Since  $U_{\supp(\nu)} = \emptyset$, this verifies  $\blk(\nu) \subseteq \blk(\eta' \bowtie \rho)$ by  Lemma \ref{lem:B-nesting}.

			Now assume $\nu$ has no unmarked base curves. Since $\nu$ is non-maximal, this means $U_{\supp(\nu)}$ must be non-empty. Since $\blk(\nu) \subseteq \blk(\eta) \cap \blk(\eta')$, we must have $U_{\supp(\nu)} \subseteq V_0$. 
			If $U_{\supp(\nu)}$ contains an annular component of $V_0$, then every vertex of $\blk(\nu)$ would have to be equal to or joined by an edge to the core curve of that component. However, this cannot be the case since $\blk(\nu)$ is not a join. Thus, $U_{\supp(\nu)} \subseteq V =  U_{\supp(\eta' \bowtie \rho)}$. Since $\nu$ has no unmarked base curves, Lemma \ref{lem:B-nesting} implies $\blk(\nu) \subseteq \blk(\eta' \bowtie \rho)$.\qedhere 	
		\end{itemize}
\end{proof}

\begin{prop}\label{prop:annular_C=C0}
	Let $\mu$ be a non-maximal simplex of $B(X_\alpha)$ and let $x,y$ be vertices of $\blk(\mu)$. If $x$ and $y$ are joined by a $B(W_\alpha)$-edge in $B(X_\alpha)^{+B(W_\alpha)}$ but not a $B(X_\alpha)$-edge, then there exists simplices $\nu_x$ and $\nu_y$ of $B(X_\alpha)$ so that $\mu \subseteq \nu_x,\nu_y$ and $\nu_x \cup x$ and $\nu_y \cup y$ are adjacent vertices of $B(W_\alpha)$.
\end{prop}

\begin{proof}
	Assume first that $x,y$ are joined by a $B(W_\alpha)$-edge coming from a twist edges of $B(W_\alpha)$. Then there exists $c \in X^0_\alpha$ so that $x,y \in B(c)$. 
	Because $x,y \in \blk(\mu)$, the curve $c$ must be in either $\base(\mu)$ or $\xlk(\base(\mu))$. By Lemma \ref{lem:Extending_To_Pants_decomposition}, there exists a pants decomposition $\rho$ of $\Sdot - c$  so that $\base(\mu) \cap (\Sdot - c) \subseteq \rho$ and $\rho \subseteq X_\alpha$.  For each $d \in \rho$, select a marking arc $t_d \in B(d)$ so that $ t_d$ is the element of $\markcurve(\mu) \cap B(d)$ whenever $\markcurve(\mu) \cap B(d) \neq \emptyset$. Then, $\nu_x =\nu_y = \rho \cup \{t_d: d \in \rho\} \cup c$ are the desired simplices.
	
	Now assume that $x$ and $y$ are joined by a $B(W_\alpha)$-edge coming from a flip edge of $B(W_\alpha)$. In this case, $\supp(x)$ and $\supp(y)$ fill a complexity 1 subsurface $V$ of $S$ and there exist marking arcs $t_x \in B(\supp(x))$ and $t_y \in B(\supp(y))$ so that $d_{A_{\supp(x)}}(t_x,\supp(y)) \leq 1$ and $d_{A_{\supp(y)}}(t_y,\supp(x)) \leq 1$. Such marking arc $t_x,t_y$ are vertices of the maximal simplices in $B(X_\alpha)$ determined by the flip edge of $B(W_\alpha)$. Moreover $x$ (resp. $y$) is equal to either $t_x$ or $\supp(x)$ (resp. $t_y$ or $\supp(y)$).
	
	Since $x,y \in \blk(\mu)$, every curve of  $\base(\mu)$ must be equal to or disjoint from each of $\supp(x)$ and $\supp(y)$. Hence $\base(\mu)$ does not intersect $\partial V$. There then exists a pants decomposition $\rho$ of $\Sdot - V$ so that $\base(\mu) \cap (\Sdot - V) \subseteq \rho$ and $\rho \subseteq X_\alpha$ (Lemma \ref{lem:Extending_To_Pants_decomposition}). For each $c \in \rho \cup \partial V$, select a marking arc $t_c \in B(c)$ so that $t_c$ is the element of $\markcurve(\mu) \cap B(c)$ whenever $ \markcurve(\mu) \cap B(c) \neq \emptyset.$ 
	We now define a maximal simplex $\nu_x$ depending on the options for $x$.
	\begin{itemize}
		\item If $x =\supp(x)$, then $\nu_x = \rho \cup \partial V \cup \{t_c:c \in \rho \cup \partial V\} \cup t_x$.
		\item If $x = t_x$, then $\nu_x = \rho \cup \partial V \cup \{t_c:c \in \rho \cup \partial V\} \cup \supp(x)$.
	\end{itemize} 
	We similarly define $\nu_y$. Thus, $\nu_x \bowtie x$ and $\nu_y \bowtie y$ are joined by an edge of $B(W_\alpha)$.
\end{proof}

\subsection{Reducing the geometric conditions to the non-annular case}

Despite the addition  of the marking arcs, the proof that $(B(X_\alpha),B(W_\alpha))$ satisfies the geometric conditions of a combinatorial HHS is nearly identical to the proof for $(X_\alpha, W_\alpha)$. In this subsection, we will describe how to reduce the proof for $(B(X_\alpha),B(W_\alpha))$ to the case of $(X_\alpha, W_\alpha)$.

To start,  we prove the hyperbolicity of $B(X_\alpha)^{+B(W_\alpha)}$.
\begin{lem}\label{lem:B(X)_hyperbolic}
	$B(X_\alpha)^{+B(W_\alpha)}$ is $E_\alpha$-equivariantly quasi-isometric to both $B(X_\alpha)$ and $X_\alpha$. Hence $B(X_\alpha)^{+B(W_\alpha)}$ and $B(X_\alpha)$ are both uniform quasi-trees.
\end{lem}

\begin{proof}
	The definition of $B(W_\alpha)$ implies that the inclusion of $B(X_\alpha)$ into	$B(X_\alpha)^{+B(W_\alpha)}$ is a quasi-isometry. This is because if two vertices $x$ and $y$ of $B(X_\alpha)^{+B(W_\alpha)}$ are joined by a $B(W_\alpha)$-edge and not a $B(X_\alpha)$-edge, then $x \in \mu$ and $y \in \nu$ where $\mu$ and $\nu$ are maximal simplices that are joined by an edge of $B(W_\alpha)$. However, $\mu$ and $\nu$ will have at least one base curve in common and this curve will need to be at most distance 1 from both $x$ and $y$ in $B(X_\alpha)$. Thus  $d_{B(X_\alpha)}(x,y) \leq 2$.
	
	Next we check that the inclusion of $X_\alpha$ into $B(X_\alpha)$ is also a  quasi-isometry.  Consider the map $B^0(X_\alpha) \to X^0_\alpha$ that sends $x \to \supp(x)$. Lemma \ref{lem:forget_marking} showed this map is $(1,0)$-coarsely Lipschitz, making the inclusion of $X_\alpha$ into $B(X_\alpha)$ an isometric embedding. Since each vertex of $B(X_\alpha)$ is at most 1 away from a vertex of $X_\alpha$, the two spaces are quasi-isometric.
\end{proof}

We now prove the quasi-isometric embedding and hyperbolicity of the links of simplices of $B(X_\alpha)$. We will continue to use $\blk(\mu)$ to denote the link of $\mu$ in $B(X_\alpha)$ and we will use $\bsat(\mu)$ to denote the saturation of $\mu$ in $B(X_\alpha)$. Define $B(Y_\mu)$ to  be $B(X_\alpha)^{+B(W_\alpha)} - \bsat(\mu)$ and let  $BH(\mu)$ be the subset of $B(Y_\mu)$ spanned by $\blk(\mu)$---as before we are taking the link in $B(X_\alpha)$ then taking the span in $B(Y_\mu)$. The spaces $Y_{\base(\mu)}$ and $H(\base(\mu))$ are the spaces defined in Section \ref{subsec:geom_conditions}. All graphs are given their intrinsic path metrics.

\begin{prop}\label{prop:geometric_conditions_annular_case}
	For each non-maximal simplex $\mu$ of $B(X_\alpha)$, the space $BH(\mu)$ is hyperbolic and the inclusion of $BH(\mu)$ into $B(Y_\mu)$ is a uniform quasi-isometric embedding.  Moreover, $BH(\mu)$ falls into one of three cases:
	\begin{itemize}
		\item $BH(\mu)$ has bounded diameter because $\blk(\mu)$ is a join or a single vertex.
		\item $\base(\mu) =\supp(\mu)$, $\base(\mu)$ is not a pants decomposition of $\Sdot$, and $\mu$ contains no unmarked base curves. In this case, the inclusion of $H(\base(\mu))$ into $BH(\mu)$ is a uniform quasi-isometry.
		\item $\base(\mu)$ is a pants decomposition of $\Sdot$ and $\mu$ contains exactly one unmarked base curve $c$. In this case,  $BH(\mu)$ is equal to $\C(A_c)$.
	\end{itemize}

\end{prop}

We break the proof into four cases based on $\mu$. In the first two cases, the $B(X_\alpha)$-link of $\mu$ is join, so $BH(\mu)$ has bounded diameter, automatically satisfying the hyperbolicity and quasi-isometric embedding properties. In the third case, we reduce the problem to the results we established for $H(\base(\mu))$ in the non-annular case. The final case addresses links  consisting only of marking arcs. In this case,  we use the subsurface projection to annular complexes. 
\begin{proof}[Proof of Proposition \ref{prop:geometric_conditions_annular_case}] \
	\begin{enumerate}	
		
		\item Suppose $\supp(\mu) \neq \base(\mu)$. By Lemma \ref{lem:B-links}, $\blk(\mu)$ is either a join or a single vertex.

		\item Suppose $\supp(\mu) = \base(\mu)$, $\base(\mu)$ is not a pants decomposition of $\Sdot$, and there is an unmarked base curve $c \in \base(\mu)$. In this case, every vertex  of $\blk(\mu) - B(c)$ is joined by an edge to each vertex of $B(c)$, so $\blk(\mu)$ is a join.

		\item Suppose $\supp(\mu) = \base(\mu)$, $\base(\mu)$ is not a pants decomposition of $\Sdot$,  and $\mu$ has no unmarked base curves. In this case, $\blk(\mu)$ is a copy of $\xlk(\base(\mu))$ with the $B(c)$-vertices attached to each $X_\alpha$-vertex.    
		Thus, we have an inclusion of $Y^0_{\base(\mu)}$ into $B^0(Y_\mu)$ that sends $H^0(\base(\mu))$ into $BH^0(\mu)$. Further, this map can be extended over edges as follows: if $x,y \in Y^0_{\base(\mu)}$ are joined by an $X_\alpha$-edge (i.e. they are disjoint curves), then they are joined by a $B(X_\alpha)$-edge as well. If instead $x$ and $y$ are joined by a $W_\alpha$-edge, then the proof of Proposition \ref{prop:annular_C=C0} shows that we can also find two maximal simplices that are joined by a flip edge of $B(W_\alpha)$ and contain $x$ and $y$ respectively.
		
		We argue that this inclusion is actually a quasi-isometry on all of $Y_{\base(\mu)}$ and when restricted to $H(\base(\mu))$.		
		Consider the map  $B^0(Y_\mu) \to Y^0_{\base(\mu)}$  defined by $y \to \supp(y)$. An identical argument as given in Lemma \ref{lem:forget_marking} shows that this map induces a $(1,0)$-coarsely Lipschitz map $B(Y_\mu) \to Y_{\base(\mu)}$. Thus the inclusion $Y_{\base(\mu)} \to B(Y_\mu)$ is an isometric embedding by Lemma \ref{lem:coarsel_lipschitz}. However, since each vertex of $B(Y_\mu)$ is at most 1 away from a vertex of $Y_{\base(\mu)}$,  the inclusion of $Y_{\base(\mu)}$ into $B(Y_\mu)$ is a quasi-isometry. Further,  the subset $H(\base(\mu))$ is quasi-isometric to the subset $BH(\mu)$.
		
		Since $H(\base(\mu))$ is hyperbolic  and quasi-isometrically embeds in $Y_{\base(\mu)}$ by Proposition \ref{prop:QI_embedding_tree_quided_link}, we have that $BH(\mu)$ is hyperbolic and quasi-isometrically embeds in $B(Y_\mu)$. 
		Uniformity of the constants follows from the uniformity of the constants in Proposition  \ref{prop:QI_embedding_tree_quided_link}.

		\item Suppose $\supp(\mu) = \base(\mu)$ is a pants decomposition and $\unmark(\mu) \neq \emptyset$. If $\mu$ has more than one unmarked base curve, then $\blk(\mu)$ is a join of the $B(c)$ for these unmarked base curves and hence $BH(\mu)$ is bounded.  Thus we can assume there exists exactly one unmarked base curve $c$. Hence,  $\blk(\mu) = B(c)$ and $BH(\mu)$ is the graph with vertex set $B(c)$ and an edge between to vertices $t$ and $t'$ if $d_{A_c}(t,t')= 1$.  Since $B(c)$ is the vertex set of $\C(A_c)$, this means $BH(\mu) = \C(A_c)$  which is hyperbolic.

		We claim that every vertex $y \in B(Y_\mu) - BH(\mu)$ has $\supp(y)$ intersecting $c$. First, we know $\supp(y) \neq c$ as $\supp(y) = c$ would imply that $y \in B(c) = BH^0(\mu)$. If $\supp(y)$ was disjoint from $c$, then we could find a simplex $\nu$ of $B(X_\alpha)$ that contains $y$ and has $\blk(\nu) = B(c) = \blk(\mu)$. To do this, start with  a pants decomposition of $\Sdot -c$ that contains $\supp(y)$ and uses only curves from $X_\alpha$. Then select marking arcs for each of these curve so that $y$ is a marking arc if $y \in B(\supp(y))$. Since the existence of such a simplex $\nu$ would imply $y \not \in B(Y_\mu)$, it must be the case that $\supp(y)$ intersects $c$.
		
		Since the curve $\supp(y)$ intersects $c$ for each $y \in B(Y_\mu) - BH(\mu)$, the subsurface projection 	
		$\pi_{A_c}(\supp(y))$ is non-empty. Thus,  we can define $\psi_{\mu} \colon B^0(Y_\mu) \to BH(\mu) =\C(A_c)$ by $\psi_\mu(y) = \pi_{A_c}(\supp(y))$ when $y\not\in BH^0(\mu)$ and $\psi_\mu(y) = y$ when $y \in BH^0(\mu) = \C^0(A_c)$.  
		
		To prove $\psi_\mu$ is uniformly coarsely Lipschitz, we uniformly bound the distance between a pair of vertices of $B(Y_\mu)$ that are joined by an edge. If two curves $x,y \in B^0(Y_\mu) - BH^0(\mu)$ are joined by an edge of $B(Y_\mu)$ and $\supp(x) \neq \supp(y)$, then $i(\supp(x),\supp(y)) \leq 4$. By construction of the annular complex, $d_{A_c}(\supp(x),\supp(y))$ is uniformly bounded in terms of the intersection number $i(\supp(x),\supp(y))$. Thus $d_{BH(\mu)}(\psi_\mu(x),\psi(\mu))$ is uniformly bounded. If $y \in B^0(Y_\mu) - BH^0(\mu)$ and $x \in BH^0(\mu)$ are joined by an edge of $B(Y_\mu)$, then $d_{A_c}(\supp(y),x) \leq 1$ because they are joined by a $B(W_\alpha)$-edge coming from a flip move.   The map $\psi_{\mu}$ is therefore uniformly coarsely Lipschitz, proving $BH(\mu)$ is quasi-isometrically embedded in $B(Y_\mu)$.\qedhere
	\end{enumerate}
\end{proof}

\subsection{$E_\alpha$ is an HHG}

We now conclude our proof that $E_\alpha$ is a hierarchically hyperbolic group.

\begin{thm}\label{thm:E_is_an_HHG}
	The pair $(B(X_\alpha), B(W_\alpha))$ is a combinatorial HHS and the group $E_\alpha$ acts on $B(X_\alpha)$ with finitely many orbits of links of simplices. In particular, $E_\alpha$ is a hierarchically hyperbolic group.
\end{thm}

\begin{proof}
	We have proved the requirements for $(B(X_\alpha), B(W_\alpha))$ to be a combinatorial HHS in the following results:
	\begin{itemize}
		\item If $\mu_1,\dots,\mu_n$ is a sequence of simplices of $B(X_\alpha)$ so that $\lk(\mu_1) \subsetneq \dots \subsetneq \lk(\mu_n)$, then there is a corresponding sequence of properly nested subsurface of $\Sdot$ by Lemma \ref{lem:B-nesting} and Remark \ref{rem:nested_subsurfaces}. Hence, $n$ is bounded in terms of the complexity of $\Sdot$, which is determined by $S$.
		\item The hyperbolicity of $B(X_\alpha)^{+B(W_\alpha)}$ (Item \ref{CHHS:hyp_X}) was shown in Lemma \ref{lem:B(X)_hyperbolic}.
		\item The hyperbolicity and quasi-isometric embedding conditions (Item \ref{CHHS:geom_link_condition}) were verified in Proposition \ref{prop:geometric_conditions_annular_case}.
		\item The combinatorial conditions (Items \ref{CHHS:comb_nesting_condition} and \ref{CHHS:C=C_0_condition}) where shown in Propositions \ref{prop:annular_nesting_condition} and \ref{prop:annular_C=C0} respectively. Note, the condition proved is Proposition \ref{prop:annular_nesting_condition} is strictly stronger than what is required by the definition of a combinatorial HHS, since $\diam(BH(\mu)) >2$ implies $\blk(\mu)$ is not a join.
	\end{itemize}
	This implies that $B(W_\alpha)$ is connected and a hierarchically hyperbolic space. We now verify the additional requirements from Theorem \ref{thm:CHHS_are_HHS} for $E_\alpha$ to be a hierarchically hyperbolic group.
	
	Because the action of $E_\alpha$ on  $B(X_\alpha)$ is induced by the action (up to isotopy) of $E_\alpha$ on $\Sdot$, the action of $E_\alpha$ on the maximal simplices of $B(X_\alpha)$ induces the same action on $B(W_\alpha)$ as the action induced by $E_\alpha$ acting on $\Sdot$. It remains to verify that the action of $E_\alpha$ on $B(W_\alpha)$ is metrically proper and cobounded plus that there are finitely many $E_\alpha$-orbits of links of simplices of $B(X_\alpha)$.

	\textbf{Coboundedly.} Because $E_\alpha$ acts coboundedly on $W_\alpha$ it suffices to prove that for any $\mu,\nu \in B^0(W_\alpha)$, if $d_{W_\alpha}(\base(\mu),\base(\nu)) \leq 1$, then there exist $\phi,\psi \in E_\alpha$ so that $$d_{B(W_\alpha)}(\phi(\mu),\psi(\nu))$$ is uniformly bounded. 
	
	Suppose  $\mu,\nu \in B^0(W_\alpha)$ with $d_{W_\alpha}(\base(\mu),\base(\nu)) \leq 1$.  Recall, there exists universal number $D \geq 1$ (independent of $S$) so that  for each curve $b$ on $\Sdot$, the quotient of the annular complex $\C(A_b)$ by the subgroup generated by the Dehn twist around $b$ has diameter $D$. Thus, if $\base(\mu) = \base(\nu)$, then there exists  $\phi \in \mcg(S;z)$ so that $\phi$ is product of Dehn twists around curves in $\base(\mu)$, and  for each $b \in \base(\mu)$,  the marking arc of $\phi(\mu)$ for $b$ is $\C(A_b)$-distance $D$ from the marking arc of $\nu$ for $b$.  Further, $\phi \in E_\alpha$ since it is the product of Dehn twist about curves in $X_\alpha$. Hence,  by performing at most  $D$ twist moves on each curve of $\nu$, we can bound $d_{B(W_\alpha)}(\phi(\mu),\nu)$ by $D$ times the number of curves in a pants decomposition of $\Sdot$ as desired. 
	
	If instead $\base(\mu)$ is joined by an edge of $W_\alpha$ to $\base(\nu)$, then there must exist $c \in \base(\mu)$ and  $d \in \base(\nu)$ so that $(\base(\mu) - c) \cup d = \base(\nu)$ and $i(c,d) \leq 4$.  This implies $ \pi_{A_c}(d) \subseteq \C^0(A_c) =B(c)$ and $\pi_{A_d}(c) \subseteq \C^0(A_d) =B(d)$. 
	
	 Similarly to the previous case, there exists $\phi,\psi \in \mcg(S;z)$ so that $\phi(\mu)$ and $\psi(\nu)$ satisfy:
	\begin{enumerate}
		\item \label{item:In_E} $\base(\phi(\mu)) = \base(\mu)$ and  $\base(\psi(\nu)) = \base(\nu)$  because $\phi$ and $\psi$ are  products of Dehn twists about curves in $\base(\mu)$ and $\base(\nu)$ respectively. Since  each of these Dehn twists are around a curve in $X_\alpha$, $\phi$ and $\psi$ are both elements of $E_\alpha$.
		\item \label{item:twist}For each $b \in \base(\mu) - c = \base(\nu) - d$, the marking arc of $\phi(\mu)$ for $b$ is $\C(A_b)$-distance $D$ from the marking arc of $\psi(\nu)$ for $b$.
		\item \label{item:flip}The marking arc of $\phi(\mu)$ for $c$ is $\C(A_c)$-distance  $D$ from $\pi_{A_c}(d)$ and the marking arc of $\psi(\nu)$ for $d$ is $\C(A_d)$-distance $D$ from $\pi_{A_d}(c)$.
	\end{enumerate}

	By performing   at most $D$ twist moves  for each curve in $\base(\phi(\mu))$ and at most $D$ twist moves in the annular complex for $d \in \base(\psi(\nu))$ we can produce two markings that differ by a flip between $c$ and $d$. Thus, we have that $d_{B(W_\alpha)}(\phi(\mu),\psi(\nu))$ is bounded above by $2D$ times the number of curves in a pants decomposition of $\Sdot$. 
	
	\textbf{Metrically properly.} Our argument uses Masur and Minsky's graph of clean markings $\mc{M}(\Sdot)$ and the subsurface projection of markings. We direct the reader to  Section \ref{sec:surfaces} to  recall the necessary definitions and properties that we will need for working with $\mc{M}(\Sdot)$. 
	
	For each $\mu \in B(W_\alpha)$, let $\cl(\mu)$ be the set of clean markings that are compatible with $\mu$; see Theorem \ref{thm:Marking_complex}.  There is a uniform bound on the diameter of $\cl(\mu)$ in $\mc{M}(S)$ depending on $\Sdot$. We first show that the map $B(W_\alpha) \to \mc{M}(\Sdot)$ by $\mu \to \cl(\mu)$ is coarsely Lipschitz.
	
	Let $\mu,\nu \in B(W_\alpha)$ be joined by an edge. The definition of the edges plus Lemma \ref{lem:subsurface_projection_coarsely_Lipschtiz} ensures that $d_V(\mu,\nu)$ is uniformly bounded for every subsurface $V$ of $\Sdot$. This implies $d_V(\cl(\mu),\cl(\nu))$ is uniformly bounded because for each $\mu' \in \cl(\mu)$ and $\nu' \in \cl(\nu)$ and for each subsurface $V$, we have \[d_V(\mu',\nu') \leq d_V(\mu',\mu) + d_V(\mu,\nu) + d_V(\nu,\nu') + 6 \leq d_V(\mu,\nu) +12. \]     By applying the second item in Theorem \ref{thm:Marking_complex}, we can therefore find  $D\geq 0$, depending only on $\Sdot$, so that for any $\mu,\nu \in B(W_\alpha)$   joined by an edge we have \[ d_{\mc{M}(S)}(\cl(\mu),\cl(\nu)) \leq D.\] Further, by possibly enlarging $D$, we can assume $\diam_{\mc{M}(S)}(\cl(\mu)) \leq D$ for each $\mu \in B(W_\alpha)$. By Lemma \ref{lem:checking_coarsely_lipschitz}, the map $B(W_\alpha) \to \mc{M}(\Sdot)$ by $\mu \to \cl(\mu)$ is $(3D,0)$-coarsely Lipschitz.
	
	Now, let $\ball_r^B(\cdot)$ and $\ball_r^\mc{M}(\cdot)$ denote the balls of radius $r$ in $B(W_\alpha)$ and $\mc{M}(S)$ respectively.  For each $r \ge 0$ and $\mu \in B(W_\alpha)$, the set $\{g \in E_\alpha : g \cdot \ball_r^B(\mu)\cap \ball_r^{B}(\mu) \neq \emptyset\}$ is a subset of  $$\{g \in E_\alpha : g \cdot \ball_{3Dr}^\mc{M}(\cl(\mu))\cap \ball_{3Dr}^{\mc{M}}(\cl(\mu)) \neq \emptyset \}.$$ Because $E_\alpha$ is a subgroup of $\mcg(S;z)$, $E_\alpha$  acts metrically properly on $\mc{M}(S)$ (Theorem \ref{thm:Marking_complex}). Thus,  $$\{g \in E_\alpha : g \cdot \ball_{3Dr}^\mc{M}(\cl(\mu))\cap \ball_{3Dr}^{\mc{M}}(\cl(\mu)) \neq \emptyset\}$$ contains finitely many elements of $E_\alpha$. This implies  $\{g \in E_\alpha : g \cdot \ball_r^B(\mu)\cap \ball_r^{B}(\mu) \neq \emptyset\}$ also contains only finitely many elements of $E_\alpha$, making the action of $E_\alpha$ on $B(W_\alpha)$  metrically proper.

	\textbf{Finitely many orbits of links.} For a simplex $\mu$ of $B(X_\alpha)$, define $\overline{\mu} = \supp(\mu) \bowtie \markcurve(\mu)$. Note, for any simplex $\mu$, we have $\base(\overline{\mu}) = \supp(\overline{\mu})  =\supp(\mu)$.  Because $\base(\overline{\mu}) = \supp(\overline{\mu})$, Lemma \ref{lem:B-nesting} implies that $\blk(\overline{\mu})$ is determined by the subsurface $U_{\supp(\overline{\mu})}$ and the unmarked base curves of $\mu$. Because there are finitely many $E_\alpha$-orbits of both curves in $X_\alpha$ and subsurfaces $V\subseteq \Sdot$ where $\partial V \subseteq X_\alpha$, there can only be a finite number of  $E_\alpha$-orbits of $B(X_\alpha)$-links of  the $\overline{\mu}$'s.  By Lemma \ref{lem:B-links}, $\blk(\mu) - \blk(\overline{\mu}) = \supp(\mu) - \base(\mu)$ for each $\mu$. Since $\supp(\mu) - \base(\mu)$ is a  multicurve of curves in $X_\alpha$, there are only finitely many $E_\alpha$-orbits of possibilities for   $\supp(\mu) - \base(\mu)$ as well. Thus, there are only finitely many $E_\alpha$-orbits of $B(X_\alpha)$-links of simplices.
\end{proof}

\section{A description of the HHG structure for $E_\alpha$}\label{sec:structure}	We now provide an explicit description of the hierarchically hyperbolic group structure that $E_\alpha$ receives from Theorem \ref{thm:E_is_an_HHG}.  We start by recalling the defining information of a hierarchically hyperbolic structure (Sections \ref{subsec:HHS_Data}).  We embrace a slightly non-standard set of notation in this subsection to draw a distinction between  curve complexes/subsurface projection maps and the hyperbolic spaces/projection maps of an abstract HHS.
Next we describe the hierarchically hyperbolic structure that Theorem \ref{thm:CHHS_are_HHS} imparts on an abstract combinatorial HHS (Section \ref{subsec:HHS_structure_for_cHHS}). 
The two subsequent subsections take this hierarchically hyperbolic structure for the space $W_\alpha$ (Section \ref{subsec:structure_for_W}) and the group $E_\alpha$ (Section \ref{subsec:structure_for_E}) and reinterpret it using the topology of the surface .  The final subsection uses this specific structure to show that some previously established results in the literature apply to $E_\alpha$ (Section \ref{subsec:applications}).

This section is largely expository and intentionally does not provide the full definition of hierarchically hyperbolic spaces and groups.  We direct the reader to \cite{BHS_HHSII} and \cite{HHS_survey} for detailed discussion of hierarchically hyperbolic spaces and to  \cite{BHMS_cHHS} for combinatorial hierarchically hyperbolic spaces.

\subsection{The defining data of an HHS structure}\label{subsec:HHS_Data}
An \emph{$E$-hierarchically hyperbolic space structure} on a geodesic metric space $\mc{X}$ is a set $\mf{S}$ indexing a collection of $E$-hyperbolic spaces $\{\mc{H}(V)\}_{V \in \mf{S}}$ and equipped with three relations:  nesting ($\nest$), orthogonality ($\perp$), and  transversality ($\trans$). For each $V \in \mf{S}$, there exists a surjective $(E,E)$-coarsely Lipschitz map $\tau_V \colon \mc{X} \to \mc{H}(V)$ called the \emph{projection to $V$}.  The relation $\nest$ is a partial order on $\mf{S}$ that contains a unique maximal element. Whenever $V \propnest U$ in $\mf{S}$ there exists a distinguished  subset $\rho_U^V$ in $\mc{H}(U)$ with diameter at most $E$. The relations $\perp$ and $\trans$ are both symmetric and anti-reflexive. Whenever $V \trans U$, there exist  distinguished  subsets $\rho_U^V \subseteq \mc{H}(U)$ and $\rho_V^U \subseteq\mc{H}(V)$ each of diameter at most $E$. When $V \propnest U$ or $V \trans U$, we call the subset $\rho_U^V$ the \emph{relative projection from $V$ to $U$}. Every pair of distinct elements of $\mf{S}$ is related by exactly one of the relations $\nest$, $\perp$, or $\trans$. To be a hierarchically hyperbolic space structure for $\mc{X}$, the set $\mf{S}$ and these relations and projection need to satisfy a number of axioms; see \cite{BHS_HHSII}  or \cite{HHS_survey} for a complete definition.

An HHS structure $\mf{S}$ for a geodesic space $\mc{X}$ is an \emph{hierarchically hyperbolic group structure} for a finitely generated group $G$ if $G$ acts metrically properly and coboundedly on $\mc{X}$ and there is an $\nest$, $\perp$, and $\trans$ preserving action of $G$ on $\mf{S}$ by bijections  satisfying:
\begin{enumerate}
	\item $\mf{S}$ has finitely many $G$-orbits;
	\item For each $V \in \mf{S}$ and $g\in G$, there exists an isometry $g_V \colon \mc{H}(V) \rightarrow \mc{H}(gV)$ satisfying the following for all $V,U \in \mf{S}$ and $g,h \in G$.
	\begin{itemize}
		\item The map $(gh)_V \colon \mc{H}(V) \to \mc{H}(ghV)$ is equal to the map $g_{hV} \circ h_V \colon \mc{H}(V) \to \mc{H}(ghV)$.
		\item For each $x \in \mc{X}$, $g_V(\tau_V(x))$ and $\tau_{gV}(g \cdot x)$ are at most $E$ far apart in $\mc{H}(gV)$.
		\item If $U \trans V$ or $U \propnest V$, then $g_V(\rho_V^U)$  and $\rho_{gV}^{gU}$ are at most $E$ far apart in $\mc{H}(gV)$.
	\end{itemize}
\end{enumerate}

\subsection{The HHS structure for combinatorial HHSs}\label{subsec:HHS_structure_for_cHHS}
Let $(X,W)$ be a combinatorial HHS. Recall the equivalence relation on simplices of $X$: $\Delta \sim \Delta'$ if $\lk(\Delta) = \lk(\Delta')$. Let $[\Delta]$ denote the $\sim$-equivalence class of the simplex $\Delta$. We adopt the convention that the empty set is a simplex of $X$ whose link is all of $X$. Let $\mf{S}$ be the set of $\sim$-equivalence classes of all non-maximal simplices of $X$, including the empty set.

This set $\mf{S}$ is the index set for the HHS structure for $W$ provided by Theorem \ref{thm:CHHS_are_HHS}.  We define $[\Delta] \nest [\Delta']$ if $\lk(\Delta) \subseteq \lk(\Delta')$ and define $[\Delta] \perp [\Delta']$ if $\lk(\Delta) \subseteq \lk(\lk(\Delta'))$---or equivalently $\lk(\Delta') \subseteq \lk(\lk(\Delta))$.  We define $[\Delta] \trans [\Delta']$  if $[\Delta] \not \perp [\Delta']$ and neither is nested into the other.
The $\nest$-maximal element of $\mf{S}$ is $[\emptyset]$ and its associated hyperbolic space is $X^{+W}$. For each $[\Delta] \in \mf{S} - \{[\emptyset]\}$, the associated hyperbolic space is $H(\Delta)$.  Note, this means the hyperbolic space associated to $[\Delta]$ has bounded diameter whenever $\lk(\Delta)$ is a join or a single vertex.

The projection map $\tau_{[\Delta]} \colon W \to H(\Delta)$ is defined as follows. Let $\sigma \in W^0$. Since $\sigma$ is a maximal simplex in $X$, there must be a vertex of $\sigma$ contained in the space $Y_\Delta = X^{+W} - \sat(\Delta)$. The space $Y_\Delta$ is shown to be  hyperbolic in Section 3 of  \cite{BHMS_cHHS}. Since $H(\Delta)$ quasi-isometrically embeds into $Y_\Delta$, there is  a coarse closest point projection $\mf{p}_\Delta \colon Y_\Delta \to H(\Delta)$.  Define $\tau_{[\Delta]} (\sigma) = \mf{p}_{\Delta} (\sigma \cap Y^0_\Delta)$. The relative projections  are similarly defined by $\rho_{[\Delta]}^{[\Delta']} = \mf{p}_\Delta (Y^0_{\Delta'} \cap Y^0_{\Delta})$ whenever $[\Delta] \propnest [\Delta']$ or $[\Delta] \trans [\Delta']$. The hyperbolic spaces and projections are well defined because $Y_\Delta =Y_{\Delta'}$ and $H(\Delta) = H(\Delta')$ whenever $\Delta \sim \Delta'$.

\subsection{The HHS structure for $W_\alpha$}\label{subsec:structure_for_W}
Let $\mf{S}_\alpha$ be the HHS structure given to $W_\alpha$  by virtue of $(X_\alpha,W_\alpha)$ being a combinatorial HHS. We give an interpretation of the relations and hyperbolic spaces for elements of $\mf{S}_\alpha$ using the topology of subsurfaces  and curves of $\Sdot$.

By Lemma \ref{lem:top_CHHS_dictionary}, there is a bijection between links of simplices $\mu$ in $X_\alpha$ and the complimentary subsurfaces $U_\mu$. Thus, we can identify the index set $\mf{S}_\alpha$ with the set of subsurfaces $\{U_\mu: \mu \text{ is a non-maximal simplex of } X_\alpha\}$.     Lemma \ref{lem:top_CHHS_dictionary} also says that the nesting of links $\lk(\mu) \subseteq \lk(\mu')$ is equivalent to the containment of subsurfaces $U_\mu \subseteq U_{\mu'}$.  Lemma \ref{lem:disjoint_iff_orthogonal} below verifies that  $\lk(\mu) \subseteq \lk(\lk(\mu'))$ if and only if $U_\mu$ and $U_{\mu'}$ are disjoint. Thus the nesting, orthogonality, and transversality of the elements of $\mf{S}_\alpha$ respectively correspond to the containment, disjointness, and overlapping of the subsurfaces in $\{U_\mu: \mu \text{ is a non-maximal simplex of } X_\alpha\}$.

The hyperbolic spaces associated to subsurfaces come in two types. If $\mu$ is of subsurface type, then the hyperbolic space $H(\mu)$ is the  modified curve graph $\C'(U_\mu)$ given in Definition \ref{defn:C'(S)}.  If $\mu$ is of tree guided type, then the hyperbolic space $H(\mu)$ is the subset of $\C'(\Sdot -\mu)$ spanned by the curves $c$ with $\Pi(c)$ disjoint from $\alpha$. This subset is quasi-isometric to tree.

\begin{lem}\label{lem:disjoint_iff_orthogonal}
	For any non-maximal simplices $\mu,\mu' \subseteq X_\alpha$, we have  $\lk(\mu) \subseteq \lk(\lk(\mu'))$ if and only if $U_\mu$ and $U_{\mu'}$ are disjoint
\end{lem}

\begin{proof}
	Since $\lk(\mu)$ and $\lk(\mu')$ are spanned by the curves of $X_\alpha$ that are on  $U_\mu$ and $U_{\mu'}$ respectively, if $U_\mu$ and $U_{\mu'}$ are disjoint, then  $\lk(\mu) \subseteq \lk(\lk(\mu'))$ because every curve in $\lk(\mu)$ is disjoint from every curve in $\lk(\mu')$.
	
	Conversely if $\lk(\mu) \subseteq \lk(\lk(\mu'))$, then $U_\mu$ must be disjoint from $U_{\mu'}$ since $U_\mu$ and $U_{\mu'}$ are filled by the curves in $\lk(\mu)$ and $\lk(\mu')$ respectively.
\end{proof}

\subsection{The HHG structure for $E_\alpha$}\label{subsec:structure_for_E} 
Let $\mf{BS}_\alpha$ denote the HHS structure for $B(W_\alpha)$ arising from $(B(X_\alpha),B(W_\alpha))$ being a combinatorial HHS. By Theorem \ref{thm:E_is_an_HHG}, $\mf{BS}_\alpha$ is also an HHG structure for $E_\alpha$. Let $\mf{BS}^T_\alpha$ denote the set of equivalence classes of simplices $\mu \subseteq B(X_\alpha)$ with $\supp(\mu) = \base(\mu)$.  We start by describing in topological terms the HHS data for equivalence classes of simplices in $\mf{BS}^T_\alpha$. At the end we will describe why the other simplices can essentially be ignored when working with the HHG structure $\mf{BS}_\alpha$ for $E_\alpha$.

Define $C_\mu$ to be the disjoint union the set of annuli $\{A_c : c \in \unmark(\mu)\}.$  Lemma \ref{lem:B-nesting} implies that $\blk(\mu) = \blk(\mu')$ if and only if $U_{\supp(\mu)} = U_{\supp({\mu'})}$ and $C_\mu = C_{\mu'}$ . Thus, the set   $\mf{BS}^T_\alpha$ can be identified with the set of all pairs of subsurfaces $(U_{\supp(\mu)},C_\mu)$.  Lemma \ref{lem:B-nesting} showed that $ \blk(\mu) \subseteq \blk(\mu')$, and hence the $\nest$ relation, is equivalent to the following ``dictionary nesting'' of subsurfaces: \[(U_{\supp(\mu)},C_\mu) \nest (U_{\supp(\mu')},C_{\mu'}) \iff U_{\supp(\mu)} \subseteq U_{\supp(\mu')} \text{ and } C_\mu \subseteq U_{\supp(\mu')} \cup C_{\mu'}.\]
Lemma \ref{lem:annular_disjoint_iff_orthogonal} below shows that  $\blk(\mu) \subseteq \blk(\blk(\mu'))$ if and only if $U_{\supp(\mu)} \cup C_\mu$ is disjoint from $U_{\supp(\mu')} \cup C_{\mu'}$.  Hence \[ (U_{\supp(\mu)},C_\mu) \perp (U_{\supp(\mu')},C_{\mu'}) \iff U_{\supp(\mu)} \cup C_\mu \text{ and } U_{\supp(\mu')} \cup C_{\mu'} \text{ are disjoint.}\]

Proposition \ref{prop:geometric_conditions_annular_case} showed that there are three possibilities for the hyperbolic space $BH(\mu)$ when $[\mu] \in \mf{BS}_\alpha$. When $U_{\supp(\mu)} \cup C_\mu$ is not connected, then  $BH(\mu)$ has bounded diameter since $\blk(\mu)$ is a join. If $C_\mu = \emptyset$ and $U_{\supp(\mu)}$ is connected, then  $BH(\mu)$  is quasi-isometric to $H(\supp(\mu))$, which in turn falls into one of the two case described in the previous section. Finally, when $U_{\supp(\mu)} = \emptyset$ and $C_\mu$ is an single annulus $A_c$, then $BH(\mu)$ is the annular complex $\C(A_c)$.

Finally consider a simplex $\mu$ where $\base(\mu) \neq \supp(\mu)$. For such a simplex, define $\overline{\mu} = \supp(\mu) \bowtie \markcurve(\mu)$. Then $\blk(\mu) =\blk(\overline{\mu}) \bowtie (\supp(\mu) - \base(\mu))$. Thus, the hyperbolic spaces associated to equivalence classes of such simplices will have finite diameter.  Further, for any simplex $\mu'$ if $\blk(\mu') \subseteq \blk(\mu)$, then either $\blk(\mu')$ is a join or $\blk(\mu') \subseteq \blk(\overline{\mu})$. With these two facts, one can essentially ignore the elements of $\mf{BS}_\alpha$ with $\base(\mu) \neq \supp(\mu)$  when working with the hierarchically hyperbolic structure of $B(W_\alpha)$ (or $E_\alpha$); Lemma \ref{lem:E_unbounded} in the next section is an example of this philosophy. In fact, one could use these two facts to prove that $\mf{BS}_\alpha^T$ itself an HHS structure for $B(W_\alpha)$. We have forgone this work as we have no need for it.

\begin{lem}\label{lem:annular_disjoint_iff_orthogonal}
	Let $\mu,\mu' \subseteq B(X_\alpha)$ be non-maximal simplices with $\supp(\mu) = \base(\mu)$ and $\supp(\mu') = \base(\mu')$.  Then  $\blk(\mu) \subseteq \blk(\blk(\mu'))$ if and only if $U_{\supp(\mu)} \cup C_\mu$ and $U_{\supp(\mu')} \cup C_{\mu'}$ are disjoint.
\end{lem}

\begin{proof}
We start by describing $\blk(\blk(\eta))$ when $\eta \subseteq B(X_\alpha)$ is a simplex with $\supp(\eta) = \base(\eta)$. By Lemma \ref{lem:B-links}, $\blk(\eta)$ is spanned by the set of vertices:
\[ \left\{x \in B(X_\alpha):  \supp(x) \in \xlk(\supp(\eta))\right\} \cup \bigcup_{c  \in \unmark(\eta)} B(c). \tag{$\star$} \label{tag:links}\]
	 Let $y \in \blk(\blk(\eta))$. There are two facts we shall need about $y$.
	 \begin{enumerate}[I.]
	 	\item \label{item:disjiont_support}  $\supp(y)$ is disjoint from $U_{\supp(\eta)}$: to start,  $\supp(y) \not\in \xlk(\supp(\eta))$, since $y \not\in \blk(\eta)$. Thus, the only way for $y$ to joined by edge to each vertex of $\blk(\eta)$, is for $\supp(y)$ to be disjoint from  every curve in $\xlk(\supp(\eta))$. Since $U_{\supp(\eta)}$ is defined to be the subsurface filled by the curves in $\xlk(\supp(\eta))$, we must have that $y$ is disjoint from $U_{\supp(\eta)}$.  
	 	\item \label{item:unmarked_curves} for each unmarked base curve $c$ of $\eta$, either $y = c$ or $\supp(y)$ is disjoint from $c$: if $c$ is an unmarked base curve of $\eta$, then $B(c) \subseteq  \blk(\eta)$. Hence, $y$ must be joined by an edge to each vertex of $B(c)$. The only way for that to happen is if $y = c$ or $\supp(y)$ is disjoint from $c$.
	 \end{enumerate}

	Assume $\blk(\mu) \subseteq \blk(\blk(\mu'))$. Because of (\ref{tag:links}),  $U_{\supp(\mu)} \cup C_\mu$ is precisely the subsurface filled by the set of curves $\{\supp(y) : y \in \blk(\mu)\}$. By Item \ref{item:disjiont_support}, $\supp(y)$ is disjoint from $U_{\supp(\mu')}$ for all $ y \in \blk(\mu)$.  Hence $U_{\supp(\mu)} \cup C_\mu$ is also disjoint from $U_{\supp(\mu')}$.
	
	To show that  $U_{\supp(\mu)} \cup C_\mu$ is disjoint from $C_{\mu'}$, we will prove that every curve in the set $\{\supp(y) : y \in \blk(\mu)\}$ is disjoint from every unmarked base curve of $\mu'$. This implies $U_{\supp(\mu)} \cup C_\mu$ is disjoint from $C_{\mu'}$, since the unmarked base curves of $\mu'$ fill $C_{\mu'}$. By Item \ref{item:unmarked_curves}, if $y \in \blk(\mu)$ and $c'$ is an unmarked base curve of $\mu'$, then either $y = c'$ or $\supp(y)$ is disjoint from $c'$.  Since $\supp(\mu) = \base(\mu)$,  (\ref{tag:links}) implies that $B(c') \subseteq \blk(\mu)$ if $y = c'$. However, $B(c') \subseteq \blk(\mu')$ by (\ref{tag:links}), implying  $B(c') \not\subseteq \blk(\blk(\mu'))$. Since $\blk(\mu) \subseteq \blk(\blk(\mu'))$, this means $y \neq c'$ and hence $\supp(y)$ must be disjoint from $c'$ as desired.
	
	Now assume that $U_{\supp(\mu)} \cup C_\mu$ and $U_{\supp(\mu')} \cup C_{\mu'}$ are disjoint and let $y \in \blk(\mu)$. Hence we have either $\supp(y) \subseteq U_{\supp(\mu)}$ or $\supp(y) \subseteq C_\mu$. Thus,  $\supp(y)$ is disjoint from both $U_{\supp(\mu')}$ and $C_{\mu'}$. Since $\supp(y)$ is disjoint from  $U_{\supp(\mu')}$, $\supp(y)$ is  disjoint from every curve in $\xlk(\supp(\mu'))$ and hence $y$ joined be an edge to every vertex of $$\{x \in B^0(X_\alpha) : \supp(x) \in \xlk(\supp(\mu'))\}.$$ 
	Similarly,  $\supp(y)$ begin disjoint from $C_{\mu'}$ implies that $\supp(y)$ is disjoint from each unmarked base curve of $\mu'$. But this implies $y$ is joined by an edge to each vertex in $B(c')$ for each $c' \in \unmark(\mu')$. By (\ref{tag:links}), this makes $y \in \blk(\blk(\mu'))$
\end{proof}

\subsection{Some applications of the HHS structure}\label{subsec:applications}
We now use the HHG structure $\mf{BS}_\alpha$ to prove the remaining statements from Theorem \ref{intro_thm:quasi-tree} of the introduction.  These results rely on two facts: (1)  $X_\alpha^{+W_\alpha}$ is a quasi-tree (2) a minor modification of $\mf{BS}_\alpha$ produces an HHG structure of $E_\alpha$ that has an additional property called \emph{unbounded products} originally defined by Abbott, Behrstock, and Durham \cite{ABD}.  The definition we give below of unbounded products  is equivalent to the original definition in the setting of hierarchically hyperbolic groups and avoids having to describe  unneeded additional aspects of the theory of hierarchically hyperbolic spaces. The proof of the equivalence is a straight forward application of the distance formula in an HHG \cite[Theorem 4.5]{BHS_HHSII} and the fact that $G$ acts on $\mf{S}$ with finitely many orbits.

\begin{defn}
	A hierarchically hyperbolic group $(G,\mf{S})$ has \emph{unbounded products} if for all non-$\nest$-maximal $U \in \mf{S}$, whenever there exists $V \nest U$ with $\diam(\mc{H}(V)) = \infty$, there also exists  $Q \perp U$ so that $\diam(\mc{H}(Q)) =\infty$.
\end{defn}

Given any hierarchically hyperbolic group $(G,\mf{S})$, Abbott, Behrstock, and Durham provide an explicit construction of  a new hierarchically hyperbolic group structure $\mf{T}$ with unbounded products for $G$. In certain cases, this new structure maintains the same $\nest$-maximal hyperbolic space.

\begin{thm}[{\cite[Theorem 3.7]{ABD}}]\label{thm:ABD}
	Let $(G,\mf{S})$ be a hierarchically hyperbolic group. There exists a hierarchically hyperbolic groups structure $\mf{T}$ for $G$ with unbounded products. Moreover, if for  every non-$\nest$-maximal $U \in \mf{S}$ with  $\diam(\mc{H}(U)) = \infty$ there is $V \in \mf{S}$ so that $V \perp U$ and $\diam(\mc{H}(V)) = \infty$, then the $\nest$-maximal hyperbolic space of $\mf{T}$ can be taken to be the $\nest$-maximal hyperbolic space of $\mf{S}$.
\end{thm}

\begin{remark}
	While the moreover clause of Theorem \ref{thm:ABD} is not given in the statement in  \cite{ABD}, it is explicit in Abbott, Behrstock, and Durham's construction of the hyperbolic spaces for the new structure $\mf{T}$. 
\end{remark}

Our HHG structure $\mf{BS}_\alpha$ for $E_\alpha$ does not have unbounded products itself, but  the next lemma verifies that $\mf{BS}_\alpha$ does satisfy the hypotheses of the moreover clause of Theorem \ref{thm:ABD}. This allows us to use  Abbott, Behrstock, and Durham's construction without changing the $\nest$-maximal hyperbolic space.

\begin{lem}\label{lem:E_unbounded}
	For every $[\mu] \in \mf{BS}_\alpha$, if $\mu \neq \emptyset$ and $\diam(BH(\mu)) = \infty$,  then there exist $[\nu] \in \mf{S}_\alpha$ with $ [\nu] \perp [\mu]$ and $\diam(BH(\nu)) =\infty$.
\end{lem}

\begin{proof} 
	Let $\mu$ be a non-maximal, non-empty simplex of $\mf{BS}_\alpha$ so that $\diam(BH(\mu)) = \infty$. If $\supp(\mu) \neq \base(\mu)$, then Lemma \ref{lem:B-links} makes $\blk(\mu)$ a join. Since this would imply that $\diam(BH(\mu)) \leq 2$, we know $\supp(\mu) = \base(\mu)$. 
	
	Recall, $C_\mu$ is the disjoint union of  annuli whose core curves are the unmarked base curves of $\mu$.   By Lemma \ref{lem:B-links}, $U_{\supp{(\mu)}} \cup C_\mu$ is filled by the set of curves $\{\supp(y) : y \in \blk(\mu)\}$. Hence, if  $U_{\supp{(\mu)}} \cup C_\mu$ is disconnected, then $\blk(\mu)$ would have to be a join. Since this would imply $BH(\mu)$ has bounded diameter, we know $U_{\supp{(\mu)}} \cup C_\mu$ is connected.
	
	Since $U_{\supp{(\mu)}} \cup C_\mu$ is connected, we have two cases: either $C_\mu = \emptyset$ and $U_{\supp(\mu)}$ is connected or $U_{\supp(\mu)} = \emptyset$ and $C_\mu$ is an single annulus. In either case, there is a curve $c \in X^0_\alpha$ that is disjoint from $U_{\supp{(\mu)}} \cup C_\mu$ and is not the core curve of $C_\mu$ in the case $U_{\supp(\mu)} = \emptyset$.  Lemma \ref{lem:Extending_To_Pants_decomposition} implies that there exists a simplex $\sigma$ of $B(X_\alpha)$ so that
	\begin{itemize}
		\item $c \in \base(\sigma)$;
		\item  $\base(\sigma)$ is a pants decomposition of $\Sdot$;
		\item $c$ is the only unmarked base curve of $\sigma$.
	\end{itemize}
	The $B(X_\alpha)$-link of $\sigma$ is exactly $B(c)$, and hence $BH(\sigma)$ has infinite diameter as it is a quasi-isometric to a line (Proposition \ref{prop:geometric_conditions_annular_case}). Further, $[\sigma] \in \mf{BS}^T_\alpha$, $U_{\supp(\sigma)} = \emptyset$, and $C_\sigma$ is precisely the annulus with core curve $c$. Since $C_\sigma$ is disjoint from $U_{\supp{(\mu)}} \cup C_\mu$ and $U_{\supp(\sigma)} = \emptyset$, we have $[\sigma] \perp [\mu]$ by Lemma \ref{lem:annular_disjoint_iff_orthogonal}. 
\end{proof}

Abbott, Behrstock, and Durham prove the following results about HHGs with unbounded products; see \cite{ABD} for the relevant definitions. 

\begin{thm}[\cite{ABD}]\label{thm:largest_action}
	Let $(G,\mf{T})$ be a hierarchically hyperbolic group with unbounded products and let $T$ be the $\nest$-maximal element of $\mf{T}$.
	\begin{itemize}
		\item The action of $G$ on $\mc{H}(T)$ is the largest cobounded acylindrical  action of $G$ on a  hyperbolic space.
		\item A subgroup  $H<G$ is  stable  if and only if the orbit map of $H$ into $\mc{H}(T)$ is a quasi-isometric embedding.
	\end{itemize}
\end{thm}

Cordes, Charney, and Sisto recently provided a complete topological characterization of the Morse boundaries of certain groups \cite{CCS_omega_Morse_boundary}. While they do not tackle hierarchically hyperbolic groups directly, their methods are readily adapted to  certain hierarchically hyperbolic groups with unbounded products. In the appendix, we use their techniques along with the work of Abbott, Behrstock, and Durham to establish the following; see the appendix or \cite{CCS_omega_Morse_boundary} for the definition of the Morse boundary and $\omega$-Cantor space.

\begin{thm}\label{thm:Morse_boundary}
	Let $(G,\mf{T})$ be a hierarchically hyperbolic group with unbounded products and let $T$ be the $\nest$-maximal element of $\mf{T}$. If $\mc{H}(T)$ is a quasi-tree and $G$ is not hyperbolic, then the Morse boundary of $G$ is homeomorphic to an $\omega$-Cantor space.	
\end{thm}

\begin{proof}
	This is  Corollary \ref{cor:HHG_with_omega_boundary} in the appendix.
\end{proof}

Combining these results, we can provide a proof of the remaining claims of Theorem \ref{intro_thm:quasi-tree}.

\begin{cor}\
	\begin{enumerate}
		\item\label{item:action} The action of $E_\alpha$ on $X_\alpha$ is the largest cobounded acylindrical action of $E_\alpha$ on a hyperbolic space.
		\item\label{item:stable} A subgroup  $H<E_\alpha$ is  stable if and only if the orbit map of $H$ into $X_\alpha$ is a quasi-isometric embedding. In particular, every stable subgroup of $E_\alpha$ is virtually free.
		\item\label{item:Morse} The Morse boundary of $E_\alpha$ is an $\omega$-Cantor set.
	\end{enumerate}
\end{cor}

\begin{proof}
	Combining Lemma \ref{lem:E_unbounded} with Theorem \ref{thm:ABD}, there is an HHG structure $\mf{T}$ for $E_\alpha$ where the $\nest$-maximal hyperbolic space of $\mf{T}$ is the $\nest$-maximal hyperbolic space of $\mf{BS}_\alpha$. Lemma \ref{lem:B(X)_hyperbolic} showed that the $\nest$-maximal element of $\mf{BS}_\alpha$ is $E_\alpha$-equivariantly quasi-isometric to the quasi-tree $X_\alpha$. 
	The conclusions then follow from Theorems \ref{thm:largest_action} and \ref{thm:Morse_boundary} plus the fact that any finitely generated group whose orbit map in a quasi-tree is a quasi-isometric embedding is virtually free.  
\end{proof}

\begin{appendices}

\section{Appendix: HHGs and $\omega$-Cantor set boundaries}

\theoremstyle{plain}
\newtheorem{lemappend}{Lemma}
\renewcommand*{\thelemappend}{A.\arabic{lemappend}}
\newtheorem{thmappend}[lemappend]{Theorem}
\newtheorem{claimappend}[lemappend]{Claim}
\newtheorem{corappend}[lemappend]{Corollary}

\theoremstyle{definition}
\newtheorem{defnappend}[lemappend]{Definition}
\newtheorem{remappend}[lemappend]{Remark}
\newtheorem{exappend}[lemappend]{Example}

Charney, Cordes, and Sisto  recently gave the first topological descriptions of the Morse boundaries of non-hyperbolic groups \cite{CCS_omega_Morse_boundary}. They proved the Morse boundaries of most right-angled Artin groups and  graph manifold groups are homeomorphic to a specific  limit of Cantor spaces that they call an $\omega$-Cantor space. This appendix describes how their techniques can be straight forwardly  extended to prove the same result for a broader class of groups. The main motivation is to show that certain hierarchically hyperbolic groups will also have $\omega$-Cantor space boundaries. However, as the full power of hierarchical hyperbolicity will not be needed for our proof, we will instead work in the simpler setting of \emph{Morse detectable} groups introduced by the author, Spriano, and Tran \cite{Morse-local-to-global}.

\subsection*{Notation}
Throughout, $G$ will be a finitely generated group with a fixed finite generating set. Whenever we discuss $G$ as a metric space, we are referring to the path metric on the Cayley graph of $G$ with respect to this finite generating set. We will denote the Morse boundary of $G$ by $\partial_M G$ and the Gromov boundary of  a hyperbolic space $X$ by $\partial X$. Throughout, $I$ will denote a closed interval of the real line. We remark that we are intentionally glib about the  theory of Morse boundaries and direct the reader to Cordes' excellent survey on the topic for more details \cite{Cordes_Survey}.

\subsection*{Morse Boundaries}

Given a function $N \colon [1,\infty) \times [0,\infty) \to [0,\infty)$ we say that a geodesic $\gamma \colon I \to G$ is $N$-Morse if every $(k,c)$-quasi-geodesic $\alpha$ with endpoints on $\gamma$ is contained in the $N(k,c)$-neighborhood of $\gamma$.  We call the function $N$ the \emph{Morse gauge} of $\gamma$. Let $\mc{M}$ denote the set of Morse gauges of all geodesics in $G$. The set $\mc{M}$ has a partial order where $N \leq N'$ if $N(k,c) \leq N'(k,c)$ for all $k\geq 1$ and $c \geq 0$.

For each $N \in \mc{M}$, we let $G^N$ denote the set of elements $g\in G$ so that there exists a $N$-Morse geodesic connecting the identity $e$ to $g$. Cordes and Hume proved that the restriction of the metric of $G$ to $G^N$ makes $G^N$ into a  $\delta$-hyperbolic metric space, where $\delta$ is determined by $N$ \cite[Proposition 2.3]{CordesDurham2017}. Let $\partial_M^N G$ denote the Gromov boundary of $G^N$ with $e$ as the base point. By definition, if $N \leq N'$, then $G^N \subseteq G^{N'}$ and $\partial_M^N G \subseteq \partial_M^{N'} G$.

The \emph{Morse boundary} of $G$ is the direct limit of the $\partial_M^N G$. That is, $$\partial_M G = \varinjlim_\mathcal{M} \partial_M^N G. $$

A consequence of this definition is that each strata $\partial_M^N G$ is a compact subset of the Morse boundary.

\begin{lemappend}[{\cite[Theorem 3.14]{Cordes_Hume_boundary}}]\label{lem:strata_compact}
	For each Morse gauge $N$, $\partial_M^N G$ is a compact subset of $\partial_M G$.
\end{lemappend}

An important property of the Morse boundary is that the Gromov boundaries of stable subgroup embed into the Morse boundary of the entire group; see \cite{CordesDurham2017} or \cite{DT_stability} for the definition of a stable subgroup.

\begin{lemappend}[{\cite[Corollary 2.12]{CordesDurham2017}}]\label{lem:stable_embedding}
	 If $H$ is a stable subgroup of $G$, then $H$ is hyperbolic and the inclusion of $H$ into $G$ induces a continuous injection of $\partial H$ into $\partial_M G$.
\end{lemappend}

Charney, Cordes, and Sisto defined an $\omega$-Cantor space as the direct limit of countably many Cantor spaces where each Cantor space has empty interior in the next Cantor space; see \cite{CCS_omega_Morse_boundary} for full details. They prove that all $\omega$-Cantor sets are homeomorphic and give the following sufficient condition for the Morse boundary of a group to be an $\omega$-Cantor set. 

\begin{thmappend}[{\cite[Theorem 1.4]{CCS_omega_Morse_boundary}}]\label{thm:sufficent_omega_Morse}
	If $\partial_M G$ is totally disconnected, $\sigma$-compact, and contains a Cantor space, then $\partial_M G$ is homeomorphic to either a Cantor space or an $\omega$-Cantor space. Moreover,  $\partial_M G$ is a Cantor space if and only if $G$ is virtually free.
\end{thmappend}

\subsection*{Morse detectability} 
Our proof that certain HHGs have $\omega$-Cantor space boundaries only needs the following \emph{Morse detectable} property.

\begin{defnappend}[Morse Detectable]\label{defn:Morse_detectable}
	A finitely generated group $G$ is \emph{Morse detectable} if there exists a hyperbolic space $X$ and a coarsely Lipschitz map $\pi \colon G \to X$ so that  for any $(k,c)$-quasi-geodesic $\gamma \colon I \to G$ the following hold.
	\begin{enumerate}
		\item  For each $N \in \mc{M}$, if $\gamma$ is $N$-Morse, then there exists $\lambda\geq 1$ depending on $N$, $k$, and $c$, so that $\pi \circ \gamma$ is a $(\lambda,\lambda)$-quasi-geodesic of $X$.
		\item For each $\lambda \geq 1$, if $\pi \circ \gamma$ is an $(\lambda,\lambda)$-quasi-geodesic, then  $\gamma$ is $N$-Morse for some $N$ determined by $k$, $c$, and $\lambda$.
	\end{enumerate}
	If the space $X$ is a quasi-tree, then we say $G$ is \emph{Morse detectable in a quasi-tree}.
\end{defnappend}

Abbott, Behrstock, and Durham prove that unbounded products imply a hierarchically hyperbolic groups is also Morse detectable.

\begin{thmappend}[{\cite[Corollary 6.2]{ABD}}]\label{thm:ABD_implies_Morse_detectable}
	If $(G,\mf{S})$ is a hierarchically hyperbolic group with unbounded products, then $G$ is Morse detectable in the $\nest$-maximal hyperbolic space of $\mf{S}$.
\end{thmappend}

Morse detectability allows us to conclude that  the Morse boundary is $\sigma$-compact and continuously embeds in the Gromov boundary of the detecting space.
\begin{lemappend}\label{lem:boundary_properties}
	If $G$ is Morse detectable in the space $X$, then
	\begin{enumerate}
		\item  $\partial_M G$ is $\sigma$-compact;
		\item the map $\pi \colon G \to X$ induces a continuous injection $\partial \pi \colon \partial_M G \to \partial X$.
	\end{enumerate}	
\end{lemappend}

\begin{proof}
	\textbf{$\sigma$-compact.} For each $\lambda \in \mathbb{N}$, define $\partial_M^\lambda G$ to be the  subset of points in $\partial_M G$ that are represented by a geodesic ray $\gamma \colon I \to G$ so that $\pi \circ \gamma$ is a $(\lambda,\lambda)$-quasi-geodesic ray in $X$. We claim that the closures of the $\partial_M^\lambda G$ give an exhaustion of $\partial_M G$ by compact sets.
	
	Because $G$ is Morse detectable, for each $\lambda \in \mathbb{N}$,  there exists $N \in \mc{M}$ so that  $\partial_M^\lambda G$ is contained in some $\partial_M^N G$ for $N$ determined by $\lambda$. Thus the closure of each $\partial_M^\lambda G$ is compact since  each $\partial_M^N G$ is compact (Lemma \ref{lem:strata_compact}). On the other hand, Morse detectability also implies that for each Morse gauge $N \in \mc{M}$, there is some $\lambda \in \mathbb{N}$ so that $\partial_M^N G \subseteq \partial_M^\lambda G$. Thus, $\partial_M G$ is contained in the union of the closures of the $\partial_M^\lambda G$.
	
	\textbf{Continuous embedding.} The definition of Morse detectable make the restriction of $\pi$ to $G^N$ a quasi-isometric embedding with constants depending ultimately only on $N$. Since quasi-isometric embeddings of hyperbolic spaces extend to continuous injections on their boundaries, we have that $\pi$ induces a continuous injection $\partial_M^N G \to \partial X$ for each $N \in \mc{M}$. Since $\partial_M G$ is the direct limit of the $\partial_M^N G$'s, this induces a continuous injection $\partial \pi \colon \partial_M G \to \partial X$.
\end{proof}

When our group is Morse detectable in a quasi-tree, we can verify Charney, Cordes, and Sisto's condition for the boundary to be an $\omega$-Cantor set. 

\begin{thmappend}\label{thm:morse_detectable_in_tree}
	Suppose $\partial_M G$ is non-empty. If $G$ is Morse detectable in a quasi-tree and not hyperbolic, then $\partial_M G$ is homeomorphic to an $\omega$-Cantor space.
\end{thmappend}	

\begin{proof}
	Let $X$ be the quasi-tree that detects the Morse quasi-geodesics of $G$. By Lemma \ref{lem:boundary_properties}, $\partial_M G$ is $\sigma$-compact. Since the Gromov boundary of a quasi-tree is totally disconnected, the continuous injection $\partial \pi \colon \partial_M G \to \partial X$ forces $\partial_M G$ to be totally disconnected. Therefore, the only requirement of Theorem \ref{thm:sufficent_omega_Morse} that we still need to verify is that $\partial_M G$ contains a Cantor space.
	
	Corollary 4.9 of \cite{Morse-local-to-global} says any Morse detectable group that has non-empty Morse boundary is either virtually cyclic or contains a stable free subgroup of rank at least 2. The desired Cantor subspace thus comes from the Gromov boundary of this stable free subgroup  by Lemma \ref{lem:stable_embedding}.
\end{proof}

As a corollary of Theorem \ref{thm:morse_detectable_in_tree}, we show that the Morse boundaries of some HHGs are $\omega$-Cantor spaces. This proves Theorem \ref{thm:Morse_boundary} from the main text.

\begin{corappend}\label{cor:HHG_with_omega_boundary}
	Let $(G,\mf{S})$ be a non-hyperbolic hierarchically hyperbolic group with non-empty Morse boundary. If $\mf{S}$ has unbounded products and the $\nest$-maximal hyperbolic space of $\mf{S}$ is a quasi-tree, then $\partial_M G$ is homeomorphic to an $\omega$-Cantor space.
\end{corappend}

\begin{proof}
	By Theorem \ref{thm:ABD_implies_Morse_detectable}, the group $G$ is Morse detectable with the $\nest$-maximal hyperbolic space in $\mf{S}$. Thus, the conclusion is a special case of Theorem \ref{thm:morse_detectable_in_tree}.
\end{proof}

We conclude by highlighting  another hierarchically hyperbolic group whose Morse boundary is an $\omega$-Cantor space because  of Corollary \ref{cor:HHG_with_omega_boundary}.

\begin{exappend}
	Let $H$ be the handlebody group of a genus $2$ handlebody. Che sser proved $H$ has a hierarchically hyperbolic group structure with unbounded products \cite{Chesser_genus_2_handlebody}. Further, Chesser shows that the $\nest$-maximal space in this structure is a quasi-tree because it is the contact graph of a CAT(0) cube complex (they also prove this space is quasi-isometric to the disk graph of the genus 2 handlebody). Thus, $\partial_M H$ is an $\omega$-Cantor space. In the case of genus 3 or larger, the handlebody group is known to not be hierarchically hyperbolic \cite{Hamenstadt_Hensel_Handlebody} although the question of whether it is or is not Morse detectable is open.
\end{exappend}

\end{appendices}

\bibliography{Bib_Extensions}{}

\begin{thebibliography}{BKMM12}

\bibitem[ABD21]{ABD}
Carolyn Abbott, Jason Behrstock, and Matthew~G. Durham.
\newblock Largest acylindrical actions and stability in hierarchically
  hyperbolic groups. {W}ith an appendix by {D}aniel {B}erlyne and {J}acob
  {R}ussell.
\newblock {\em Trans. Amer. Math. Soc. Ser. B}, 8:66--104, 2021.

\bibitem[ABO19]{ABO_largest_action}
Carolyn Abbott, Sahana~H. Balasubramanya, and Denis Osin.
\newblock Hyperbolic structures on groups.
\newblock {\em Algebr. Geom. Topol.}, 19(4):1747--1835, 2019.

\bibitem[Aou13]{Aougab_Uniform_hyperbolicity}
Tarik Aougab.
\newblock Uniform hyperbolicity of the graphs of curves.
\newblock {\em Geom. Topol.}, 17(5):2855--2875, 2013.

\bibitem[BBKL20]{BBKL_undistorted_pA}
Mladen Bestvina, Kenneth Bromberg, Autumn~E. Kent, and Christopher~J.
  Leininger.
\newblock Undistorted purely pseudo-{A}nosov groups.
\newblock {\em J. Reine Angew. Math.}, 760:213--227, 2020.

\bibitem[BDM09]{behrstock_Drutu_Mosher_Thickness}
Jason Behrstock, Cornelia Dru\c{t}u, and Lee Mosher.
\newblock Thick metric spaces, relative hyperbolicity, and quasi-isometric
  rigidity.
\newblock {\em Math. Ann.}, 344(3):543--595, 2009.

\bibitem[Beh06]{Behrstock_Thesis}
Jason~A. Behrstock.
\newblock Asymptotic geometry of the mapping class group and {T}eichm\"uller
  space.
\newblock {\em Geom. Topol.}, 10:1523--1578, 2006.

\bibitem[BH99]{BH}
Martin~R. Bridson and Andr\'e Haefliger.
\newblock {\em Metric spaces of non-positive curvature}, volume 319 of {\em
  Grundlehren der Mathematischen Wissenschaften [Fundamental Principles of
  Mathematical Sciences]}.
\newblock Springer-Verlag, Berlin, 1999.

\bibitem[BHMS20]{BHMS_cHHS}
Jason Behrstock, Mark Hagen, Alexandre Martin, and Alessandro Sisto.
\newblock A combinatorial take on hierarchical hyperbolicity and applications
  to quotients of mapping class groups.
\newblock arXiv:2005.00567, 2020.

\bibitem[BHS17a]{BHS_HHS_AsDim}
Jason Behrstock, Mark~F. Hagen, and Alessandro Sisto.
\newblock Asymptotic dimension and small-cancellation for hierarchically
  hyperbolic spaces and groups.
\newblock {\em Proc. Lond. Math. Soc. (3)}, 114(5):890--926, 2017.

\bibitem[BHS17b]{BHS_HHSI}
Jason Behrstock, Mark~F. Hagen, and Alessandro Sisto.
\newblock Hierarchically hyperbolic spaces {I}: curve complexes for cubical
  groups.
\newblock {\em Geom. Topol.}, 21(3):1731--1804, 2017.

\bibitem[BHS19]{BHS_HHSII}
Jason Behrstock, Mark~F. Hagen, and Alessandro Sisto.
\newblock {Hierarchically hyperbolic spaces II: combination theorems and the
  distance formula}.
\newblock {\em Pacific J. Math.}, 299:257--338, 2019.

\bibitem[BHS21]{BHS_HHS_Quasiflats}
Jason Behrstock, Mark~F. Hagen, and Alessandro Sisto.
\newblock Quasiflats in hierarchically hyperbolic spaces.
\newblock {\em Duke Math. J.}, 170:909--996, 2021.

\bibitem[BKMM12]{BKMM_Rigidity}
Jason Behrstock, Bruce Kleiner, Yair Minsky, and Lee Mosher.
\newblock Geometry and rigidity of mapping class groups.
\newblock {\em Geom. Topol.}, 16(2):781--888, 2012.

\bibitem[BM08]{BM_MCG_Rank}
Jason~A. Behrstock and Yair~N. Minsky.
\newblock Dimension and rank for mapping class groups.
\newblock {\em Ann. of Math. (2)}, 167(3):1055--1077, 2008.

\bibitem[Bow06]{Bowditch_Guessing_Geodesic}
Brian~H. Bowditch.
\newblock Intersection numbers and the hyperbolicity of the curve complex.
\newblock {\em J. Reine Angew. Math.}, 598:105--129, 2006.

\bibitem[Bow14]{bowditch_uniform_hyperbolicity}
Brian~H. Bowditch.
\newblock Uniform hyperbolicity of the curve graphs.
\newblock {\em Pacific J. Math.}, 269(2):269--280, 2014.

\bibitem[BR20]{BR_Combination}
Federico Berlai and Bruno Robbio.
\newblock A refined combination theorem for hierarchically hyperbolic groups.
\newblock {\em Groups Geom. Dyn.}, 14(4):1127--1203, 2020.

\bibitem[Bro03]{Brock_WeilPetersson}
Jeffrey~F. Brock.
\newblock The {W}eil-{P}etersson metric and volumes of 3-dimensional hyperbolic
  convex cores.
\newblock {\em J. Amer. Math. Soc.}, 16(3):495--535, 2003.

\bibitem[CCS23]{CCS_omega_Morse_boundary}
Ruth Charney, Matthew Cordes, and Alessandro Sisto.
\newblock Complete topological descriptions of certain {M}orse boundaries.
\newblock {\em Groups Geom. Dyn.}, 17(1):157--184, 2023.

\bibitem[CD17]{CordesDurham2017}
Matthew Cordes and Matthew~G. Durham.
\newblock Boundary convex cocompactness and stability of subgroups of finitely
  generated groups.
\newblock {\em International Mathematics Research Notices}, 2019(6):1699--1724,
  2017.

\bibitem[CH17]{Cordes_Hume_boundary}
Matthew Cordes and David Hume.
\newblock Stability and the {M}orse boundary.
\newblock {\em J. Lond. Math. Soc. (2)}, 95(3):963--988, 2017.

\bibitem[Che22]{Chesser_genus_2_handlebody}
Marissa Chesser.
\newblock Stable subgroups of the genus 2 handlebody group.
\newblock {\em Algebr. Geom. Topol.}, 22(2):919--971, 2022.

\bibitem[Cor19]{Cordes_Survey}
Matthew Cordes.
\newblock A survey on {M}orse boundaries and stability.
\newblock In {\em Beyond hyperbolicity}, volume 454 of {\em London Math. Soc.
  Lecture Note Ser.}, pages 83--116. Cambridge Univ. Press, Cambridge, 2019.

\bibitem[CRS14]{CRS_Uniform_hyperbolicity}
Matt Clay, Kasra Rafi, and Saul Schleimer.
\newblock Uniform hyperbolicity of the curve graph via surgery sequences.
\newblock {\em Algebr. Geom. Topol.}, 14(6):3325--3344, 2014.

\bibitem[DDLS]{DDLS_Veech_groups_II}
Spencer Dowdal, Matthew~G. Durham, Christopher~J. Leininger, and Alessandro
  Sisto.
\newblock {Extensions of Veech groups II: Hierarchical hyperbolicity and
  quasi-isometric rigidity}.
\newblock arXiv:2111.00685.

\bibitem[DT15]{DT_stability}
Matthew~G. Durham and Samuel~J. Taylor.
\newblock Convex cocompactness and stability in mapping class groups.
\newblock {\em Algebr. Geom. Topol.}, 15(5):2839--2859, 2015.

\bibitem[Dur16]{Durham_Combinatorial_Teich}
Matthew~G. Durham.
\newblock The augmented marking complex of a surface.
\newblock {\em Journal of the London Mathematical Society}, 94(3):933, 2016.

\bibitem[EMR17]{EMR_Teich_Rank}
Alex Eskin, Howard~A. Masur, and Kasra Rafi.
\newblock Large-scale rank of {T}eichm\"uller space.
\newblock {\em Duke Math. J.}, 166(8):1517--1572, 2017.

\bibitem[FM02]{Farb_Mosher_convex_cocompact}
Benson Farb and Lee Mosher.
\newblock Convex cocompact subgroups of mapping class groups.
\newblock {\em Geom. Topol.}, 6:91--152, 2002.

\bibitem[Ham]{Hamenstadt_extensions_of_surface_groups}
Ursula Hamenstadt.
\newblock Word hyperbolic extensions of surface groups.
\newblock arXiv:math/0505244, 2005.

\bibitem[HH21]{Hamenstadt_Hensel_Handlebody}
Ursula Hamenst\"{a}dt and Sebastian Hensel.
\newblock The geometry of the handlebody groups {II}: {D}ehn functions.
\newblock {\em Michigan Math. J.}, 70(1):--, 2021.

\bibitem[HHP20]{HHP_Coarse_Helly_and_HHG}
Thomas Haettel, Nima Hoda, and Harry Petyt.
\newblock The coarse helly property, hierarchical hyperbolicity, and
  semihyperbolicity.
\newblock arXiv:2009.14053, 2020.

\bibitem[HMS24]{Extra_Large_Artin}
Mark Hagen, Alexandre Martin, and Alessandro Sisto.
\newblock Extra-large type {A}rtin groups are hierarchically hyperbolic.
\newblock {\em Math. Ann.}, 388(1):867--938, 2024.

\bibitem[HPW15]{HPW_Uniform_Hyperbolicity}
Sebastian Hensel, Piotr Przytycki, and Richard C.~H. Webb.
\newblock 1-slim triangles and uniform hyperbolicity for arc graphs and curve
  graphs.
\newblock {\em J. Eur. Math. Soc. (JEMS)}, 17(4):755--762, 2015.

\bibitem[HS20]{Hagen_Susse_factor_system}
Mark~F. Hagen and Tim Susse.
\newblock On hierarchical hyperbolicity of cubical groups.
\newblock {\em Israel J. Math.}, 236(1):45--89, 2020.

\bibitem[KL07]{Kent_Leininger_Survey}
Autumn Kent and Christopher~J. Leininger.
\newblock Subgroups of mapping class groups from the geometrical viewpoint.
\newblock In {\em In the tradition of {A}hlfors-{B}ers. {IV}}, volume 432 of
  {\em Contemp. Math.}, pages 119--141. Amer. Math. Soc., Providence, RI, 2007.

\bibitem[KL08a]{Kent_Leininger_convex_cocompactness}
Autumn~E. Kent and Christopher~J. Leininger.
\newblock Shadows of mapping class groups: capturing convex cocompactness.
\newblock {\em Geom. Funct. Anal.}, 18(4):1270--1325, 2008.

\bibitem[KL08b]{Kent_Leininger_convergence}
Autumn~E. Kent and Christopher~J. Leininger.
\newblock Uniform convergence in the mapping class group.
\newblock {\em Ergodic Theory Dynam. Systems}, 28(4):1177--1195, 2008.

\bibitem[KLS09]{KLS_Trees_and_MCG}
Autumn~E. Kent, Christopher~J. Leininger, and Saul Schleimer.
\newblock Trees and mapping class groups.
\newblock {\em J. Reine Angew. Math.}, 637:1--21, 2009.

\bibitem[KR14]{Kapovich_Rafi_Hyp_Free_Factors}
Ilya Kapovich and Kasra Rafi.
\newblock On hyperbolicity of free splitting and free factor complexes.
\newblock {\em Groups Geom. Dyn.}, 8(2):391--414, 2014.

\bibitem[Man05]{ManningBottleneck}
Jason~F. Manning.
\newblock Geometry of pseudocharacters.
\newblock {\em Geom. Topol.}, 9(2):1147--1185, 2005.

\bibitem[MM99]{MMI}
Howard~A. Masur and Yair~N. Minsky.
\newblock Geometry of the complex of curves. {I}. {H}yperbolicity.
\newblock {\em Invent. Math.}, 138(1):103--149, 1999.

\bibitem[MM00]{MMII}
Howard~A. Masur and Yair~N. Minsky.
\newblock Geometry of the complex of curves. {II}. {H}ierarchical structure.
\newblock {\em Geom. Funct. Anal.}, 10(4):902--974, 2000.

\bibitem[MS12]{Mj_Pranab_Metric_bundles}
Mahan Mj and Pranab Sardar.
\newblock A combination theorem for metric bundles.
\newblock {\em Geom. Funct. Anal.}, 22(6):1636--1707, 2012.

\bibitem[Ota96]{Thurston_fibering_over_circle}
Jean-Pierre Otal.
\newblock Le th\'{e}or\`eme d'hyperbolisation pour les vari\'{e}t\'{e}s
  fibr\'{e}es de dimension 3.
\newblock {\em Ast\'{e}risque}, (235):x+159, 1996.

\bibitem[Raf07]{Rafi_Combinatorial_Teich}
Kasra Rafi.
\newblock A combinatorial model for the {T}eichm\"uller metric.
\newblock {\em Geom. Funct. Anal.}, 17(3):936--959, 2007.

\bibitem[RS]{Robbio_Spriano_Combination}
Bruno Robbio and Davide Spriano.
\newblock Hierarchical hyperbolicity of hyperbolic-2-decomposable groups.
\newblock arXiv:2007.13383.

\bibitem[RST22]{Morse-local-to-global}
Jacob Russell, Davide Spriano, and Hung~Cong Tran.
\newblock The local-to-global property for {M}orse quasi-geodesics.
\newblock {\em Math. Z.}, 300(2):1557--1602, 2022.

\bibitem[RST23]{RST_Quasiconvexity}
Jacob Russell, Davide Spriano, and Hung~C. Tran.
\newblock Convexity in hierarchically hyperbolic spaces.
\newblock {\em Algebr. Geom. Topol.}, 23(3):1167--1248, 2023.

\bibitem[Sis19]{HHS_survey}
Alessandro Sisto.
\newblock What is a hierarchically hyperbolic space?
\newblock In {\em Beyond hyperbolicity}, volume 454 of {\em London Math. Soc.
  Lecture Note Ser.}, pages 117--148. Cambridge Univ. Press, Cambridge, 2019.

\bibitem[Vok22]{vokes}
Kate~M. Vokes.
\newblock Hierarchical hyperbolicity of graphs of multicurves.
\newblock {\em Algebr. Geom. Topol.}, 22(1):113--151, 2022.

\end{thebibliography}
\bibliographystyle{alpha}	
\end{document}